\documentclass[preprint,11tpt]{elsarticle}

\evensidemargin \oddsidemargin
\advance \topmargin by -\headheight
\usepackage{lineno,hyperref}
\modulolinenumbers[5]

\usepackage{graphicx}
\usepackage{psfrag}

\usepackage{amsmath}
\usepackage{amsfonts}
\usepackage{amssymb}

\usepackage{rotating}
\usepackage[pdftex]{color}
\begin{document}

\begin{frontmatter}
\title{A numerical study of the spectral properties of Isogeometric collocation matrices for
acoustic wave problems}

\author{Elena Zampieri\fnref{footnote1}}
\address{Department of Mathematics,Universit\`a di Milano,Via Saldini 50,20133 Milano,Italy.
E-mail:{\tt elena.zampieri@mat.unimi.it}.}
\fntext[footnote1]{This work was partially supported by
the European Research Council through the FP7 Ideas Consolidator Grant
\emph{HIGEOM} n. 616563, by the Italian Ministry of Education, University and Research (MIUR) through the ``Dipartimenti di Eccellenza Program 2018-22 -
Dept. of Mathematics,University of Pavia'', and by
the Istituto Nazionale di Alta Matematica (INdAM - GNCS), Italy.}

\author{Luca F. Pavarino\fnref{footnote1}}
\address{Department of Mathematics,Universit\`a di Pavia,Via Ferrata 5,27100 Pavia,Italy.
E-mail:{\tt luca.pavarinoi@unipv.it}.}

\begin{abstract}
This paper focuses on the spectral  properties of the mass and stiffness matrices for acoustic wave problems discretized with Isogeometric analysis (IGA) collocation methods in space and Newmark methods in time. Extensive numerical results are reported for the eigenvalues and condition numbers of the acoustic mass and stiffness matrices in the reference square domain with Dirichlet, Neumann and absorbing boundary conditions, varying the polynomial degree $p$, mesh size $h$, regularity $k$, of the IGA discretization and the time step $\Delta t$ and parameter $\beta$ of the Newmark method. Results on the sparsity of the matrices and the eigenvalue distribution with respect to the degrees of freedom {\textsf { d.o.f.}} and the number of    nonzero entries {\textsf  {nz}} are also reported.
The results are comparable with the available spectral estimates for IGA  Galerkin matrices associated to the Poisson problem with Dirichlet boundary conditions, and in some cases  the IGA collocation results are better than the corresponding IGA Galerkin estimates.  
\end{abstract}


\end{frontmatter}

%\linenumbers

\section{Introduction}
Since its introduction in \cite{HuCoBa05}, Isogeometric analysis (IGA) has generated a large amount of work and interest in various fields involving the numerical solution of partial differential equations. Several studies have shown the advantages of the IGA approach in many problems and applications, see, e.g., \cite{BaBeCoHuSa06,ABHRS,CoHuBa:book,BeBuSanVa14},  and the references therein. IGA employs basis functions associated with $B$-splines and non-uniform rational $B$-splines (NURBS) to   discretize both the problem domain, as in computer aided design (CAD) systems, and the solution space of the differential problem. This approach allows a simple representation of the problem geometry and at the same time yields high order methods with respect to standard $p$- and $hp$- refinements, where $p$ is the polynomial degree of the IGA basis  functions and $h$ is the mesh size. An additional $k$-refinement, where $k \leq p-1$ is the global regularity of the IGA basis functions, provides highly continuous numerical solutions and better accuracy than in the case of classical FEM $p$-refinement.

The application of NURBS-based IGA to acoustic  and elastic wave propagation problems has initially being carried out using standard Galerkin approaches.
More recently, IGA collocation variants have been investigated, with the aim of dealing with sparser mass and stiffness  matrices than those arising from IGA Galerkin techniques. IGA collocation has also the additional advantage of reducing the global computational cost, 
%because of the easier and less expensive construction of the 
since collocation matrices require only one function evaluation per collocation point, independently of $p$; see \cite{ABHRS2,CRBH2006,djq2015,ehhr2018,HRS09,ipfs2014,kntdlh15,zdq2017}.

In our previous work \cite{pzIGA19}, we have considered IGA Galerkin and explicit Newmark approximations of the acoustic wave equation with absorbing boundary conditions, whereas in \cite{pzIGA21} we have extended the study to IGA collocation and implicit Newmark schemes. Several numerical results have shown the convergence and stability properties of these schemes, and we have presented a detailed comparison between IGA Galerkin and IGA Collocation with respect to the space and time discretization parameters, CPU time, sparsity of matrices  and degrees of freedom. 
Since  both the IGA Galerkin and  collocation mass matrices becomes fuller for increasing $p$ and $k$,  the main difference between explicit and implicit IGA Newmark schemes is related to the stability bounds for the time step, rather than to the cost of the solution of the linear systems at each temporal instant. The theoretical analysis in \cite{pzIGA19} is confined to the stability properties of the IGA Galerkin Newmark method and additionally  it is only partially based on proven results. In fact,  there is still a lack of theoretical spectral bounds for IGA matrices in the literature, and most of the known estimates regarding eigenvalues and conditioning of the mass and stiffness matrices are conjectures.  
For the above reasons, it follows that in the framework of propagation problems a detailed experimental analysis is of interest both to fill the gaps of the theoretical analysis and for the investigation of efficient solution of the linear system at each time step of the time advancing scheme, possibly involving preconditioning techniques.
Among other relevant works,
\cite{iga-sem2018} presents a methodical numerical comparison between Spectral Element Method and NURBS-based IGA Galerkin when applied to the Poisson problem, analyzing the convergence, computational costs,  and conditioning with respect to $h$ and $p$, for minimal ($k=0$) and maximal ($k=p-1$) regularity of the IGA basis functions. 
In \cite{loli2022}, the authors study the  condition number of IGA Galerkin mass matrices and propose efficient preconditioners for the related linear systems, focusing in particular on $k$-refinement. In \cite{ehhr2018}, the  spectral properties of a semi-discrete predictor–multicorrector method are investigated for the   case of a 1D pure Dirichlet IGA problem. 

In this paper, we consider the approximation of acoustic wave problems with absorbing boundary conditions  based on IGA collocation at Greville points in space and Newmark advancing schemes in time, both explicit and implicit. We observe  that the implementation of absorbing boundary conditions is mathematically equivalent to the more common Robin boundary conditions, as it involves a linear combination of the values of the function and of its normal derivative at the collocation points on the domain boundary. 
Differently from our previous works \cite{pzIGA19} and  \cite{pzIGA21}, where we focused on 
%theoretical and numerical results related to  
the convergence and stability properties of IGA Galerkin and collocation methods, 
%with respect to the discretization parameters $p$, $h$, $k$, and $\Delta t$. 
in this paper we focus instead on the spectral  properties of the mass and stiffness matrices for acoustic wave problems discretized with IGA collocation in space and Newmark methods in time. We present a numerical study of the behaviour of the eigenvalues and condition numbers of the mass and stiffness matrices in the reference square domain with Dirichlet, Neumann and absorbing boundary conditions, varying the polynomial degree $p$, mesh size $h$, regularity $k$, time step $\Delta t$ and parameter $\beta$ of the Newmark method. We also   report some results on the sparsity of the matrices and the eigenvalue distribution with respect to the degrees of freedom {\textsf { d.o.f.}} and the number of    nonzero entries {\textsf  {nz}}.
In order to provide a simple basis of comparison, we recall some bounds and estimates that are available for IGA  Galerkin matrices associated to the Poisson problem with Dirichlet boundary conditions, see \cite{gahalaut2012},  \cite{garoni2014} and \cite{iga-sem2018}. Our results show that the same estimates hold for the condition numbers of the  IGA collocation matrices considered in this paper, and in some cases  the collocation matrices satisfy even better bounds.  

The rest of the paper is organized as follows. The acoustic wave model problem and its mathematical analysis are introduced in Section 2, and its approximation by IGA collocation in space and Newmark methods in time in Section 3. In Section 4, we give a brief overview of eigenvalue and condition number estimates  for the  IGA Galerkin approximation of the Poisson problem. Finally, in Section 5, we present several numerical tests on the behaviour of eigenvalues and condition numbers of IGA collocation mass and stiffness matrices with different types of boundary conditions, varying all the discretization parameters, and comparing the results with the ones reported for the IGA Galerkin case.

\section{The  model problem and mathematical analysis}
\label{acoustic_waves}

We consider the reference square $\Omega=[0,1]\times [0,1]$ in the plane with boundary $\partial \Omega$, while $(0,T)$ is the temporal interval, with $T$ real and positive. Any point of $\Omega$ is denoted by $\mathbf{x}$ and  $t$  represents the time. The acoustic wave problem (see e.g., Junger and Feit \cite{jf86} and Ihlenburg \cite{ihlenburg}) is as follows:

\begin{equation}
\frac {\partial ^2 u}{\partial t^2}({\bf x},t) -c_0\Delta  u({\bf x},t) =f({\bf x},t)  \quad {\rm
in}\ \Omega
\times (0,T),
\label{pde}
\end{equation}
\noindent with    initial conditions
\begin{equation}
u({\bf x},0)=u_0({\bf x}) , \qquad \frac{\partial u}{\partial t} ({\bf x},0)= u_1({\bf x})
\qquad   {\rm in}\  \Omega ,
\label{init_cond}
\end{equation}
\noindent where $u$ is the unknown acoustic pressure, $c_0$ is the acoustic wave propagation velocity, $f$ is the source term, $u_0$ and $u_1$ are the initial pressure and velocity, respectively.

\noindent Standard  Dirichlet or Neumann boundary conditions can be prescribed on $\partial \Omega$:
\begin{equation}
u({\bf x},t)=\Phi ({\bf x},t) \quad  {\rm on} \ \Gamma_D \times (0,T), \qquad \frac
{\partial u}{\partial {\bf n}}({\bf x},t) =\Psi  ({\bf x},t)
\quad {\rm on}\  \Gamma_N   \times (0,T).
\label{bound_cond}
\end{equation}

\noindent where $\Phi$ and $\Psi$ are the prescribed pressure and velocity on $\Gamma_D$ and $\Gamma_N$, respectively, and $\mathbf{n}$ is the outward boundary normal unit vector. We recall that 
when $\Psi=0$ the homogeneous Neumann condition represents a free surface in correspondence of which full reflection occurs, whereas the case $\Phi\neq 0$ arises when the medium is subjected to a source located at the portion of the boundary $\Gamma_N$. Nevertheless, Dirichlet conditions are less usual when dealing with acoustic wave problems, since the physical solution is rarely assigned at some portion of the boundary. Other common boundary conditions in the simulation of wave propagation through an unbounded domain are the so called absorbing boundary conditions (ABCs for brevity), involving a truncation of the infinite original domain. On artificial boundaries of the novel finite domain suitable boundary conditions are then imposed in order to eliminate or reduce as much as possible spurious wave reflections.
Since exact transmitting boundary conditions are non-local neither in space nor in time, and consequently not useful for numerical implementation, several ABCs have been introduced in the literature with the aim to make the boundary transparent to outgoing and opaque to ingoing waves (see e.g., Clayton and Engquist \cite{ce77} and Engquist and Majda \cite{em79}). Here we consider natural first-order ABCs based on first spatial and temporal partial derivatives only \cite{mur81}:
\begin{equation}
\frac{1}{\sqrt{c_0}}\frac{\partial u}{\partial t}  ({\bf x},t)+ \frac{\partial u}{\partial {\bf n}} ({\bf x},t) =0
\qquad {\rm on}\ \Gamma_{AB}\times (0,T) ,
\label{ABC} 
\end{equation}
where $\Gamma_{AB}$ is the artificial boundary where ABCs are enforced and, in the most general case,  $\partial \Omega=\Gamma_D\cup \Gamma_N \cup \Gamma_{AB}$.

%To introduce the weak formulation of (\ref{pde})-(\ref{ABC}), we assume %that for each $t\in (0,T)$, $f\in L^2(\Omega\times (0,T)),
%\  \Phi \in H^{\frac {1}{2}}(\partial \Omega\times (0,T))$,
%$\Psi \in  L^{2}(\partial \Omega\times (0,T))$,
%$u_0 \in H^1(\Omega)$ and $u_1 \in L^2(\Omega)$, and introduce the %symmetric bilinear form and inner products:
%\[
%a(u,v)=c_0\int_{\Omega}\nabla u \cdot \nabla v \ {\rm d}{\bf x},\qquad
%(f,v)= \int_{\Omega} fv\  {\rm d}{\bf x},
%  \qquad
% <\Psi,v>_{\Gamma_{(\cdot)}}=\int_{\Gamma_{(\cdot)}} \Psi v \ {\rm d} s,
%\] where $\Gamma_{(\cdot)}={\Gamma_N}$ or %$\Gamma_{(\cdot)}={\Gamma_{AB}}$.
%Then, the weak formulation reads as follows: Find $u:(0,T) \to
%H^1(\Omega)$, such that for $ a.e. \ t \in (0,T), \ u(t)=\Phi(t) \ {\rm %on} \ \Gamma_D$ and
%\begin{equation}
%\left(\frac {\partial ^2 u} {\partial t^2}, v\right)+a(u,v)
%+{\sqrt{c_0}}<\frac {\partial  u} {\partial t},v>_{\Gamma_{AB}}=
%(f,v)+<\Psi,v>_{\Gamma_N} 
%\quad \forall v \in V, \label{var_pde_ABC} \end{equation}
%\[
%<\frac {\partial  u} {\partial t},v>_{\Gamma_{(\cdot)}}=\int_{\Gamma_{(\cdot)}}
% \frac {\partial  u} {\partial t} v \ {\rm d} s, \]
%with $V= \{v \in H^1(\Omega) : v=0 \ {\rm on } \ \Gamma_D \}$.  We observe %that the partial derivative $\partial
%\cdot  / \partial t$ in (\ref{ABC}) can be inserted in (\ref{var_pde_ABC}) %as a boundary load on the boundary $\Gamma_{AB}$ associated to a Neumann %condition.  

\noindent The uniqueness of the solution and the stability of the continuous acoustic problem (\ref{pde})-(\ref{bound_cond})   can be proved
following the analysis of elasto-dynamics linear problems (see
Quarteroni, Tagliani and Zampieri \cite{qtz98}).

\noindent Higher-order ABCs have been introduced and  analyzed in the literature (see, e.g., Givoli \cite{g91}) and involve also derivatives of order greater than one in space and time, and derivatives in the tangential direction.

%---------------------------------------------------------------------------------------------------------------------
\section{Approximation of the wave problem by Isogeometric Collocation and Newmark methods}
\label{sec:approximation}

We introduce the numerical approximation of the strong form of the acoustic wave problem (\ref{pde})-(\ref{ABC}). The discretization of the space variable ${\mathbf{x}}=(x_1,x_2)$ is based on  the collocation form of Isogeometric approximation methods, whereas the time discretization is based on Newmark time advancing schemes.
%---------------------------------------------------------------------------------------------------------------------

%---------------------------------------------------------------------------------------------------------------------
\subsection{Isogeometric Analysis and Collocation Methods}\label{subsec_iga}

In this Section we briefly recall the basic notions of Non-Uniform Rational B-splines (NURBS) basis functions and introduce collocation points at which the strong form of the acoustic problem will be enforced.  Given a knot vector of non-decreasing real numbers on the reference interval 
\begin{equation}
\{\xi_1=0,...,\xi_{\nu+p+1}=1\},
\label{knot_vector}\end{equation}
we associate univariate B-spline basis functions $N^p_i$, where the integers $p$ and $\nu$ are the polynomial degree of the B-spline, and the number of basis functions and control points, respectively. 
Starting from a knot vector (\ref{knot_vector}), B-splines are built recursively starting from piecewise constant function when $p=0$, obtaining B-splines with support $(\xi_i,\xi_{i+p+1})$, $i=1,2,...,\nu$ (see, e.g. \cite{schumi})
It is known that B-spline basis functions are $C^{p-1}$-continuous if internal nodes are not repeated, whereas they are $C^k$-continuous, $k=p-\alpha$, $\alpha$ being the multiplicity of the associated knot. Moreover, when a knot has multiplicity $\alpha=p$, the basis is $C^0$-continuous, interpolating the control point at that location where the knot has multiplicity $\alpha$. From now on,  we will assume that the maximum knot
multiplicity is $p$ ensuring at least the global continuity of all considered functions. 
%As an example, in Fig. \ref{Greville_fig} we represent   basis functions associated with the open knot %vector $\boldsymbol{\xi}= \{0,0,0,0,1/6,1/3,1/2,2/3,5/6,1,1,1,1\}$, in the particular case $p=3$ and %$\nu=9$. 

%\begin{figure} %[!t]
%\centerline{\includegraphics[scale=0.5]{cub_basis_fcn.pdf}}
%\caption{Example of cubic basis functions $N^3_i$ associated to the knot vector %$\boldsymbol{\xi}= \{0,0,0,0,1/6,1/3,1/2,2/3,5/6,1,1,1,1\}$. 
%\label{basis_fig}}
%\end{figure}

\noindent Starting from the one-dimensional spline space 
\begin{equation}
  \label{eq:Sn}
  {\widehat{\mathcal{S}}_h} = \text{span}\{ N_{i}^{p}(\xi), i=1, \ldots, \nu\},
\end{equation}
we construct by tensor products multi-dimensional B-spline functions. For the  simplicity of exposition, we examine here the case of a two-dimensional domain, and B-spline of same degree $p$ in each directions. The case of higher-dimensional case and different degrees  can be dealt with analogously.
We introduce the two-dimensional parametric space $\widehat{\Omega} := (0,1)\times(0,1)$ with a   knot vector     (\ref{knot_vector}) in each direction, and   $\mathbf{C}_{i,j}$ is a net of  $\nu^2$  control points, $i,j=1,...,\nu$.  The bi-variate spline basis on $\widehat{\Omega}$ is then
$B_{i,j}^p(\xi,\eta)= N_i^p(\xi)N_j^p(\eta)$.
Similarly, the mesh of rectangular elements in the parametric space is generated in a natural way by the Cartesian product of two knot vectors $\{\xi_1=0,...,\xi_{\nu+p+1}=1\}$. Then,   
\begin{equation}
  \label{eq:Snm}
  {\widehat{\mathcal{S}}_h} = \text{span}\{ B_{i,j}^{p}(\xi,\eta) , \ i,j=1, \ldots, \nu  \}
\end{equation}  is the bi-variate spline space analogous to 
(\ref{eq:Sn}).
\noindent We recall that a rational B-spline in $\mathbb{R}^d$ is the
projection onto $d$-dimensional physical space of a polynomial
B-spline defined in ($d+1$)-dimensional homogeneous coordinate
space.  According to this, a great
variety of geometrical objects can be constructed. In
particular, all conic sections in physical space can be obtained
exactly. See
Rogers~\cite{Rogers01} and references therein for a complete discussion of these space projections.
We indicate a NURBS basis function of degree $p$ as
\begin{equation}
\label{eq:NURBS-1d-basis}
R_{i}^{p}(\xi) = \cfrac{N_{i}^{p}(\xi)\omega_i}{\sum_{\hat{i}=1}^\nu
N_{\hat{i}}^{p}(\xi)\omega_{\hat{i}}} = \cfrac{N_{i}^{p}(\xi)\omega_i}{w(\xi)} ,
\end{equation}
\noindent with $w(\xi)=\sum_{\hat{i}=1}^\nu
N_{\hat{i}}^{p}(\xi)\omega_{\hat{i}} \in {\widehat{\mathcal{S}}_h}$ a weight function. 
Analogously to the construction of B-splines, NURBS basis functions on the
two-dimensional parametric space $\widehat\Omega = (0,1)\times(0,1)$ are obtained from the bi-variate spline basis as
\begin{equation}
\label{eq:NURBS-2d-basis}
R_{i,j}^{p}(\xi,\eta) = \cfrac{B_{i,j}^{p}(\xi,\eta)\omega_{i,j}}
{\sum_{\hat{i},\hat{j}=1}^\nu 
B_{\hat{i},\hat{j}}^{p}(\xi,\eta)\omega_{\hat{i},\hat{j}}} = \cfrac{B_{i,j}^{p}(\xi,\eta)\omega_{i,j}}{w(\xi,\eta)},
\end{equation}
\noindent where $\omega_{i,j} \in \mathbb{R}$, and the denominator is the two-dimensional weight function commonly denoted   by $w(\xi,\eta)$. The continuity and support of NURBS basis functions are the same as for B-splines and NURBS spaces are the span of the basis functions (\ref{eq:NURBS-2d-basis}). 
Let us consider a \emph{single-patch} domain $\Omega$ as a NURBS region
associated with the net  
$\mathbf{C}_{i,j}$. We introduce the geometrical map
$\mathbf{F}:\widehat\Omega\rightarrow \Omega$
defined by
\begin{equation}
\label{eq:F}
\mathbf{F}(\xi,\eta)=\sum_{i,j=1}^{\nu} 
R^p_{i,j}(\xi,\eta)\mathbf{C}_{i,j}.
\end{equation}
\noindent According to the isogeometric approximation, the space of NURBS scalar
fields on the domain $\Omega$ is defined by the  isoparametric approach as the span of the  \emph{push-forward}  of the basis
functions (\ref{eq:NURBS-2d-basis})
\begin{equation}
 \label{eq:V-k-space}
 \mathcal{N}_{h} := \text{span}\{R_{i,j}^{p} \circ
 \mathbf{F}^{-1}, \text{ with } i,j=1,\ldots,
 \nu \}.
\end{equation}

\noindent Having defined the essential  elements of IGA spaces and basis functions, we can now review the IGA collocation method for the approximation of our acoustic wave problem in space, see \cite{ABHRS,ABHRS2,Schillinger}. As regards the   collocation points, several choices have been proposed in the literature, including Cauchy-Galerkin points \cite{gdl2016}, Demko abscissae \cite{demko-85}, Galerkin superconvergent points \cite{MST.2017} and Greville abscissae \cite{deboor-book}. In particular, we will employ the set of Greville collocation points in our numerical experiments. See \cite{pzIGA21} for stability and convergence properties.
Given the knot vector $\{\xi_1=0,...,\xi_{\nu+p+1}=1\}$, the corresponding Greville collocation points are 
\begin{equation}
  \label{greville-abscissae}
  \overline{\xi}_i \doteq (\xi_{i+1} + \xi_{i+2} + ... + \xi_{i+p})/p \ ,
\end{equation}
with   $\overline{\xi}_1=0$, $\overline{\xi}_{\nu}=1$, and the remaining points are in $(0,1)$. 
The tensor product 
$$
\widehat\tau_{ij} = (\overline{\xi}_i, \overline{\xi}_j) \ \in \ \overline{\big(\widehat\Omega\big)} \ ,
\ \ i,j=1,...,\nu,
$$provides the grid of collocation points $\tau_{ij} = {\bf F}(\widehat\tau_{ij})\in \Omega$.

The theoretical knowledge of the spectral properties and convergence estimates for the IGA collocation  of elliptic problems in two and three dimensions  is still an open issue. A large number of numerical tests are available in the literature regarding   convergence properties with respect to the main discretization parameters, namely the degree $p$, the mesh size $h$ and the regularity $k$ (e.g., \cite{ABHRS}, \cite{ABHRS2}, \cite{kntdlh15}, \cite{MST.2017}, \cite{Schillinger}, \cite{pzIGA21}). In Section \ref{Sec_Numerical_res} several numerical tests illustrate the spectral properties of the matrices arising from the IGA collocation methods, as well as their eigenvalue distribution in the complex plane.
 
%---------------------------------------------------------------------------------------------------------------
%---------------------------------------------------------------------------------------------------------------------
\subsection{Space discretization of the wave model problem (\ref{pde})-(\ref{ABC}) by IGA collocation methods}\label{subsec_IGA_Coll}
%---------------------------------------------------------------------------------------------------------------------
\noindent For a simpler description of the application of the  collocation problem to the wave problem, we  enumerate
the grid points $\{\tau_{ij}\in \Omega,\ i,j=1,...,p\}$ using only one index. Thus, each collocation point $\tau_{ij} $ is one-to-one identified
 to the point $P_k$ of the tensor product grid, with $k=1,...,\nu\times \nu$.
For the clarity of exposition we also introduce the disjoint set of indexes:
\noindent ${\mathcal I}_{\Omega}:=\{k | P_k\in \Omega\}$ (internal points), ${\mathcal I}_{D}:=\{k | P_k\in \Gamma_D\}$ (Dirichlet points), ${\mathcal I}_{N}:=\{k | P_k\in \Gamma_N\}$ (Neumann points),
${\mathcal I}_{AB}:=\{k | P_k\in \Gamma_{AB}\}$ (ABCs points), and ${\mathcal I}:= {\mathcal I}_{\Omega}\cup {\mathcal I}_{D} \cup {\mathcal I}_{N} \cup {\mathcal I}_{AB}$
is the set of $\nu \times \nu$ indexes of the whole mesh of collocation points.

\noindent We can now write the IGA collocation semi-discrete continuous-in-time formulation of  the acoustic problem (\ref{pde})-(\ref{ABC}). To this aim, we 
 collocate the continuous problem with   initial and  boundary conditions, at the Greville collocation points:
\begin{equation}
\frac {\partial ^2 u}{\partial t^2}(P_k,t) -c_0\Delta  u(P_k,t) =f(P_k,t),  \qquad k\in {\mathcal I}_{\Omega},\ t\in (0,T),
\label{colloc_pde}
\end{equation}
\begin{equation}
u(P_k,0)=u_0(P_k) , \qquad \frac{\partial u}{\partial t} (P_k,0)= u_1(P_k),
\qquad   {k \in}\ {\mathcal I},
\label{colloc_init_cond}
\end{equation}
\begin{equation}
 u(P_k,t) =\Phi(P_k,t),   \  k\in {\mathcal I}_{D},\ t\in (0,T), 
\qquad
   \frac{\partial u}{\partial {\bf n}}(P_k,t) =\Psi(P_k,t),
\   k \in   {\mathcal I}_{N}, \ t\in (0,T)
\label{colloc_DirNeu}
\end{equation}
\begin{equation}
\frac{1}{\sqrt{c_0}}\frac{\partial u}{\partial t}(P_k,t) + \frac{\partial u}{\partial {\bf n}}(P_k,t) =0,
\qquad k \in {\mathcal I}_{AB} , \ t\in (0,T).
\label{colloc_ABC}
\end{equation}
\noindent We observe that the semi-discrete collocation problem is equivalent to the problem of finding a vector $\mathbf{u}$ of elements $\{u_k,\ k\in \mathcal{I}\}$,
which are in correspondence with elements $\{u_{ij},\ i,j=1,...,\nu\}$. Stemming from (\ref{eq:F}) and (\ref{eq:V-k-space}), the IGA numerical solution results in
\begin{equation}
u(\mathbf{x},t)=\sum_{i,j=1}^{\nu}u_{ij}R_{ij}^{p}\circ \mathbf{F}^{-1}(\mathbf{x},t).
\label{colloc_u_sum}
\end{equation}

\noindent In order to assemble the mass and stiffness IGA matrices, we introduce the IGA collocation matrices $[D_r]$, with $r=0,\ 1,\ 2$,  accounting for $r$-th derivative, respectively, at collocation points. Precisely,
  the collocation matrices $D_0$, $D_1$ and $D_2$ are associated to identity,
 $\frac{\partial }{\partial \mathbf{n}}$ and $\Delta $ operators,  respectively.
The detailed MATLAB construction is based on the structure $\tt sp\_ eval$ of the GeoPDEs library; see, e.g., \cite{ehhr2018}, \cite{GeoPDEs} and \cite{vazquez16}. Hereinafter, we will refer to  $D_0$ as to the mass matrix and in Section \ref{Sec_Numerical_res} we will    denote it by $\cal M$.

\noindent Finally, the matrix form of the set of equations (\ref{colloc_pde})-(\ref{colloc_ABC}) can be  rewritten  
as  a system of second-order ordinary differential
equations:

\begin{equation}
\frac{\partial^2}{\partial t^2} [ D_0 \mathbf{u}(t)]_k -c_0[D_2  \mathbf{u}(t)]_k= [\mathbf{f}(t)]_k,\quad  k\in \mathcal{I}_{\Omega}
\label{system_ODE_f}
\end{equation}

\begin{equation}
[D_0 \mathbf{u}(t)]_k = [\boldsymbol{\Phi}(t)]_k,   \quad k\in \mathcal{I}_{D} ,\quad
 [D_1 \mathbf{u}(t)]_k  =  [\boldsymbol{\Psi}(t)]_k, \quad  k\in \mathcal{I}_{N}
 \label{system_ODE_DirNeu}
\end{equation}

\begin{equation}
\frac{1}{\sqrt{c_0}} \frac{\partial u}{\partial t} [D_0 \mathbf{u}(t)]_k  + [ D_1 \mathbf{u}(t)]_k = 0 ,\quad  k\in \mathcal{I}_{AB}
\label{system_ODE_g}
\end{equation}

\noindent with initial conditions
\begin{equation}
[D_0\mathbf{u}(0)]_k=[\mathbf{u}_0]_k,\quad  \frac{\partial}{\partial t}[D_0\mathbf{u}(0)]_k=[\mathbf{u}_1]_k,\quad  k\in \mathcal{I} .
\label{system_ODE_ic}
\end{equation}
\noindent Here $\mathbf{u} (t):=\{u(P_k,t),\ k\in \mathcal{I}\}$, $\mathbf{f} (t):=\{f(P_k,t),\ k\in \mathcal{I}\}$,
$ \boldsymbol{\Phi}(t):=\{\Phi(P_k,t),\ k\in \mathcal{I}_D\}$, $ \boldsymbol{\Psi}(t):=\{\Psi(P_k,t),\ k\in \mathcal{I}_N\}$,
$\mathbf{u}_0 :=\{u_0(P_k),\ k\in \mathcal{I}\}$,
$\mathbf{u}_1 :=\{u_1(P_k),\ k\in \mathcal{I}\}$, with all vectors  assigned equal to zero elsewhere.

\subsection{Time discretization of the wave model problem (\ref{pde})-(\ref{ABC}) by Newmark advancing schemes}\label{subsec_Newmark}

For the discretization of time derivatives in (\ref{system_ODE_f}), (\ref{system_ODE_g}) and (\ref{system_ODE_ic}), we introduce the   finite difference Newmark scheme \cite{New59}, that, in its general form, reads:
\begin{equation}  {\bf u}_{n+1}={\bf u}_{n}+\Delta t\ {\bf v}_{n} + (1-2\beta)\Delta t^2 {\bf a}_{n}/2 + \beta \Delta t^2{\bf a}_{n+1},\ \ \ {\bf v}_{n+1}={\bf v}_{n}+(1-\gamma)\Delta t\ {\bf a}_{n} + \gamma\Delta t {\bf a}_{n+1},\label{Newmark_uv}
\end{equation}

\noindent where ${\bf u}_{n}:=\{u(P_k,t_n),\ k\in \mathcal{I}\}$, ${\bf v}_{n}:=\{v(P_k,t_n),\ k\in \mathcal{I}\}$, ${\bf a}_{n}:=\{a(P_k,t_n),\ k\in \mathcal{I}\}$ are the  vectors of approximated displacement, velocity and acceleration, respectively, at time $t_n$, having partitioned the  interval $(0,T)$   into $N$ subintervals $[t_n, t_{n+1}]  $, with $t_0=0$, $t_N=T$, $t_{n+1}=t_n+\Delta t$, $n=0,...,N-1$ and $\Delta t=T/N$.  
Moreover, $\beta\ge 0$ and $\gamma\ge 0$ are real parameters.  It can be shown (see, e.g., \cite{wood84},\cite{wood90}) that we can eliminate the velocity and acceleration vectors and express
the Newmark method as a two-step scheme in the displacement term  $\mathbf{u}_n$ only, whose entries   give  the corresponding IGA solution at  time step $t_n$, according to (\ref{colloc_u_sum}).  The     initial vector ${\bf u}_{1}$ at the second time instant $t_1=t_0+\Delta t$, can be computed, for example, from the first one ${\bf u}_{0}$ associated to initial conditions (\ref{init_cond}) by means of a second order explicit one-step method, (e.g., an explicit two-stage  Runge-Kutta method), in order to  preserve the global accuracy of the scheme.

\noindent If we denote  by $[D_r]_{k}$ the $k$-th row of the collocation matrix $D_r$, $r=0,1,2$, by $[\mathbf{w}]_k$ the $k$-th element of a general vector $\mathbf{w}$, and apply  the Newmark scheme (\ref{Newmark_uv})  to the numerical solution of the acoustic wave IGA collocation problem  (\ref{system_ODE_f})-(\ref{system_ODE_ic}),  we obtain
   the set of recurrence relations for the displacement term  $\mathbf{u}_n$ at collocation points:

\begin{equation}\label{Newmark_Omega}
[D_0]_{k }\frac{ \mathbf{u}_{n+1}-2\mathbf{u}_{n}+\mathbf{u}_{n-1}}{\Delta t^2}-c_0
[D_2]_{k }\Big[\beta \mathbf{u}_{n+1}+\Big(\frac12-2\beta+\gamma\Big)\mathbf{u}_{n}+\Big(\frac12+\beta-\gamma\Big){\bf u}_{n-1}\Big]=
\end{equation}

$$\Big[\beta \mathbf{f}_{n+1}+\Big(\frac12-2\beta+\gamma\Big)\mathbf{f}_{n}+\Big(\frac12+\beta-\gamma\Big){\bf f}_{n-1}\Big]_k, \quad k\in \mathcal{I}_{\Omega}$$

\begin{equation}\label{Newmark_DirNeu}
[D_0]_{k }  \mathbf{u}_{n+1} = [\boldsymbol{\Phi}(t_{n+1})]_k, \quad k\in \mathcal{I}_{D},\quad
[D_1]_{k }  \mathbf{u}_{n+1} = [\boldsymbol{\Psi}(t_{n+1})]_k, \quad k\in \mathcal{I}_{N}
\end{equation}

\begin{equation}\label{Newmark_AB}
\frac{1}{\sqrt{c_0}}[D_0]_{k}\frac{ \gamma \mathbf{u}_{n+1}+(1-2\gamma)\mathbf{u}_{n}+(\gamma-1)\mathbf{u}_{n-1}}{\Delta t}+
[D_1]_{k }  \mathbf{u}_{n+1} =0, \quad k\in \mathcal{I}_{AB} .
\end{equation}

\noindent At   corner points involving   ABCs and/or Neumann boundary conditions we enforce the average of normal derivatives, whereas Dirichlet conditions   overrides Neumann or ABC ones. With a suitable generalizations of coefficients multiplying matrices $D_0$ and $D_1$ in (\ref{Newmark_AB}), we note that the approximation of ABCs is mathematically equivalent to that of the most common Robin boundary conditions, involving a linear combination of the values of a function and of its normal derivative at the collocation points of the boundary. 

%\noindent It can be proven  that the Newmark method is first-order accurate with respect to $\Delta t$ %if $\gamma \neq \frac12$, and it is second order if $\gamma=\frac12$.
% If $\beta=0$ and $\gamma=\frac12$ the Newmark scheme coincides with the \textit{Leap-Frog} method, %which in particular is  second-order accurate with respect to $\Delta t$.
%See, e.g., \cite{New59} and \cite{wood84}, for a review of the mathematical background and convergence %and stability analysis for the Newmark method. 

\noindent The scheme (\ref{Newmark_Omega})  is customarily considered explicit when $\beta=0$, even if  the matrix $D_0$   is not diagonal. 
 More generally, regardless of the parameter $\beta$,  we observe that each step of (\ref{Newmark_Omega})-(\ref{Newmark_AB}) involves  the resolution of a linear system $\mathcal{K} \mathbf{u}_{n+1}=\mathbf{\Upsilon}(t_{n+1},t_n,t_{n-1})$, where the matrix
\begin{equation}\label{global_matrix}
\mathcal{K}=[D_0]_{k\in\mathcal{I}_D}+
[D_1]_{k\in\mathcal{I}_N\cup\mathcal{I}_{AB}}+
\frac{\gamma}{ \Delta t\sqrt{c_0}}[D_0]_{k\in\mathcal{I}_{AB}}-c_0\beta [D_2]_{k \in  \mathcal{I}_{\Omega} }
+\frac{1}{ \Delta t^2}[D_0]_{k\in\mathcal{I}_{\Omega}},
\end{equation}
is full and non-symmetric both for explicit ($\beta =0$) and implicit ($\beta \neq 0$) methods.  From now on we will refer to $\cal K$ as stiffness matrix.
Finally, the right term
$
 \boldsymbol{\Upsilon}(t_{n+1},t_n,t_{n-1})$
accounts  for the values of data functions $u_0$, $u_1$, $f$, $\Phi$ and $\Psi$  at times $t_{n+1},t_n,t_{n-1}$.  

%\noindent Furthermore, the Newmark method  is unconditionally stable for $\gamma \ge \frac12$ and %$\beta\ge \frac14 (\frac 12 +\gamma)^2$.
%Then, the time step $\Delta t$ is subject to a stability condition in the case of explicit methods,
%whereas implicit methods can be unconditionally stable regardless of the time step $\Delta t$ for %suitable choices of the parameters $\beta$ and $\gamma$.

\section{Condition number estimates}
\label{cond_estimates}

The following condition number estimates for the Galerkin isogeometric mass and stiffness matrices have been proven in \cite{gahalaut2012} in the two dimensional case, independently of the $k$-regularity of the spline basis functions
\begin{eqnarray}
cond(\cal M) & \approx & c p^2 16^p, \ \ \mbox{with $c$ independent of $h$ and $p$,}\\
cond(\cal K) & \leq &c(h) p^8 16^p.
\end{eqnarray}
Additional bounds on the smallest and largest eigenvalues are proven in \cite{gahalaut2012}, and also in \cite{garoni2014} in the one dimensional case.
A numerical study of the conditioning of Galerkin isogeometric mass and stiffness matrices in $d$ dimensions has been carried out in \cite{iga-sem2018}, reporting the following more detailed and sharper estimates:
%\begin{align}
\begin{eqnarray}
& \mbox{for $k=0$ regularity: }\ \ & cond({\cal M})  \approx  p^{-d/2} 4^{pd},\label{Mk0}\\
& \mbox{for $k=p-1$ regularity: }\ \ 
& cond(\cal M)  \approx \Big\{
\begin{array}{ll}
 e^{pd}    & \mbox{if $h \leq 1/p$} \\
(e/4)^{d/h}(hp)^{-d/2}4^{pd}     & \mbox{otherwise,}\label{Mkpm1}
\end{array}
\end{eqnarray}
%\end{align}
\begin{eqnarray}
& \mbox{for $k=0$ regularity: }\ \ & cond(\cal K)  \approx  \Big\{
\begin{array}{ll}
h^{-2}p^2    & \mbox{if $h \leq (p^{2+d/2}d^{-dp})^{1/2}$} \\
p^{-d/2}4^{pd}     & \mbox{otherwise,}\label{Kk0}
\end{array}\\
& \mbox{for $k=p-1$ regularity: }\ \ 
& cond(\cal K)  \approx \Bigg\{
\begin{array}{ll}
h^{-2}p  & \mbox{if $h\leq e^{-dp/2}$} \\
pe^{pd}    & \mbox{if $e^{-dp/2} \leq h \leq 1/p$} \\
(e/4)^{d/h}p^{-d/2}h^{-d/2-1}4^{pd}     & \mbox{otherwise.}\label{Kkmp1}
\end{array}
\end{eqnarray}

The numerical results reported in the next section show that these sharper estimates
(\ref{Mk0})-(\ref{Kkmp1}) are also upper bounds for the condition number of the collocation isogeometric mass and stiffness matrices considered in this work, but in some cases the collocation matrices satisfy improved bounds. In particular, the term $4^{pd}$ appearing in these estimates seem to improve in the two-dimensional case ($d=2$) to at least $4^{\frac{3}{2}p}$, and in some cases for sufficiently small $h$ to $4^{p}$ or even 
$4^{\frac{1}{2}p}$.

\section{Numerical results}
\label{Sec_Numerical_res}
In this Section, we present a numerical study of the behaviour of the eigenvalues and condition numbers of the mass matrix $\cal M$ and stiffness matrix $\cal K$ for the acoustic wave problem in the reference square domain $\Omega=[0,1]\times[0,1]$ with different types of boundary conditions. We vary the degree $p$, regularity $k$ and mesh size $h$ for the IGA collocation method introduced in Section \ref{subsec_IGA_Coll}, and the discretization parameters $\Delta t$ and $\beta$ of the Newmark time advancing scheme introduced in Section \ref{subsec_Newmark}. We denote by {\textsf { d.o.f.}} the number of degrees of freedom of the discrete problem and by {\textsf  {nz}} the number of nonzero entries in the mass and stiffness matrices. All tests are implemented in MATLAB R2020b by the GeoPDEs  3.0 library \cite{GeoPDEs,vazquez16}, and the condition numbers are computed using the MATLAB $\tt condest$ function.

\begin{figure} %[!t]
\vspace{-10mm}
\centerline{
\begin{tabular}{cc}
\includegraphics[scale=0.5]{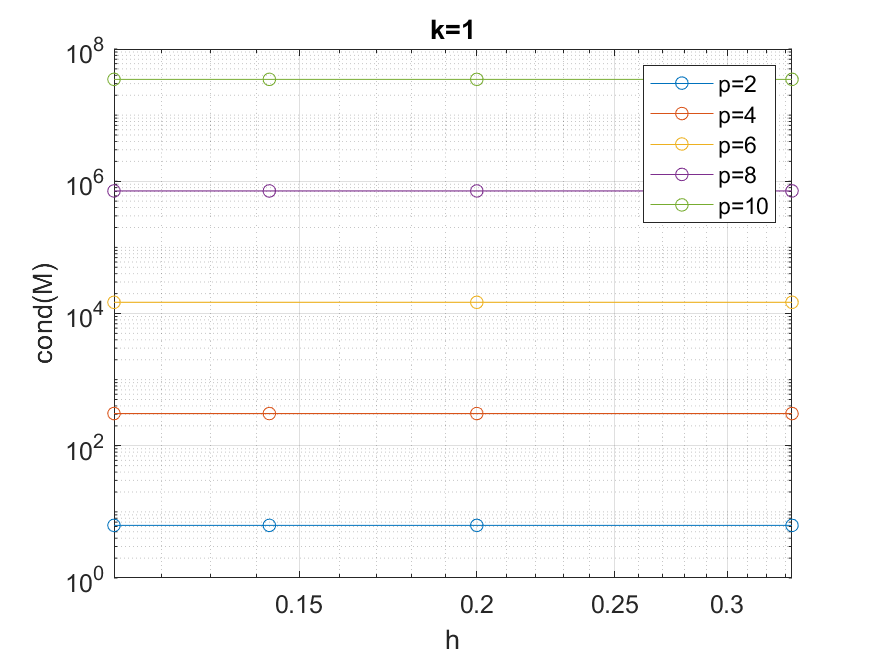} &
\includegraphics[scale=0.5]{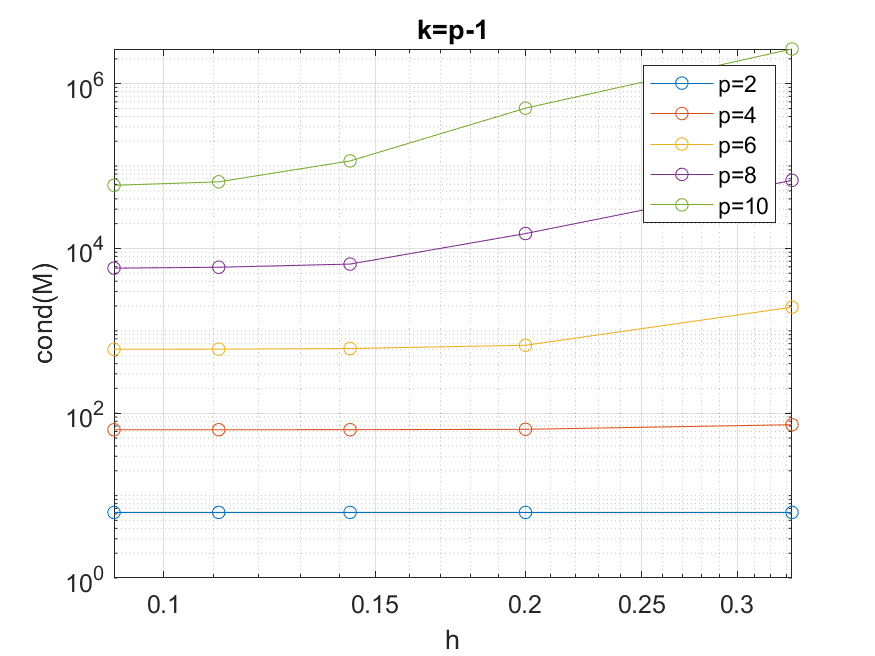} \\
\includegraphics[scale=0.5]{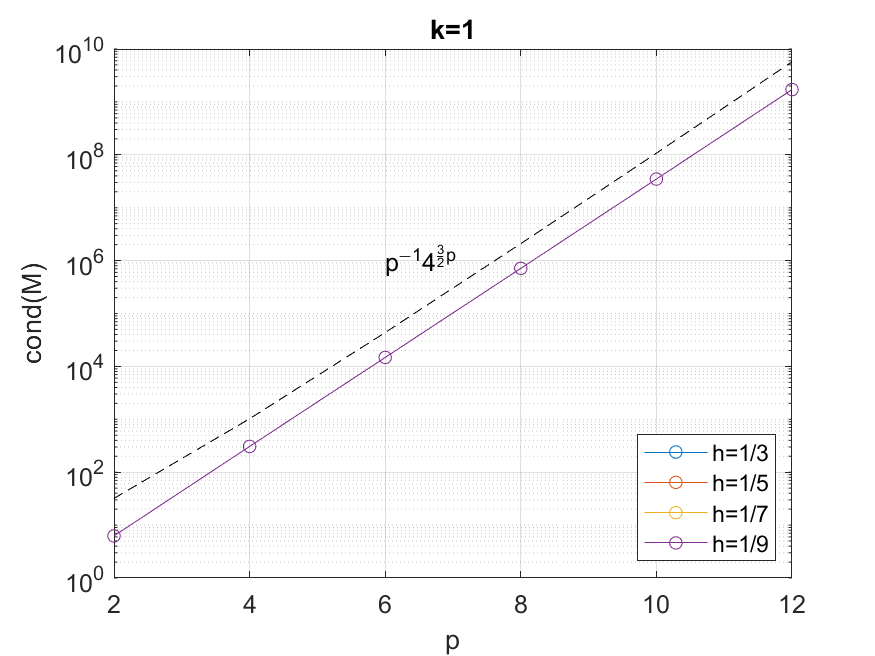} & \includegraphics[scale=0.5]{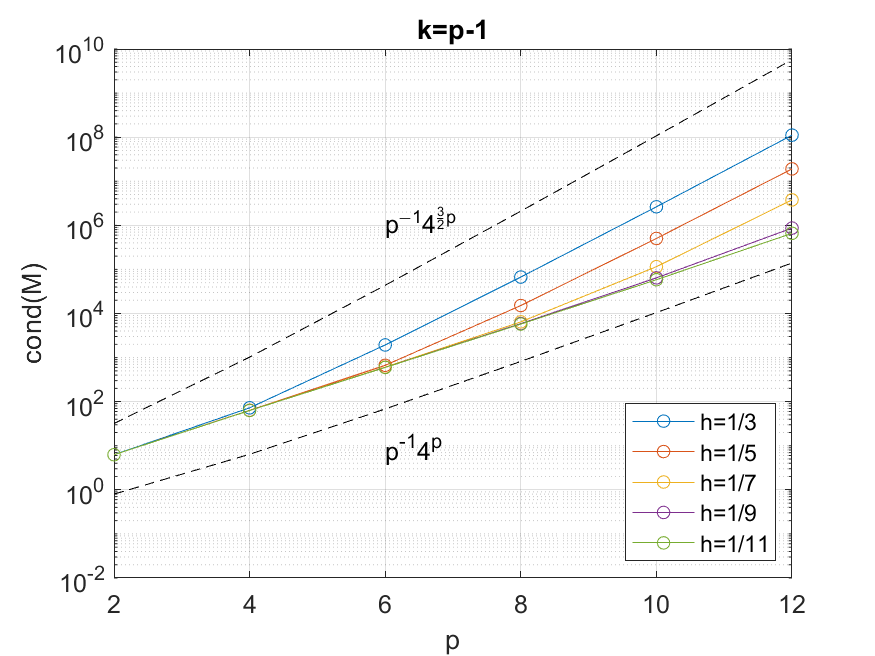} \\
\multicolumn{2}{c}{\includegraphics[scale=0.5]{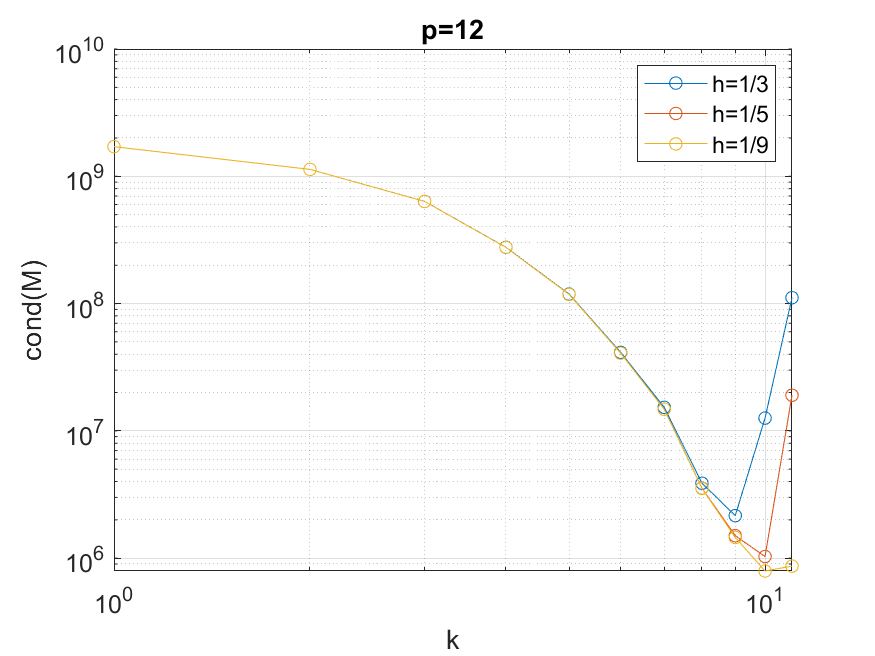}}\\
\end{tabular}
}
\caption{Condition number of the mass matrix versus: $h$ (top), for $p=2,4,6,8,10$,  $k=1$ (left) or  $k=p-1$ (right); $p$ (center), for $h=1/3,1/5,1/7,1/9$,   $k=1$ (left) or $h=1/3,1/5,1/7,1/9,1/11$,   $k=p-1$ (right);    $k$ (bottom), for $h=1/3,1/5, 1/9$, fixed $p=12$.
\label{cond_M_k}}
\end{figure}

{\bf{Eigenvalues and condition number of the mass matrix.}} 
%\label{num_res_mass}
In Fig. \ref{cond_M_k},  we report the condition numbers {\textsf{cond}}($\cal M$) versus: the mesh size $h$ (top), with five different values of degree $p$ and minimal regularity $k=1$ (left) or maximal regularity  $k=p-1$ (right); the degree $p$ (center), with four different values of mesh size $h$ and minimal regularity $k=1$ (left) or five different values of mesh size $h$ and maximal regularity  $k=p-1$ (right); the regularity $k$ (bottom), with three different values of mesh size $h$, fixed $p=12$. For the $h$-refinement (top) the condition numbers {\textsf{cond}}($\cal M$) are independent of $h$, whereas they grow slower than the estimates (\ref{Mk0})-(\ref{Mkpm1}) for $p$- refinement. For increasing $k$, fixed $p$ and $h$, {\textsf{cond}}($\cal M$) seem to decrease exponentially except when the regularity $k$ approaches the maximum value $p-1$ and {\textsf{cond}}($\cal M$) increases sharply.

\begin{figure} %[!t]
%\vspace{10mm}
\centerline{
\begin{tabular}{ccc}
\includegraphics[scale=0.40]{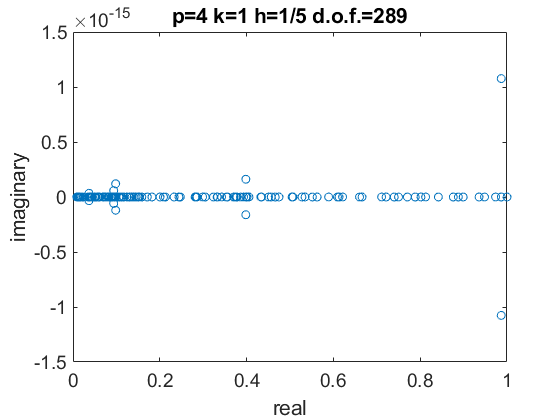} &
\hspace{-5mm}\includegraphics[scale=0.40]{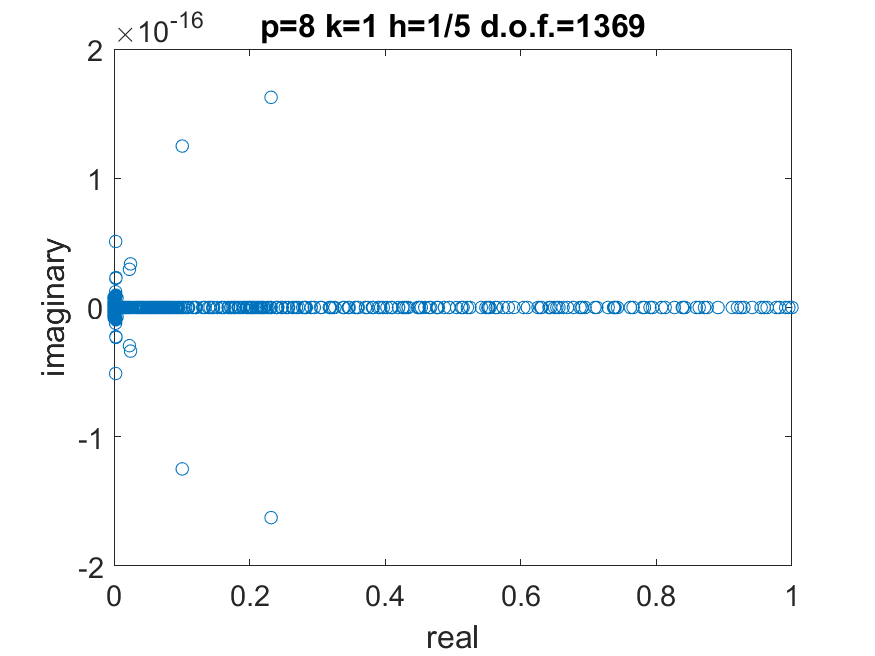} &
\hspace{-5mm}\includegraphics[scale=0.40]{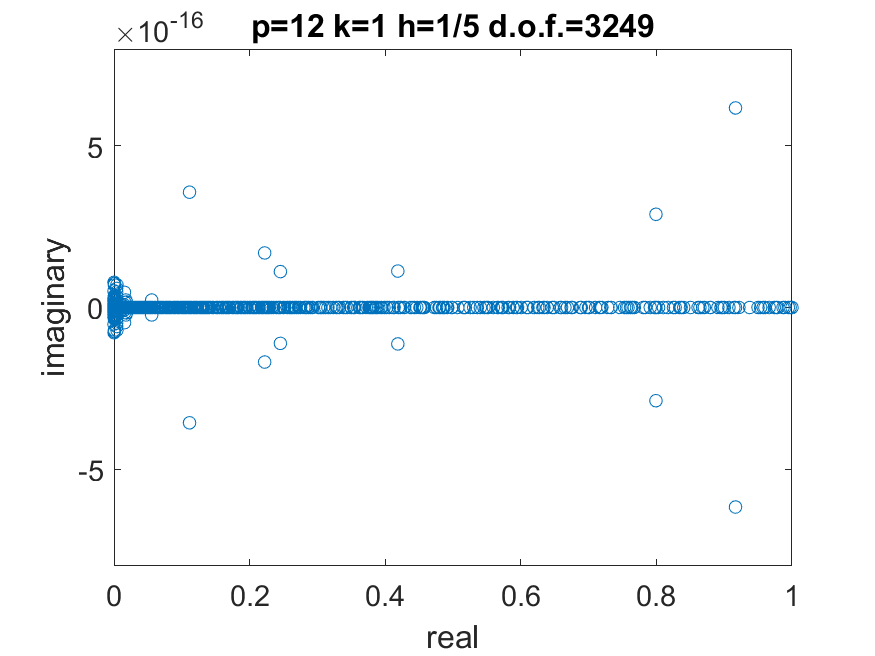} \\
\includegraphics[scale=0.40]{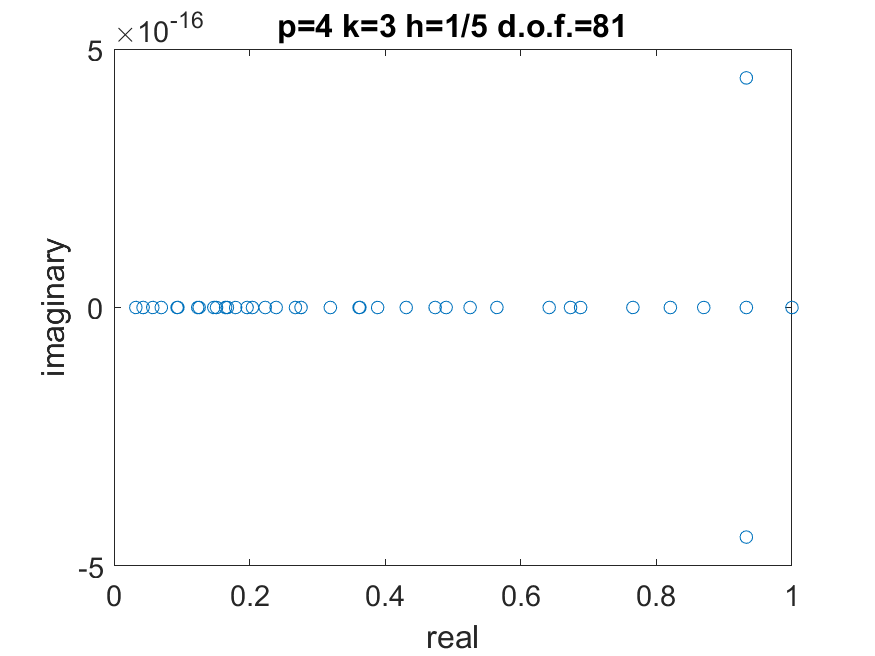} & 
\hspace{-5mm}\includegraphics[scale=0.40]{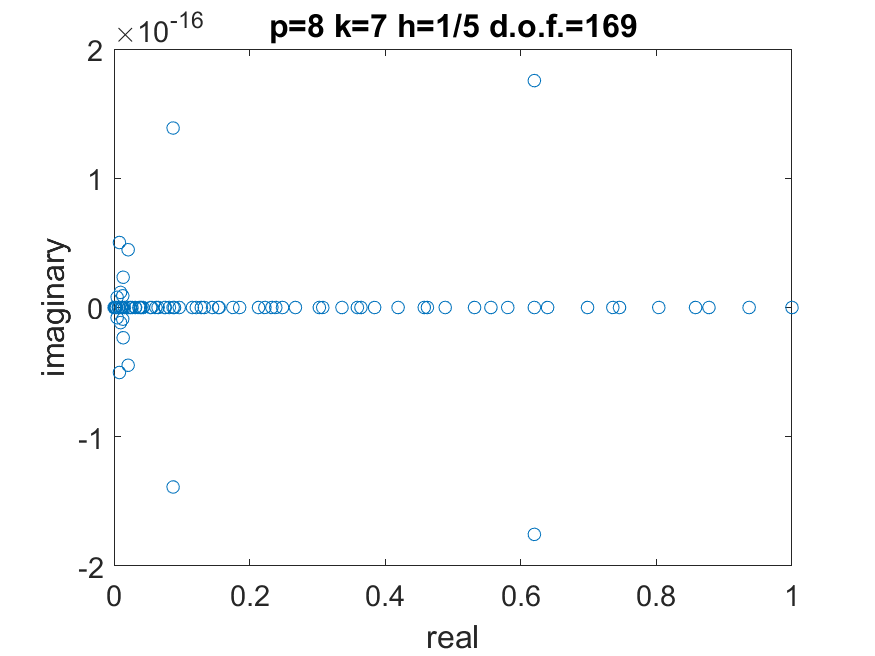} & 
\hspace{-5mm}\includegraphics[scale=0.40]{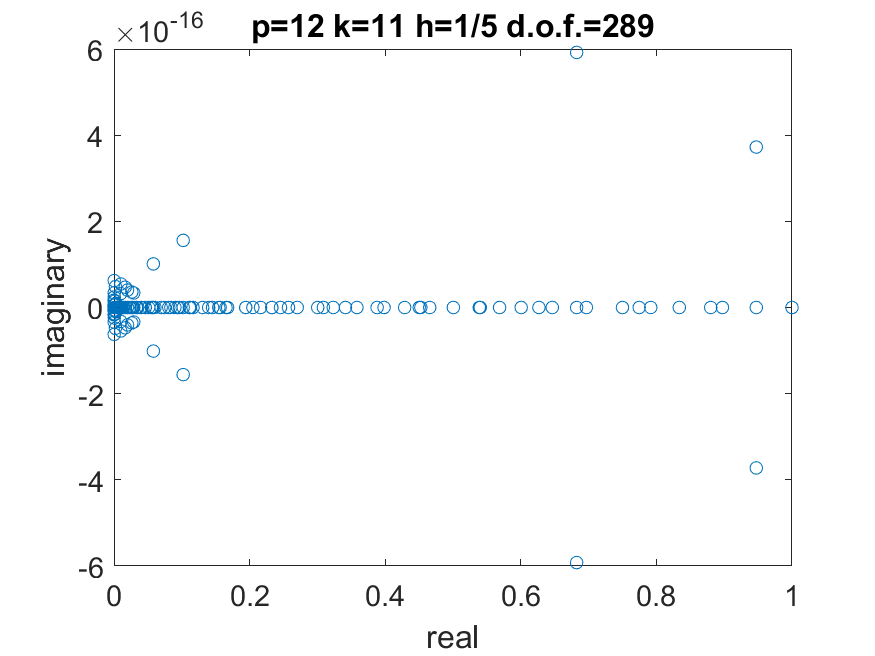}  \\
\end{tabular}
}
\vspace{-2mm}
\caption{Mass matrix eigenvalue distribution in the complex plane, for $p = 4$ (left),
 $p=8$ (center), $p=12$ (right), with $k=1$ (top) or $k=p-1$ (bottom), fixed $h=1/5$.
\label{eig_M_h5_VSp}}
\end{figure}

\begin{figure} %[!t]
\vspace{4mm}
\centerline{
\begin{tabular}{ccc}
\includegraphics[scale=0.40]{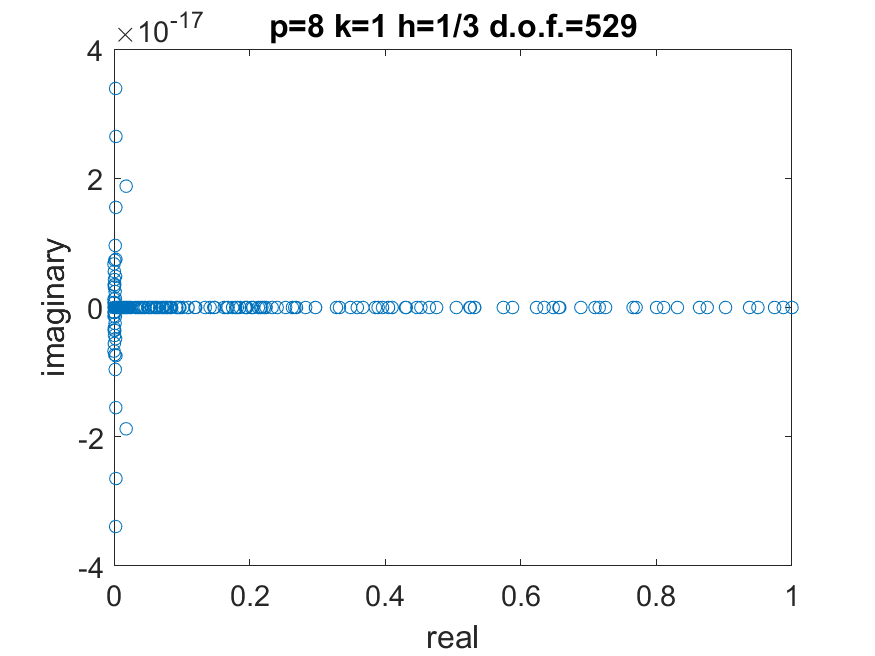} &
\hspace{-5mm}\includegraphics[scale=0.40]{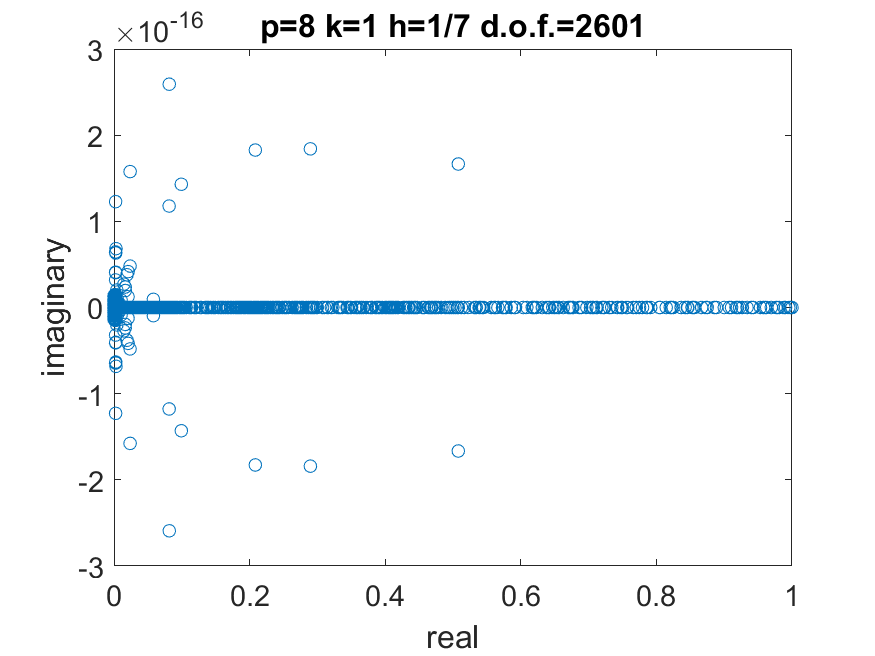} &
\hspace{-5mm}\includegraphics[scale=0.40]{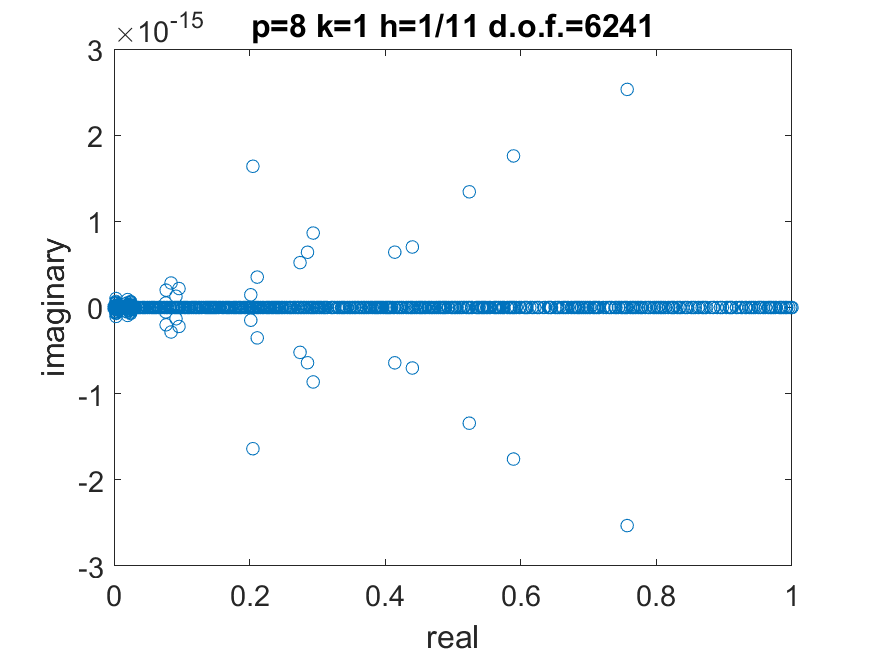} \\
\includegraphics[scale=0.40]{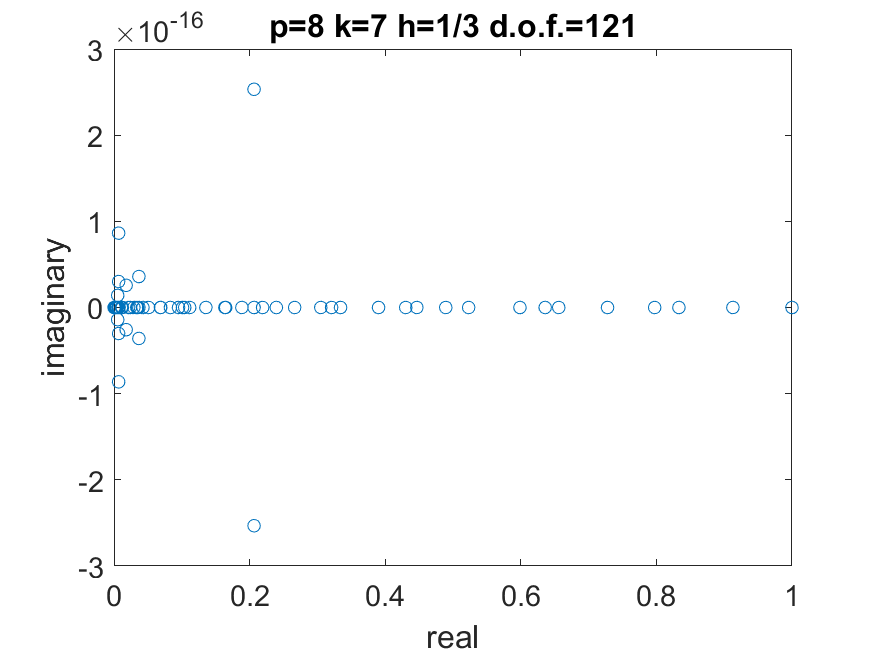} & 
\hspace{-5mm}\includegraphics[scale=0.40]{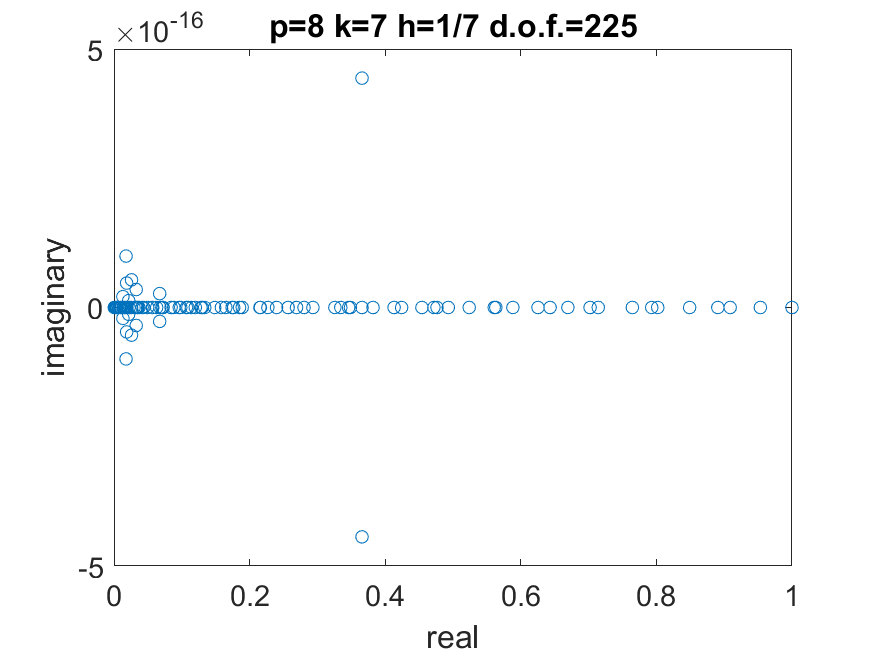} & 
\hspace{-5mm}\includegraphics[scale=0.40]{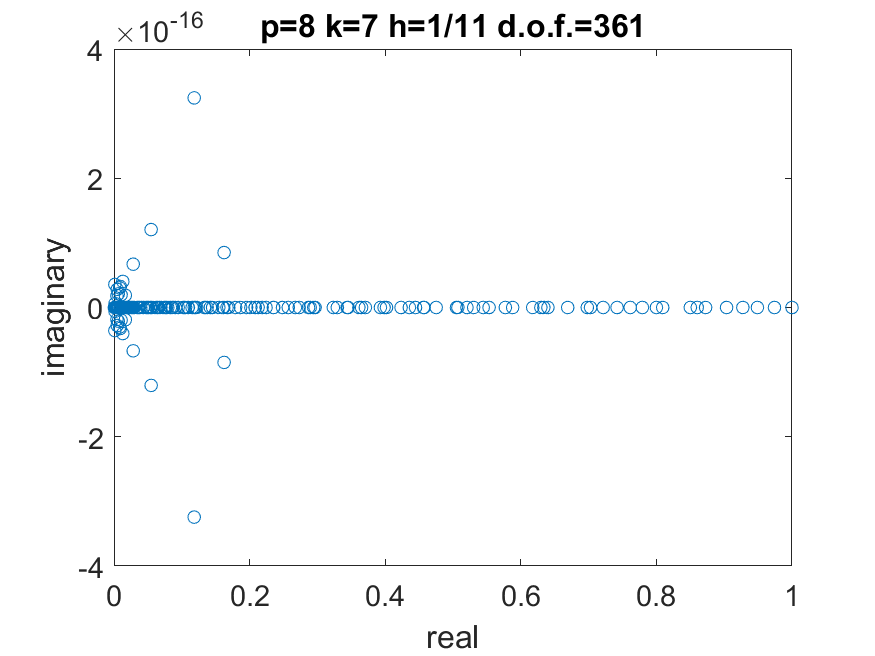}  \\
\end{tabular}
}
\vspace{-2mm}
\caption{Mass matrix eigenvalue distribution in the complex plane, for $h=1/3$ (left),
 $1/h=7$ (center), $1/h=11$ (right), with $k=1$ (top) or $k=p-1$ (bottom), fixed $p=8$.
\label{eig_M_p8_VSh}}
\end{figure}

Fig. \ref{eig_M_h5_VSp} reports the mass matrix $\cal M$ eigenvalue distribution in the complex plane for three different values of degree $p=4$ (left), $p=8$ (center), $p=12$ (right), minimal regularity $k=1$ (top) or maximal regularity $k=p-1$ (bottom), for fixed $h=1/5$, whereas in Fig. \ref{eig_M_p8_VSh} we report the eigenvalue distribution %in complex plane of mass matrix $\cal M$%
for three different values of mesh size $h=1/3$ (left), $h=1/7$ (center), $h=1/11$ (right), minimal regularity $k=1$ (top) or maximal regularity $k=p-1$ (bottom), for fixed $p=8$. In all cases eigenvalues are essentially real since their imaginary parts are of the order of machine precision.  

\begin{figure} %[!t]
%\vspace{10mm}
\centerline{
\begin{tabular}{ccc}
\includegraphics[scale=0.40]{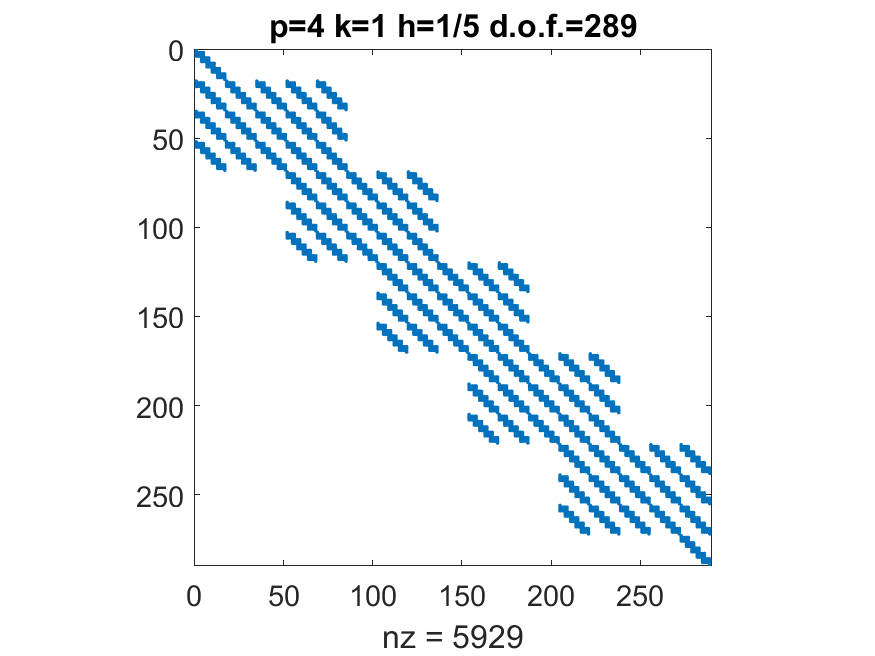} &
\hspace{-5mm}\includegraphics[scale=0.40]{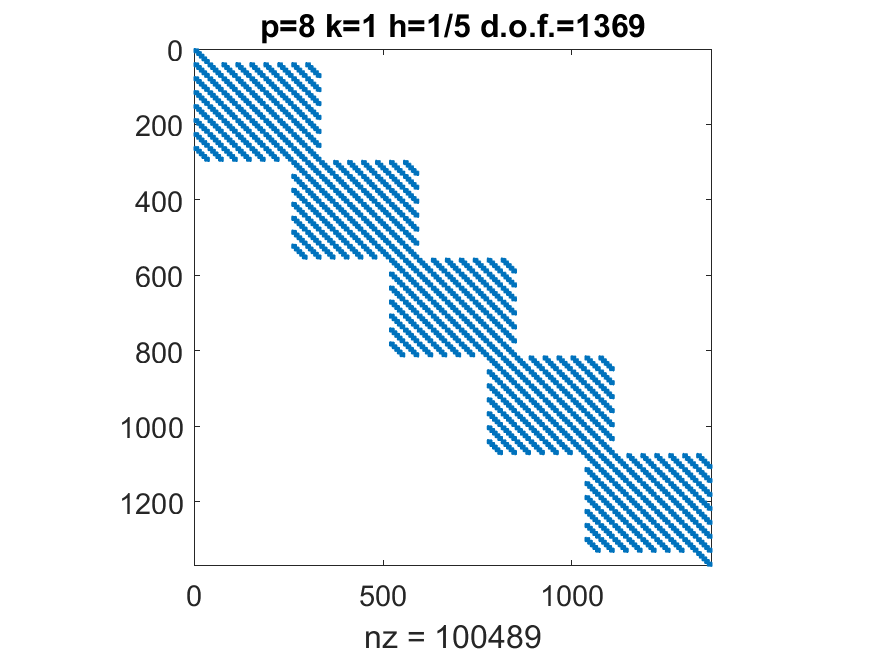} &
\hspace{-5mm}\includegraphics[scale=0.40]{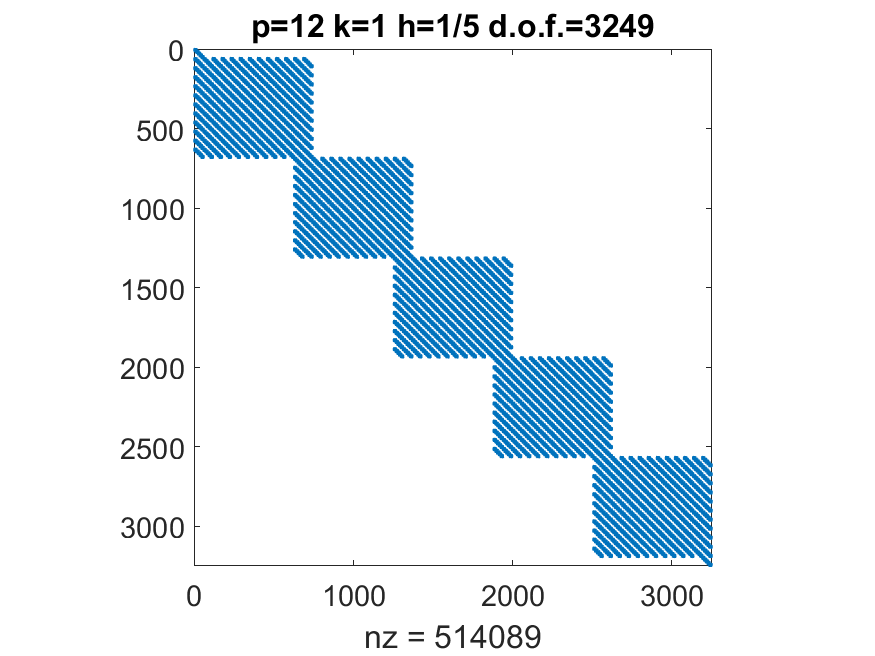} \\
\includegraphics[scale=0.40]{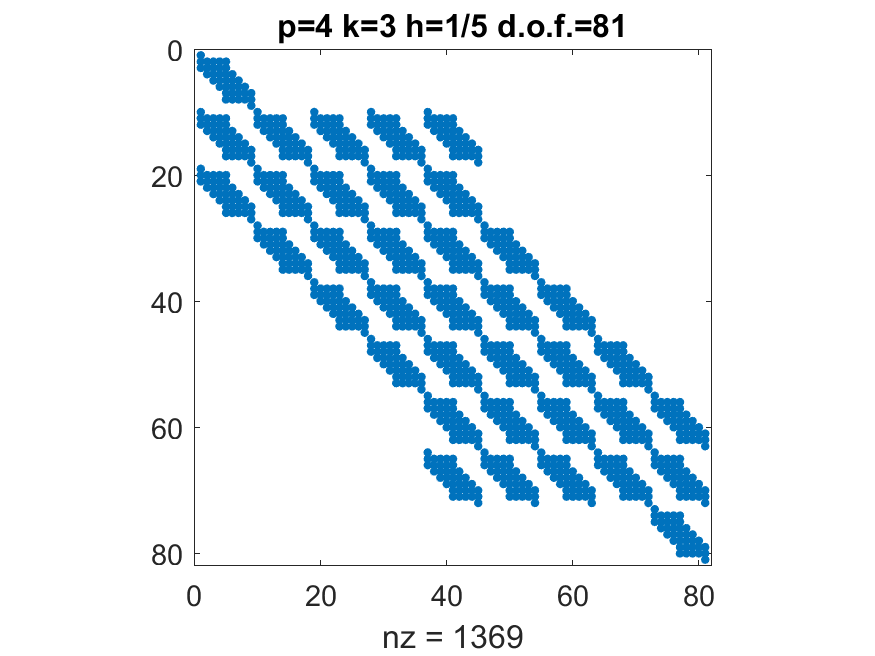} & 
\hspace{-5mm}\includegraphics[scale=0.40]{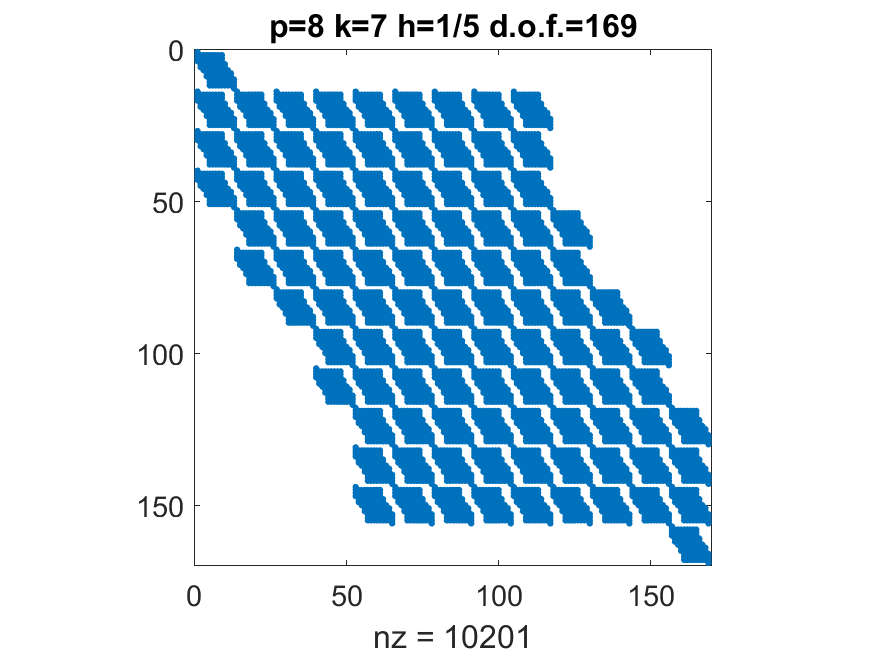} & 
\hspace{-5mm}\includegraphics[scale=0.40]{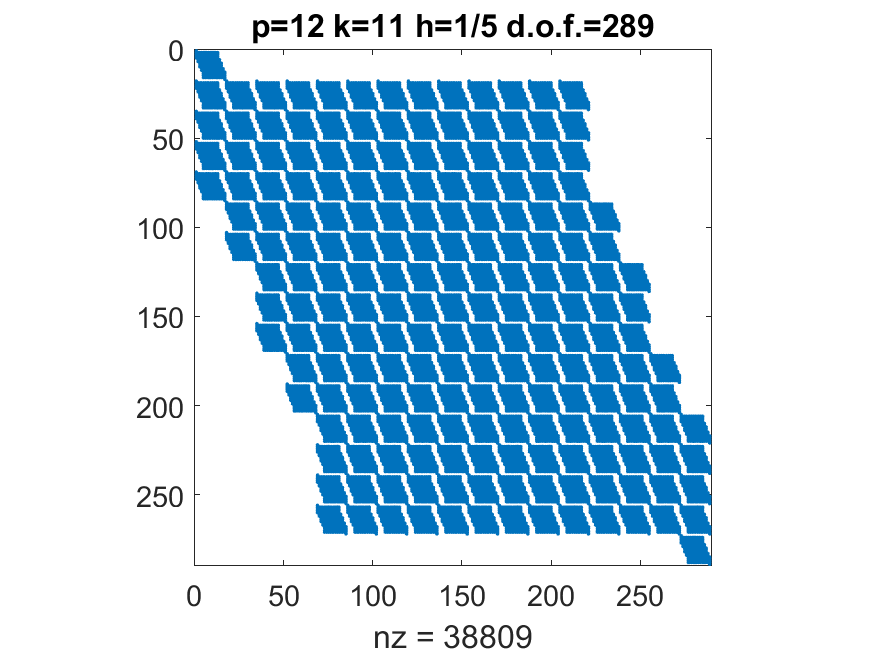}  \\
\end{tabular}
}
\vspace{-2mm}
\caption{Mass matrix sparsity pattern, for $p = 4$ (left),
 $p=8$ (center), $p=12$ (right), with $k=1$ (top) or $k=p-1$ (bottom), fixed $h=1/5$.
\label{spy_M_h5_VSp}}
\end{figure}

\begin{figure} %[!t]
%\vspace{10mm}
\centerline{
\begin{tabular}{ccc}
\includegraphics[scale=0.40]{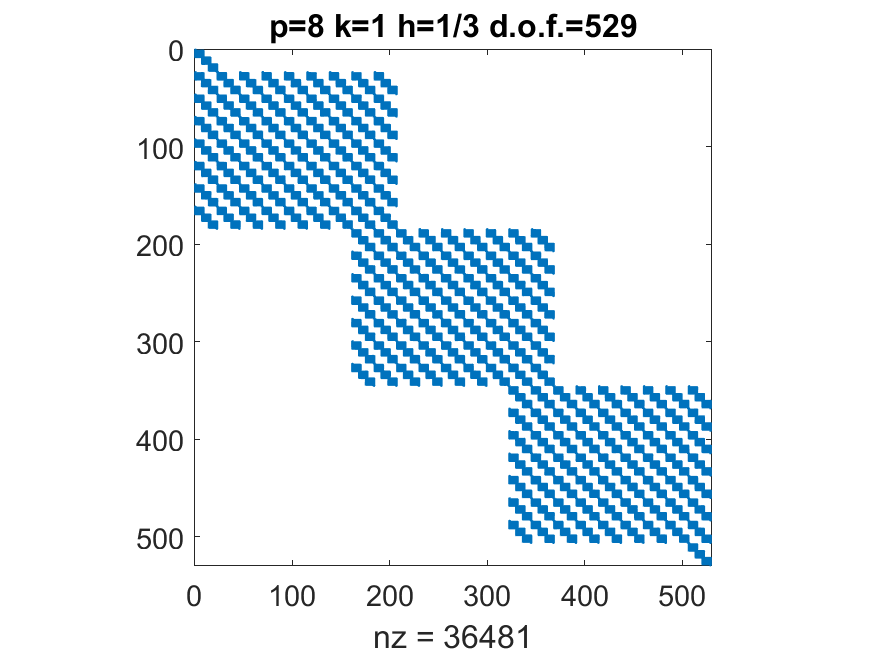} &
\hspace{-5mm}\includegraphics[scale=0.40]{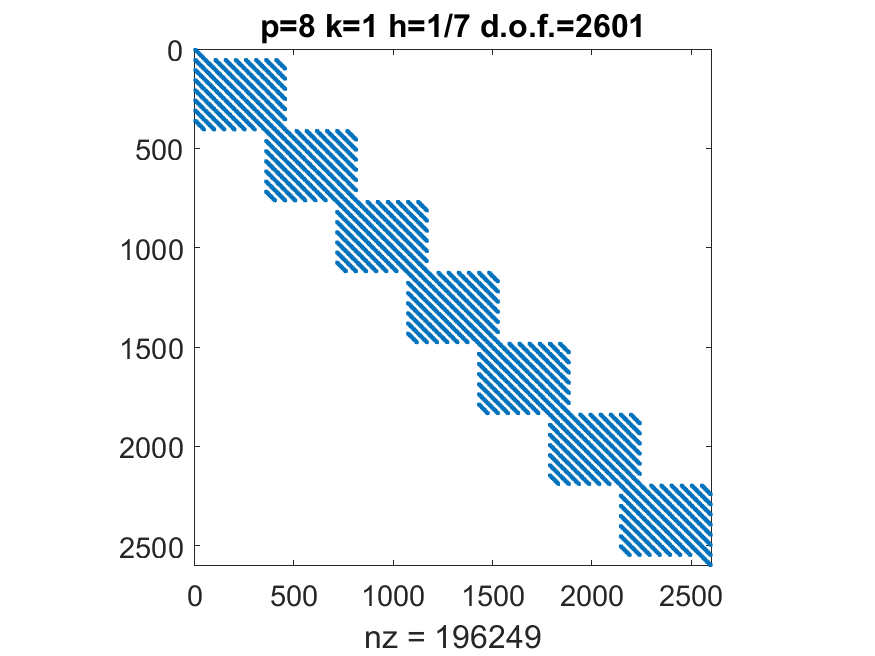} &
\hspace{-5mm}\includegraphics[scale=0.40]{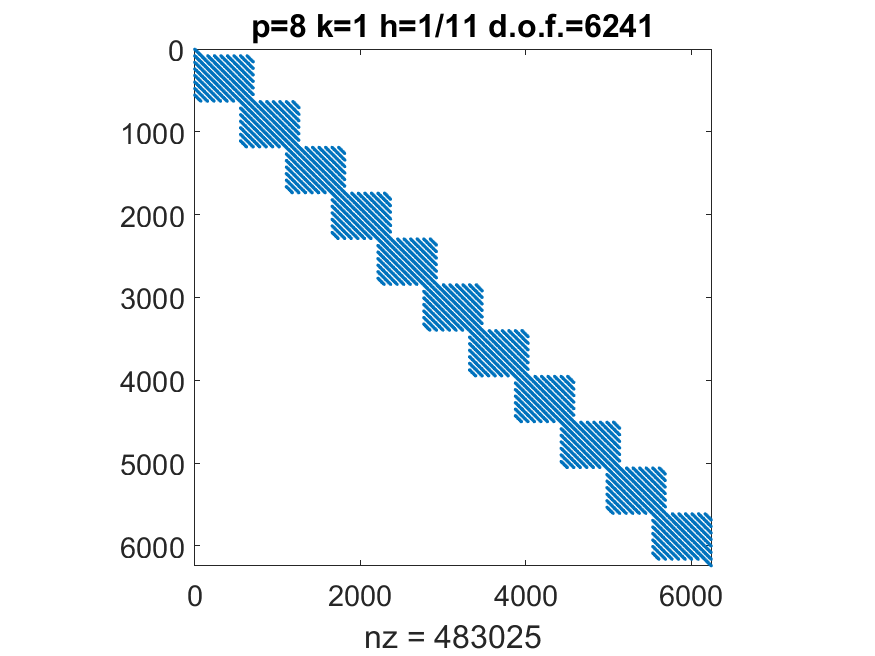} \\
\includegraphics[scale=0.40]{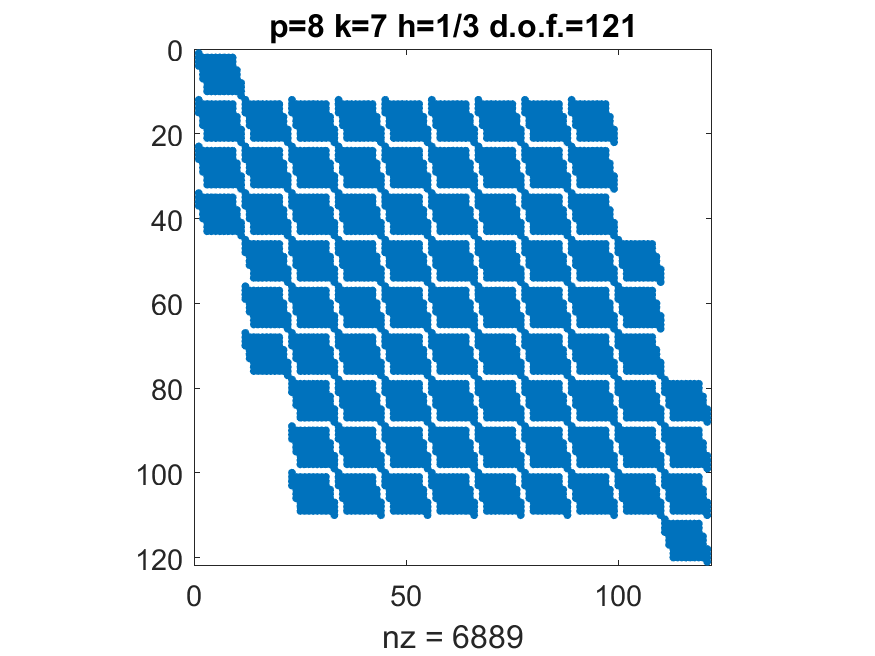} & 
\hspace{-5mm}\includegraphics[scale=0.40]{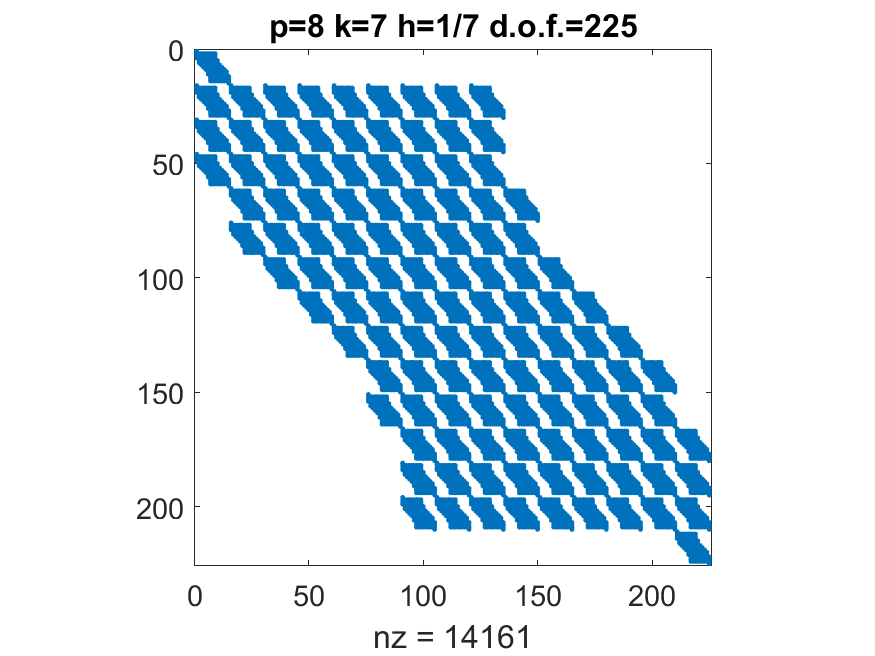} & 
\hspace{-5mm}\includegraphics[scale=0.40]{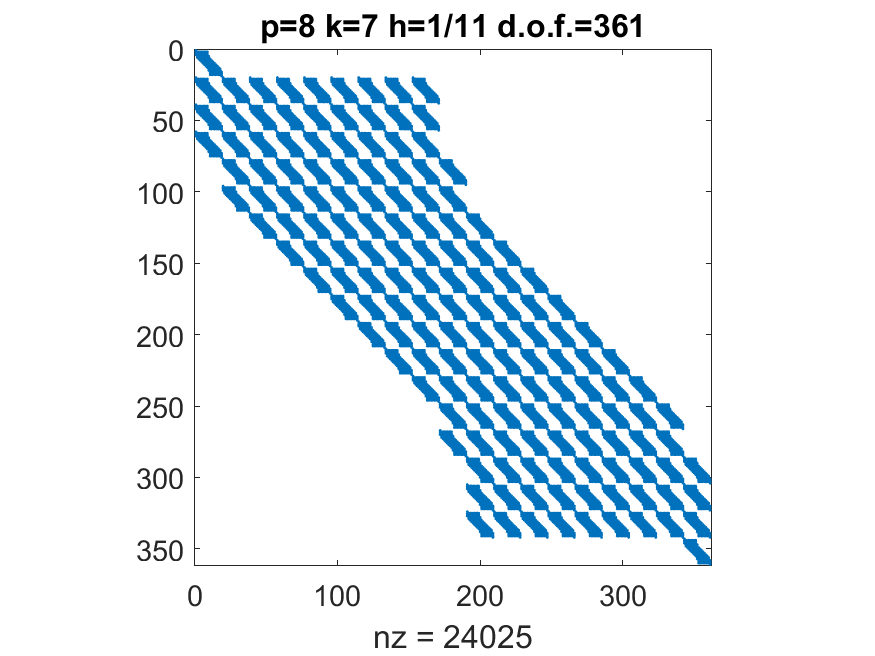} \\
\end{tabular}
}
\vspace{-2mm}
\caption{Mass matrix sparsity pattern, for $h=1/3$ (left),
 $1/h=7$ (center), $1/h=11$ (right), with $k=1$ (top) or $k=p-1$ (bottom), fixed $p=8$.
\label{spy_M_p8_VSh}}
\end{figure}

Figs. \ref{spy_M_h5_VSp} and \ref{spy_M_p8_VSh} show the sparsity pattern of the mass matrix $\cal M$  for the same values of parameters of Figs. \ref{eig_M_h5_VSp} and \ref{eig_M_p8_VSh}, respectively. 
%Matrices are block-diagonal for minimal regularity and almost full in the case of maximal regularity. 
The number of {\textsf { d.o.f.}} and nonzero elements {\textsf  {nz}} decrease considerably when the regularity is increased from minimal to maximal, and the difference grows when $p$ or $1/h$ are increased.
%are at least an order of magnitude larger than the ones of the second case.

{\bf{Eigenvalues and condition number of the stiffness matrix with Dirichlet and Neumann boundary conditions. }}
%\label{stiffDIR_NEU}

\begin{figure} %[!t]
\vspace{-10mm}
\centerline{
\begin{tabular}{cc}
\includegraphics[scale=0.485]{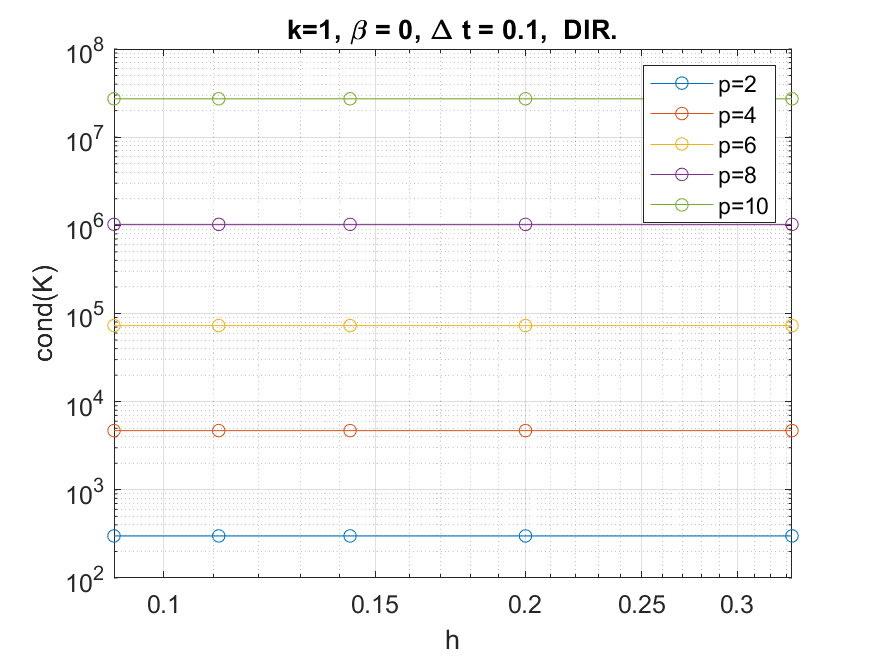} &
\includegraphics[scale=0.485]{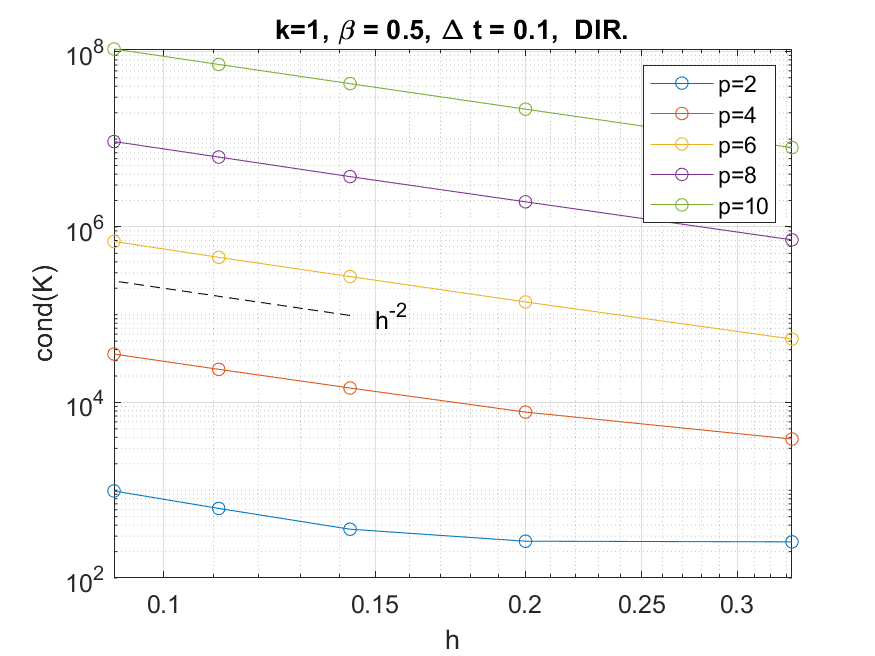} \\
\includegraphics[scale=0.485]{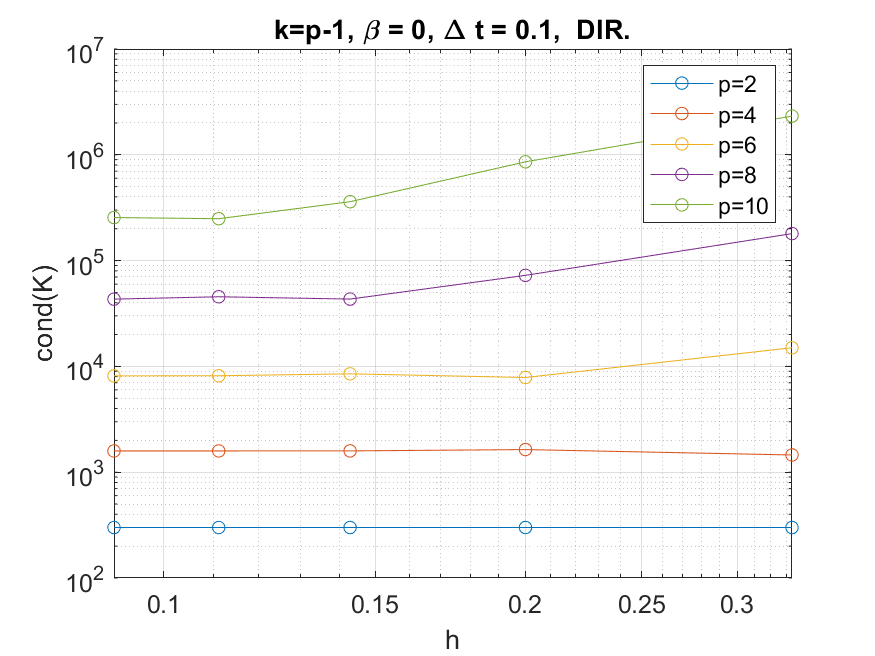} & \includegraphics[scale=0.485]{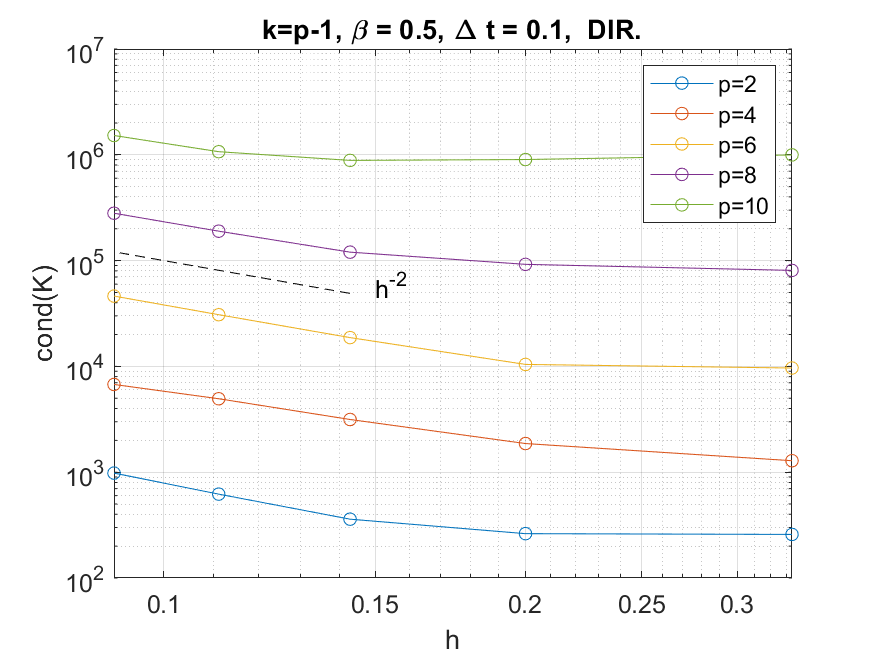} \\
\includegraphics[scale=0.485]{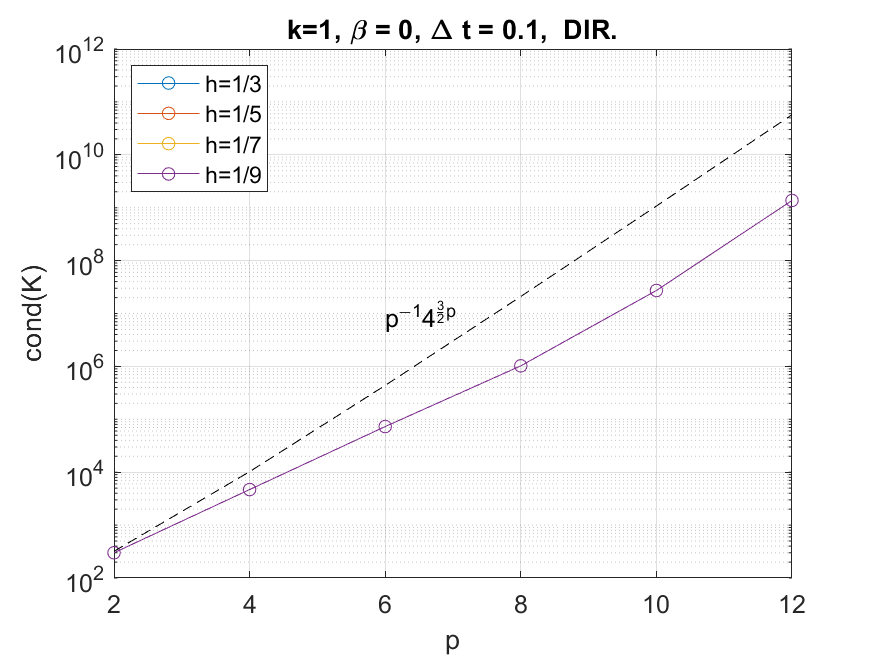} & 
\includegraphics[scale=0.485]{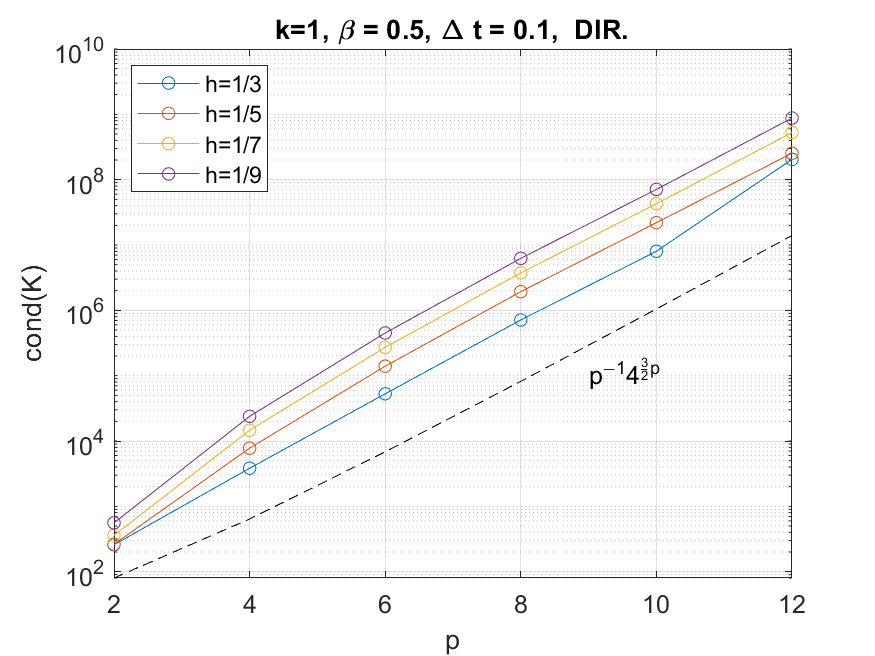} \\
\includegraphics[scale=0.485]{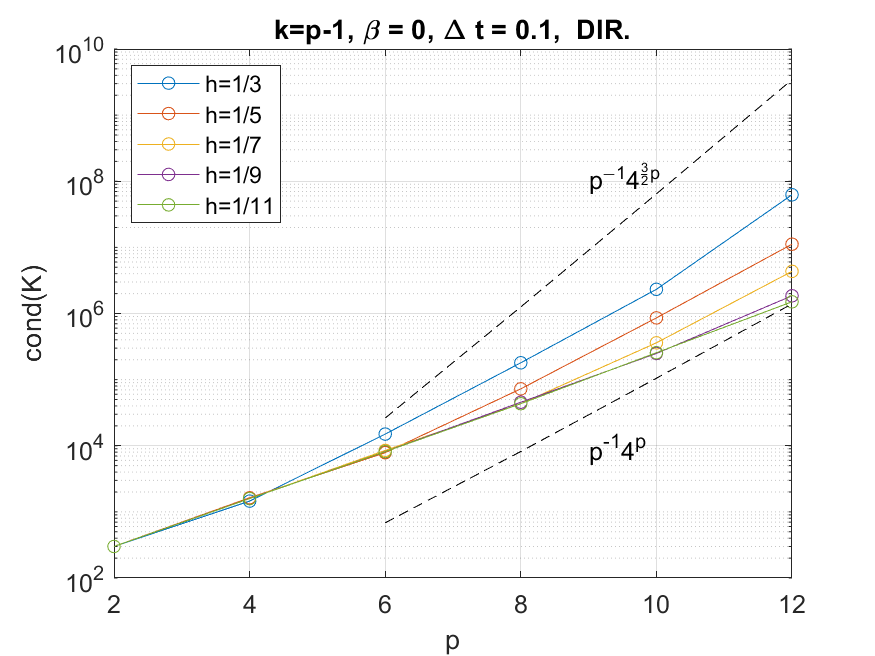} &
\includegraphics[scale=0.485]{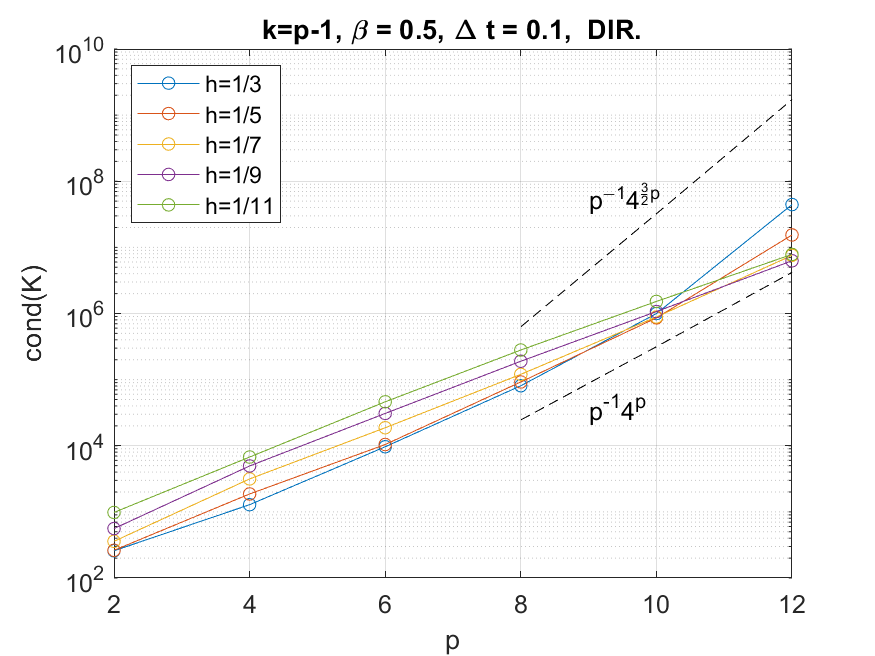} \\
\end{tabular}
}
\caption{Condition number of the stiffness matrix for the acoustic wave problem with Dirichlet boundary conditions,  $\Delta t = 0.1$, $\gamma=0.5$, $\beta = 0$ (explicit Newmark, left), $\beta = 0.5$ (implicit Newmark, right). From the top to the bottom, vs.: (1)  $h$, for $p=2, 4, 6,8, 10$,  $k=1$; (2) $h$, for $p=2,4,6,8, 10$,  $k=p-1$; (3)  $p$, for $h=1/3, 1/5, 1/7, 1/9$,   $k=1$; (4) $p$, for $h=1/3, 1/5, 1/7, 1/9, 1/11$, $k=p-1$.
\label{cond_S_dt01_DIR}}
\end{figure}

\begin{figure} %[!t]
\vspace{-10mm}
\centerline{
\begin{tabular}{cc}
\includegraphics[scale=0.485]{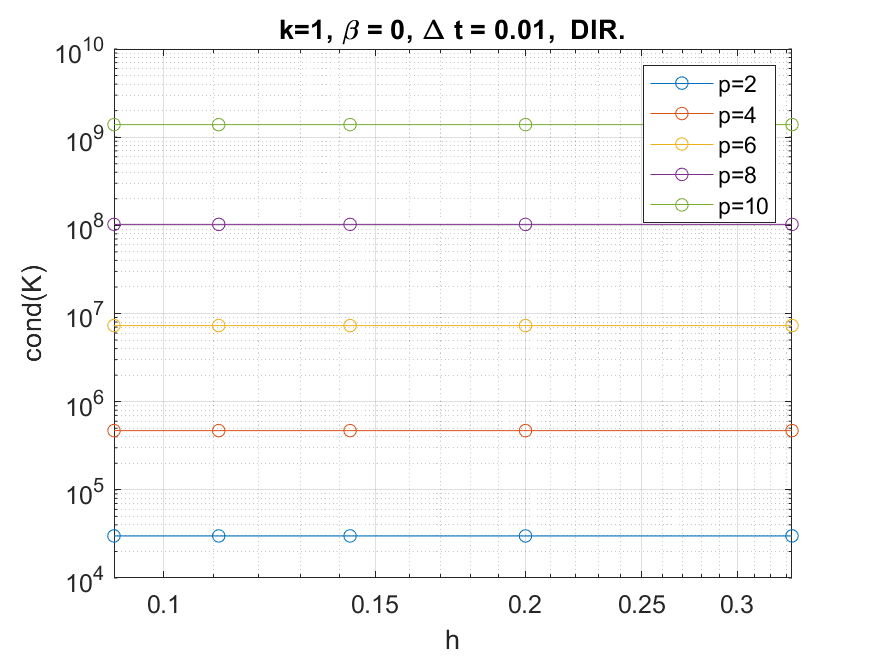} &
\includegraphics[scale=0.485]{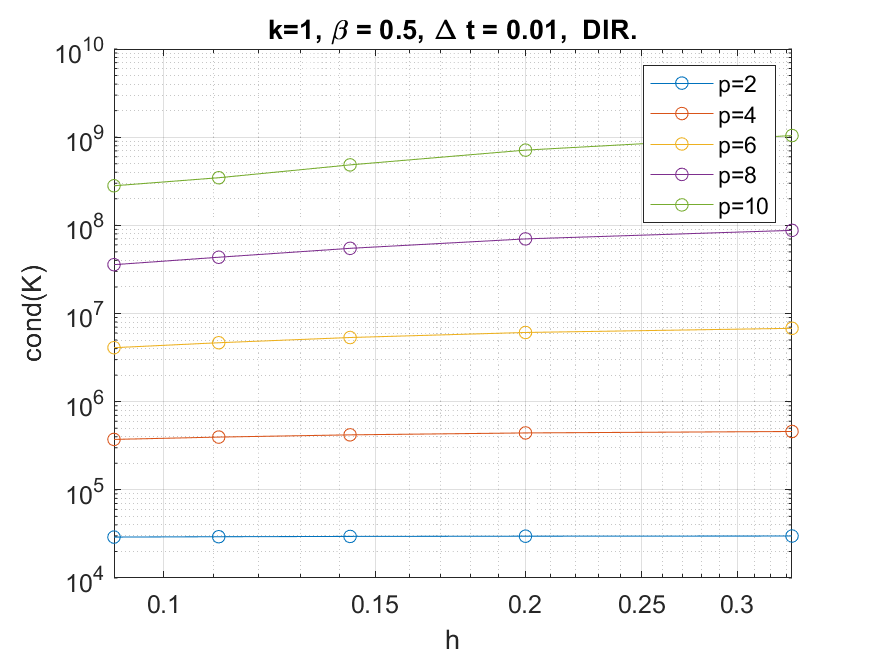} \\
\includegraphics[scale=0.485]{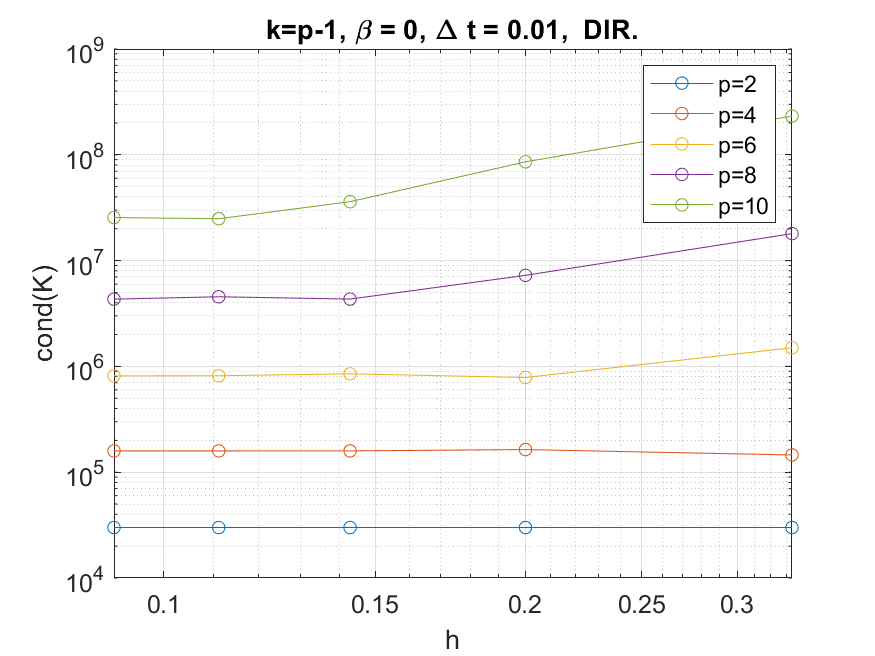} & \includegraphics[scale=0.485]{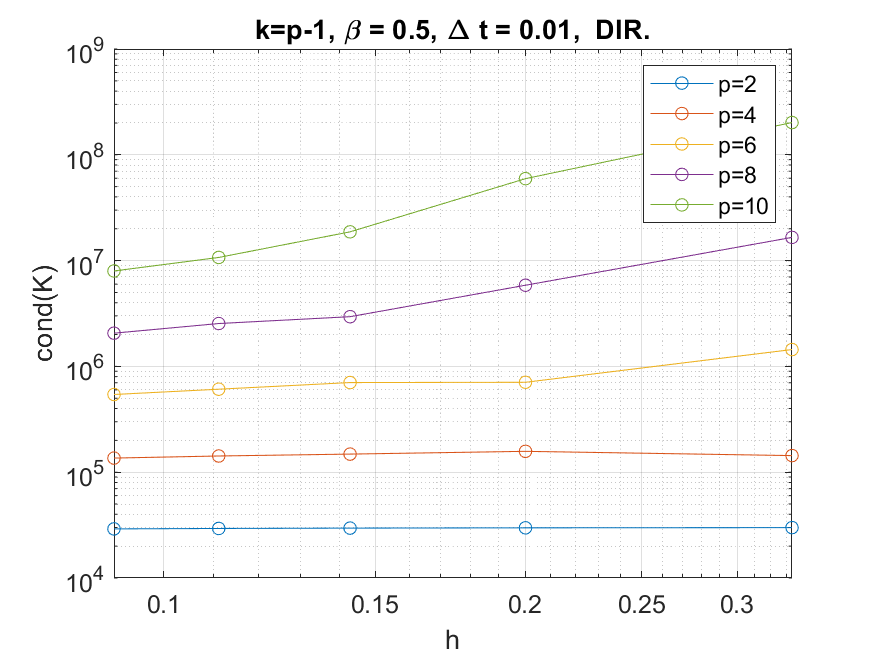} \\
\includegraphics[scale=0.485]{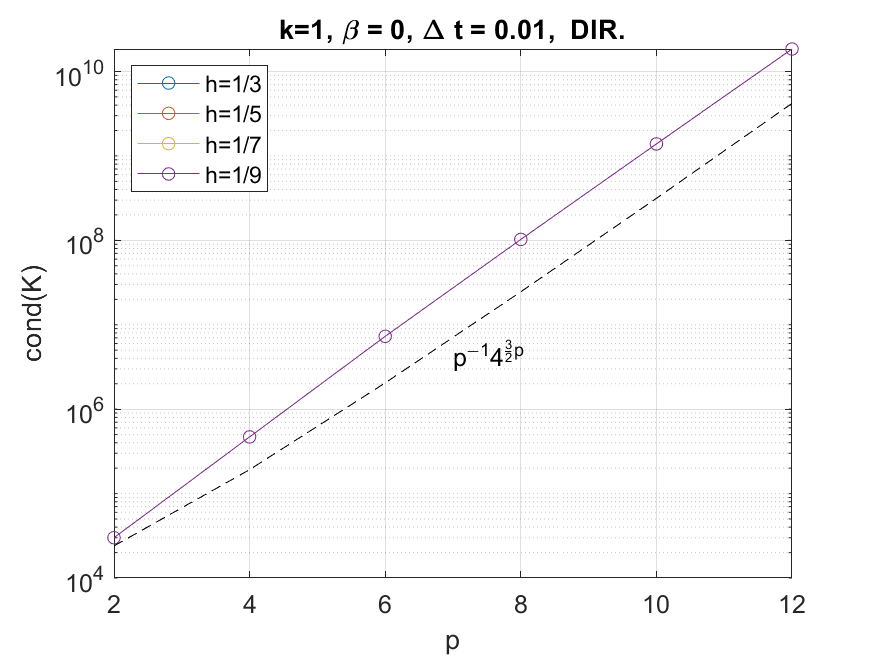} & 
\includegraphics[scale=0.485]{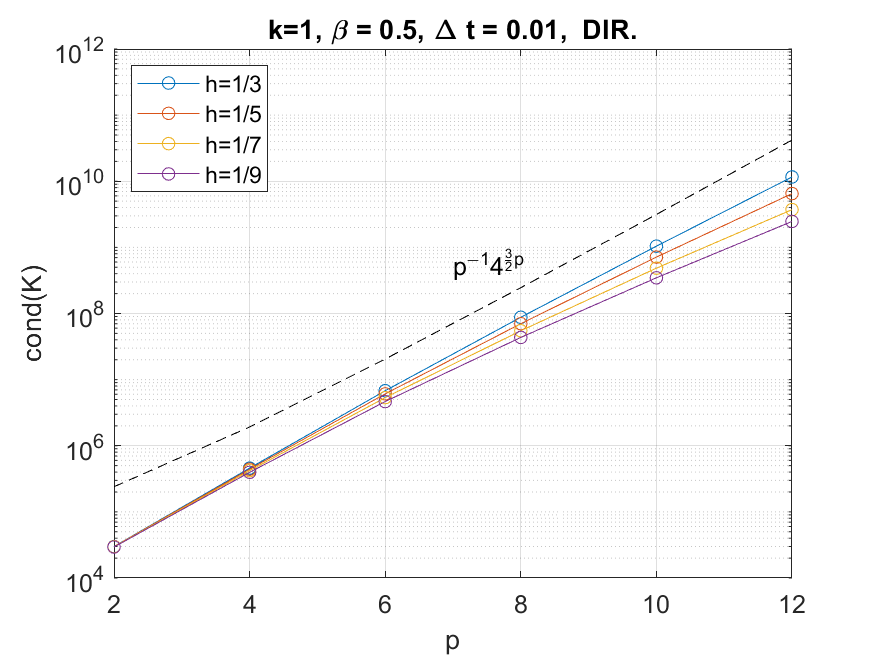} \\
\includegraphics[scale=0.485]{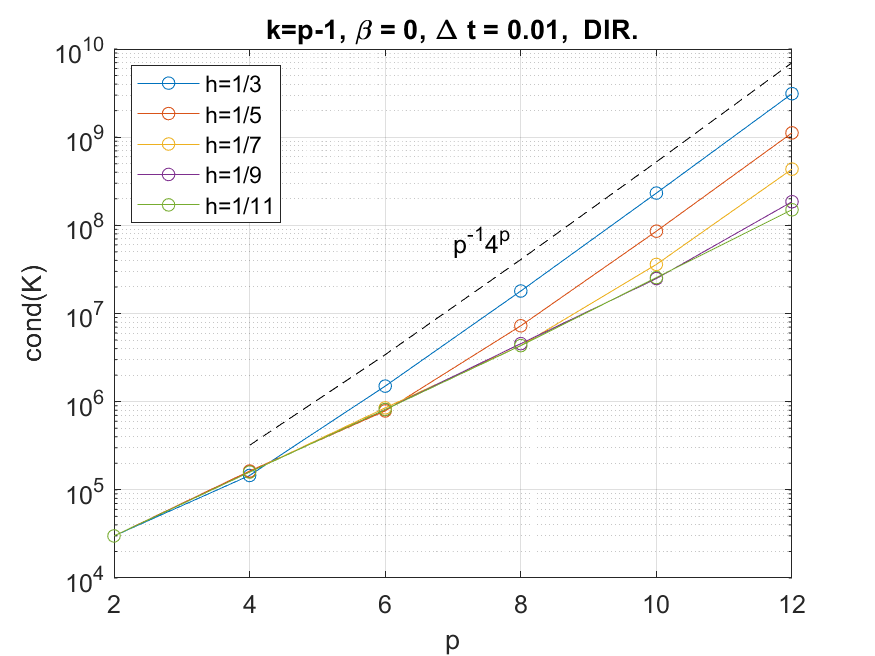} &
\includegraphics[scale=0.485]{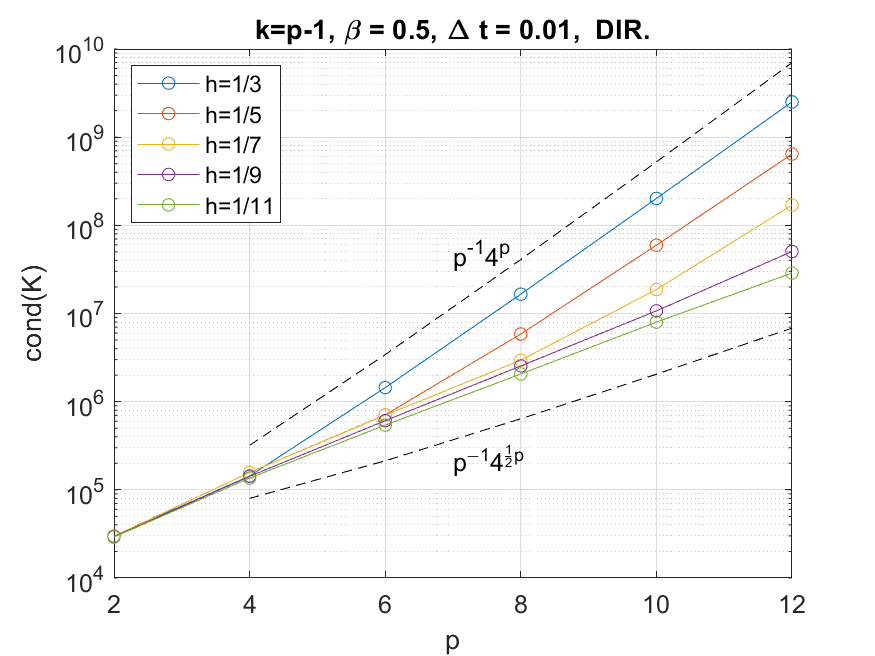} \\
\end{tabular}
}
\caption{Condition number of the stiffness matrix for the acoustic wave problem with Dirichlet boundary conditions,  $\Delta t = 0.01$, $\gamma=0.5$, $\beta = 0$ (explicit Newmark, left), $\beta = 0.5$ (implicit Newmark, right). From the top to the bottom, vs.: (1)  $h$, for $p=2,4, 6,8, 10$,  $k=1$; (2) $h$, for $p=2, 4, 6,8, 10$,  $k=p-1$; (3)  $p$, for $h=1/3, 1/5, 1/7, 1/9$,   $k=1$; (4) $p$, for $h=1/3, 1/5, 1/7, 1/9,1/11$, $k=p-1$.
\label{cond_S_dt001_DIR}}
\end{figure}

In Fig.  \ref{cond_S_dt01_DIR},  we report the condition numbers {\textsf{cond}}($\cal K$) of the stiffness matrix for the acoustic wave problem with Dirichlet boundary conditions, for $\Delta t=0.1$,   $\beta=0$ (explicit Newmark scheme, left) or  $\beta=0.5$ (implicit Newmark scheme, right), versus, from the top to the bottom:  (1)  the mesh size $h$, with five different values of degree $p$ and minimal regularity $k=1$, (2)  the mesh size $h$, with five different values of degree $p$ and maximal regularity $k=p-1$, (3) the degree $p$, with four different values of mesh size $h$ and minimal regularity $k=1$, (4) the degree $p$, with five different values of mesh size $h$ and maximal regularity $k=p-1$.
In Fig. \ref{cond_S_dt001_DIR}, we consider  the same tests as in Fig. \ref{cond_S_dt01_DIR} but with smaller time step $\Delta t=0.01$. 
The numerical results show that if $p$ is fixed, the condition numbers  {\textsf{cond}}($\cal K$)  are almost always independent  of $h$ in both Figs.  \ref{cond_S_dt01_DIR} and  \ref{cond_S_dt001_DIR}, except for the implicit scheme (right) in Fig. \ref{cond_S_dt01_DIR} where they seem to grow as $h^{-2}$,  
%with $\Delta t=0.1$, 
according to estimates (\ref{Kk0})-(\ref{Kkmp1}). For the $p$- refinement with fixed $h$,   it seems that the numerical results  are better than these estimates. Indeed, the condition numbers  {\textsf{cond}}($\cal K$) grow as $ p^{-1}4^{\frac32 p}$ in the case of minimal regularity $k=1$, whereas for maximal regularity $k=p-1$ the growth ranges between $ p^{-1}4^p$  and $p^{-1}4^{\frac32 p}$  for $\Delta t=0.1$,  between $p^{-1}4^{\frac12 p}$  and $ p^{-1}4^p$ in the implicit case with $\Delta t=0.01$, while the growth seems to be $ p^{-1}4^p$ in the explicit case with $\Delta t=0.01$.  

\begin{figure} %[!t]
\vspace{-10mm}
\centerline{
\begin{tabular}{cc}
\includegraphics[scale=0.485]{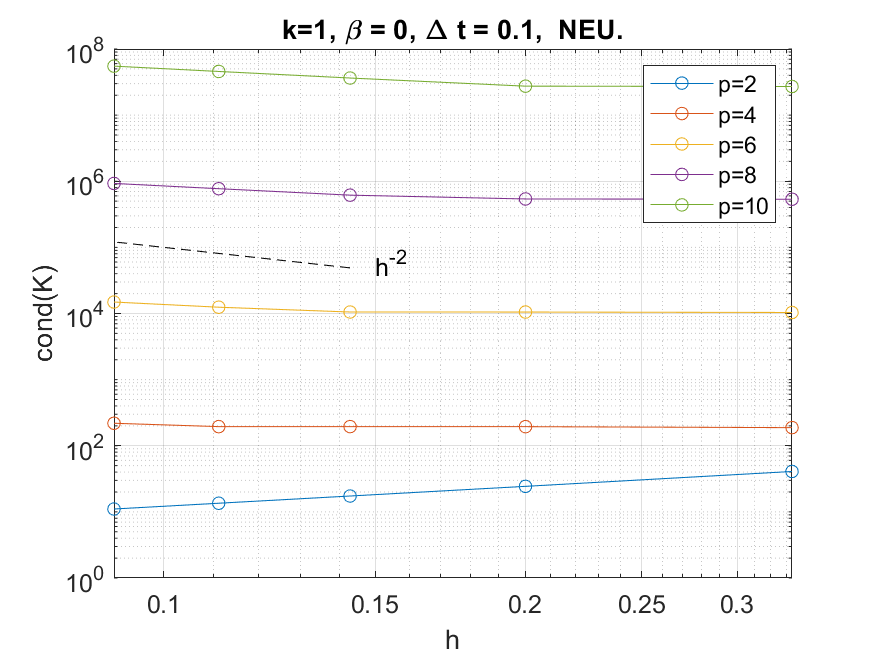} &
\includegraphics[scale=0.485]{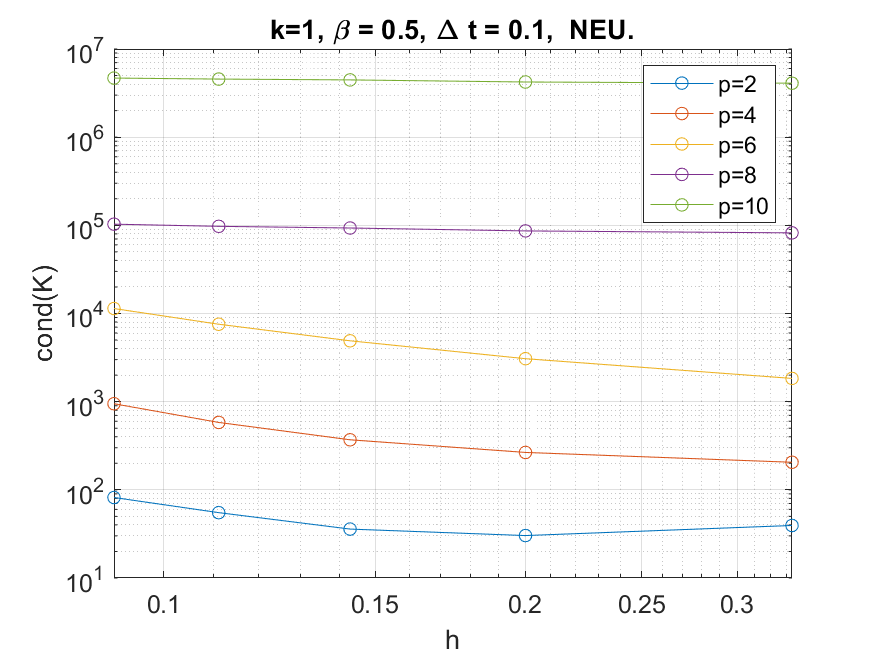} \\
\includegraphics[scale=0.485]{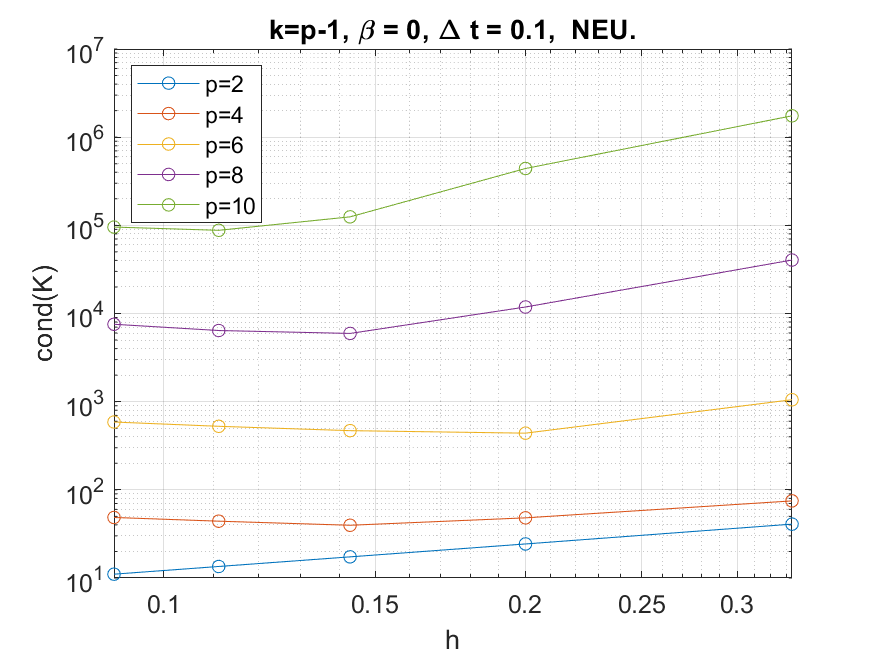} & \includegraphics[scale=0.485]{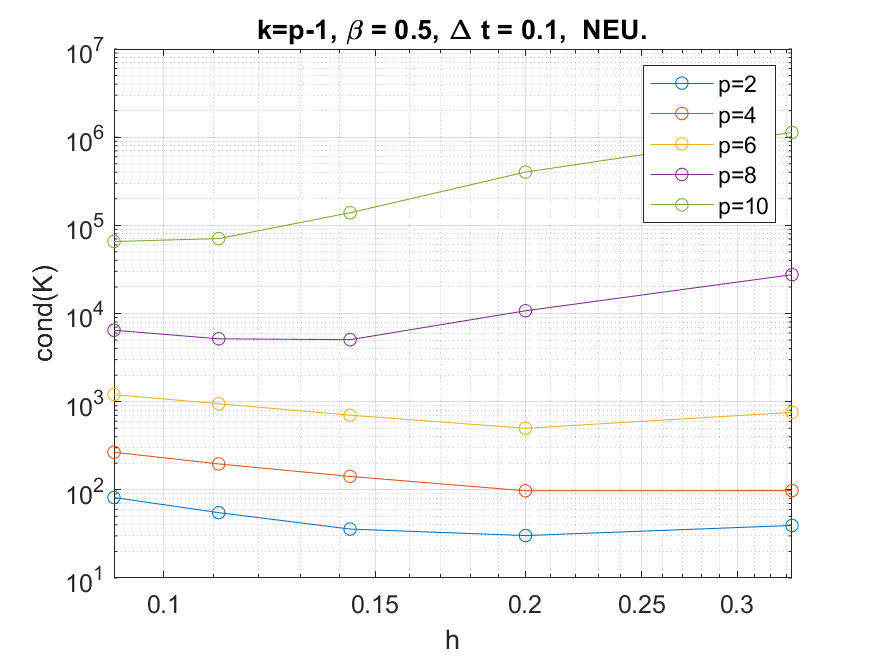} \\
\includegraphics[scale=0.485]{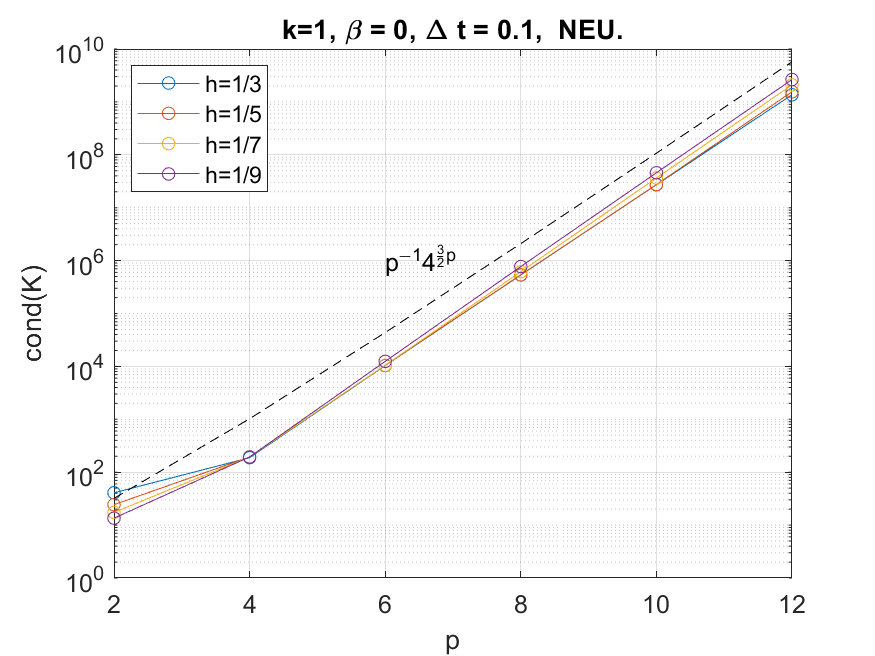} & 
\includegraphics[scale=0.485]{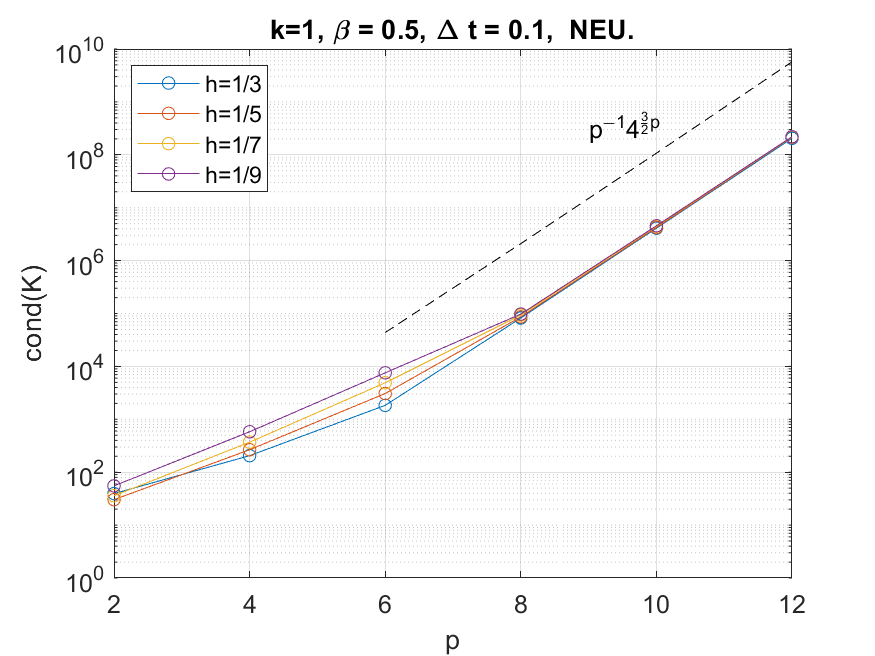} \\
\includegraphics[scale=0.485]{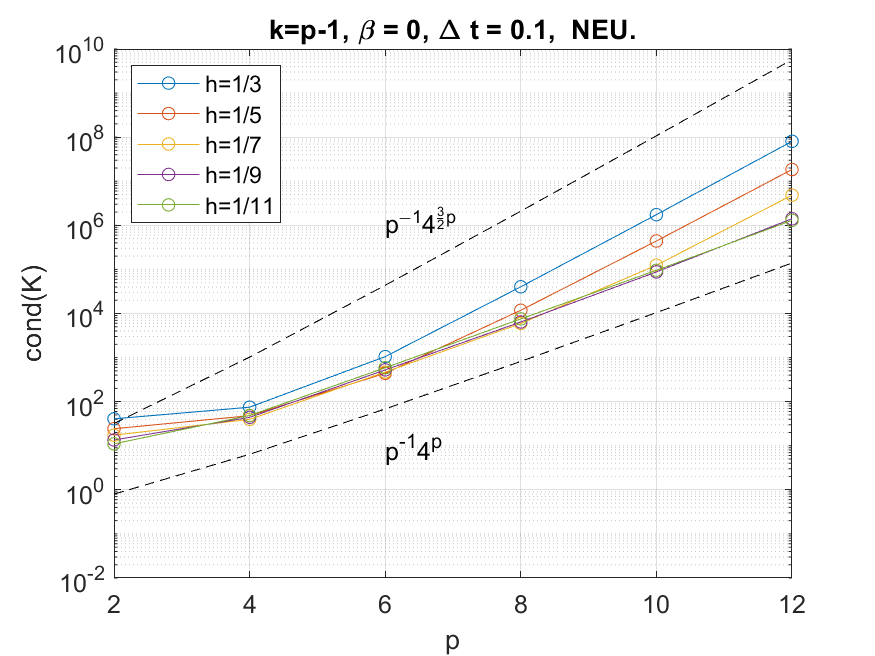} &
\includegraphics[scale=0.485]{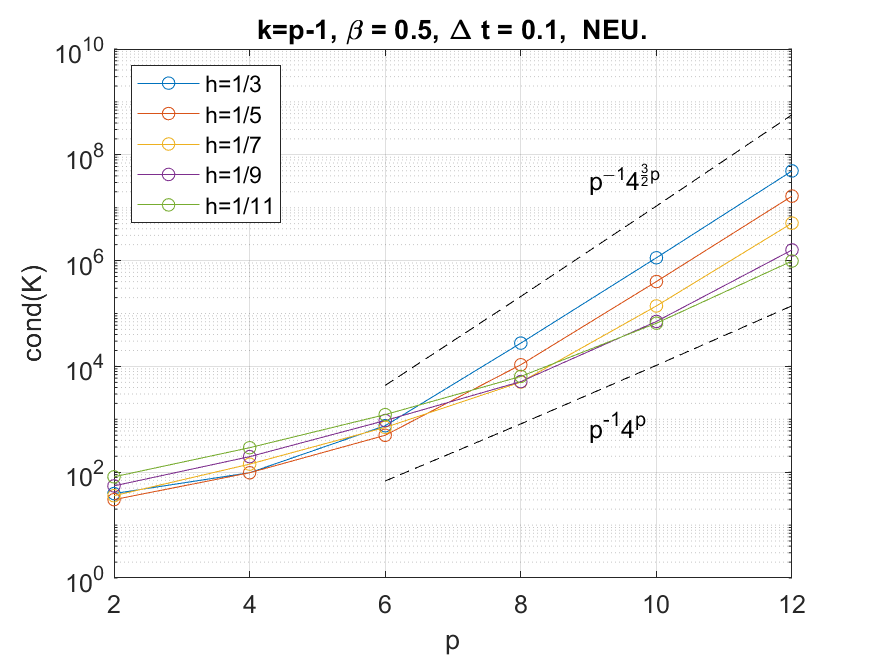} \\
\end{tabular}
}
\caption{Condition number of the stiffness matrix for the acoustic wave problem with Neumann boundary conditions,  $\Delta t = 0.1$, $\gamma=0.5$, $\beta = 0$ (explicit Newmark, left), $\beta = 0.5$ (implicit Newmark, right). From the top to the bottom, vs.: (1)  $h$, for $p=2, 4, 6,8, 10$,  $k=1$; (2) $h$, for $p=2, 4, 6,8, 10$,  $k=p-1$; (3)  $p$, for $h=1/3, 1/5, 1/7,1/9$,   $k=1$; (4) $p$, for $h=1/3, 1/5, 1/7, 1/9,1/11$, $k=p-1$. 
\label{cond_S_dt01_NEU}}
\end{figure}

\begin{figure} %[!t]
\vspace{-10mm}
\centerline{
\begin{tabular}{cc}
\includegraphics[scale=0.485]{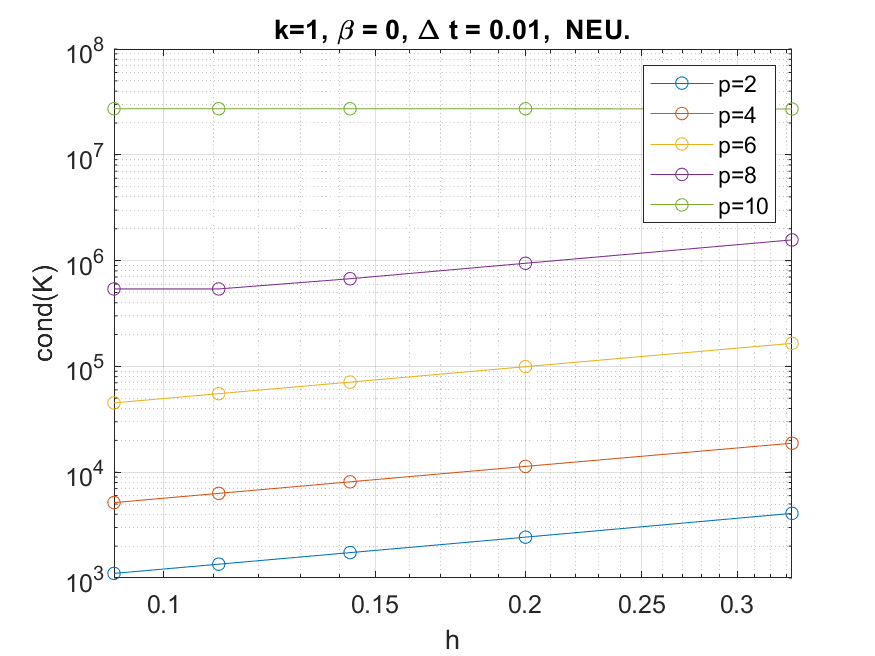} &
\includegraphics[scale=0.485]{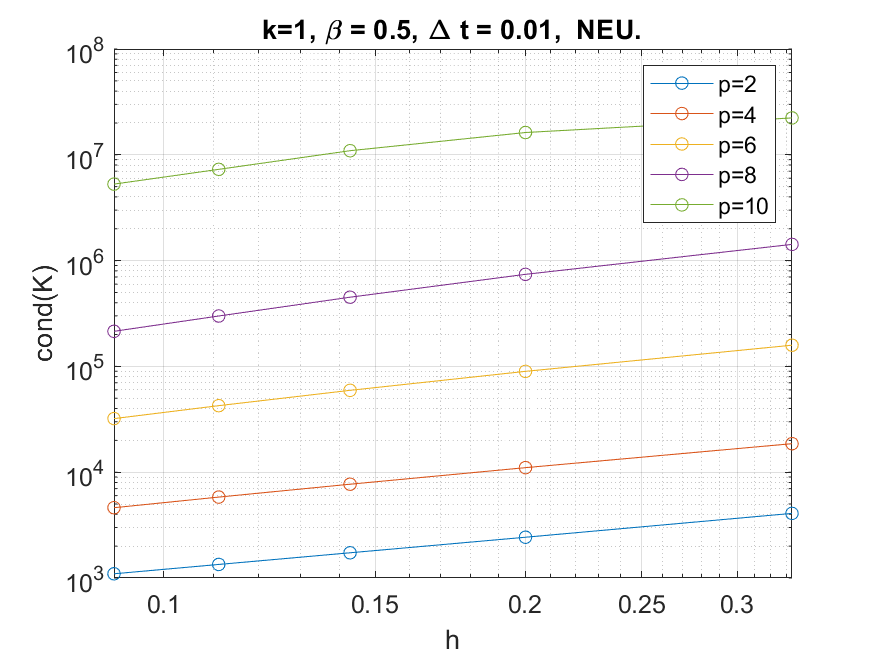} \\
\includegraphics[scale=0.485]{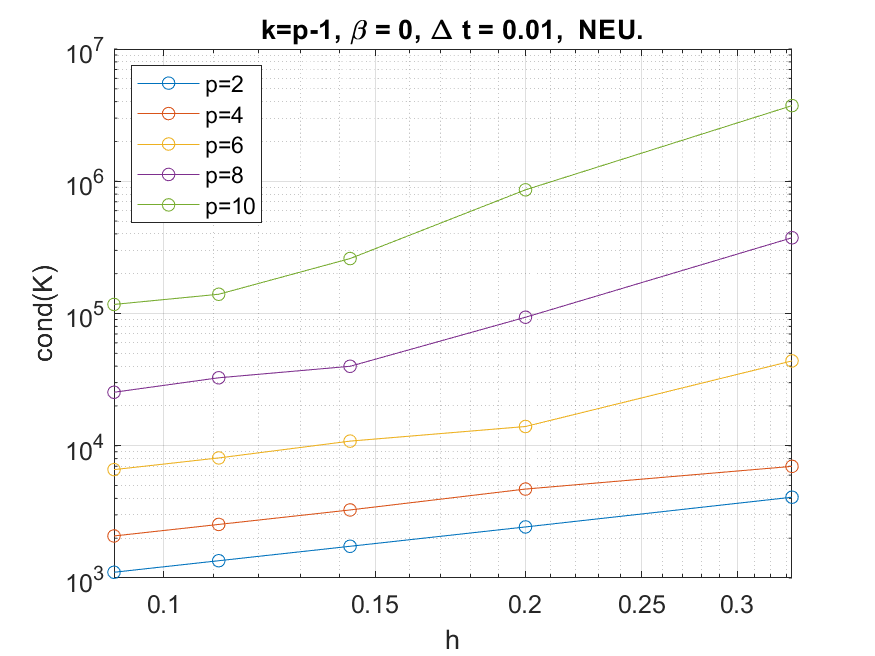} & \includegraphics[scale=0.485]{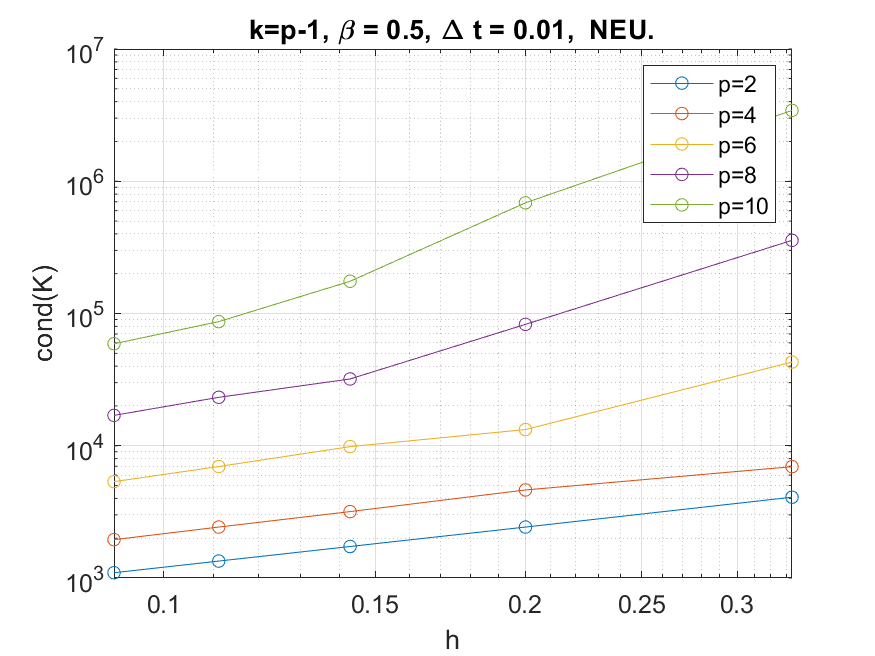} \\
\includegraphics[scale=0.485]{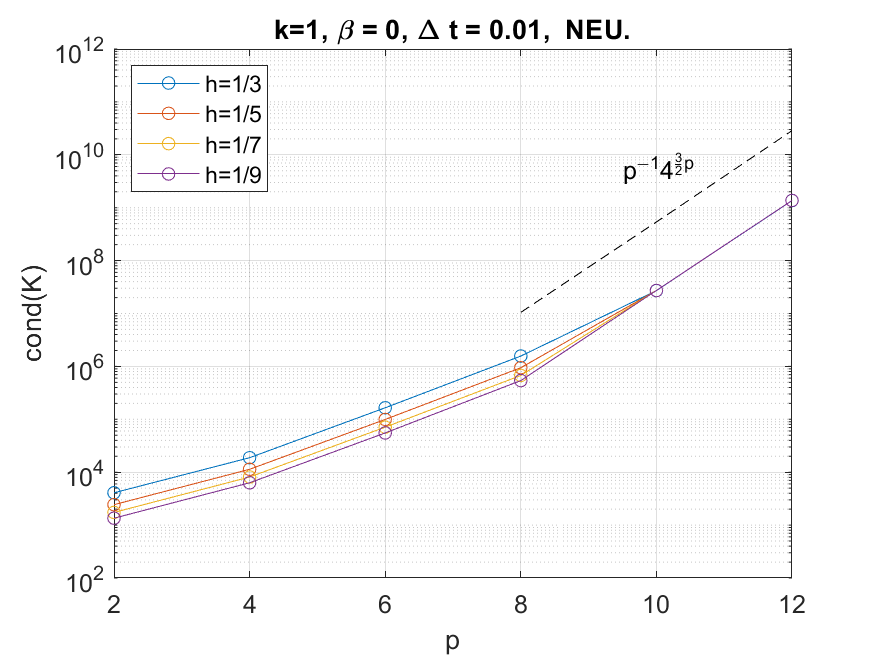} & 
\includegraphics[scale=0.485]{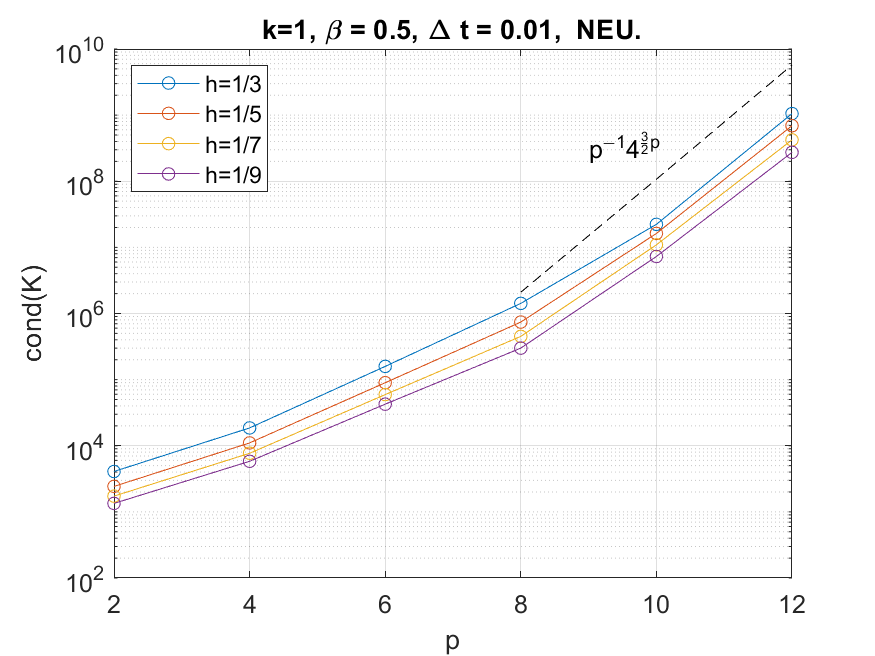} \\
\includegraphics[scale=0.485]{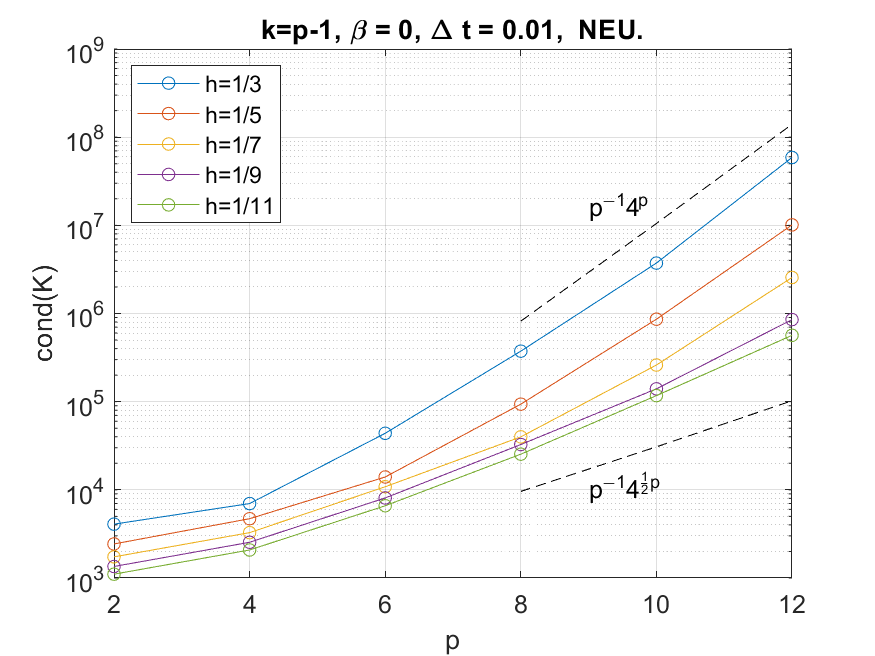} &
\includegraphics[scale=0.485]{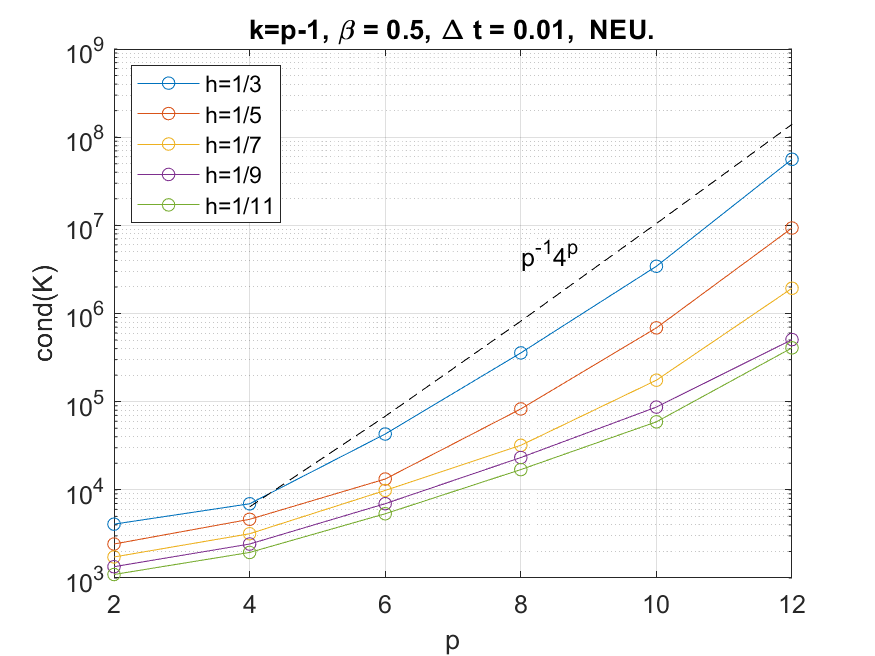} \\
\end{tabular}
}
\caption{Condition number of the stiffness matrix for the acoustic wave problem with Neumann boundary conditions,  $\Delta t = 0.01$, $\gamma=0.5$, $\beta = 0$ (explicit Newmark, left), $\beta = 0.5$ (implicit Newmark, right). From the top to the bottom, vs.: (1)  $h$, for $p=2, 4, 6,8, 10$,  $k=1$; (2) $h$, for $p=2, 4,6,8, 10$,  $k=p-1$; (3)  $p$, for $h=1/3, 1/5, 1/7, 1/9$,   $k=1$; (4) $p$, for $h=1/3, 1/5, 1/7, 1/9,1/11$, $k=p-1$. 
\label{cond_S_dt001_NEU}}
\end{figure}

Figs. \ref{cond_S_dt01_NEU} and \ref{cond_S_dt001_NEU} reports the condition numbers {\textsf{cond}}($\cal K$) of the stiffness matrix for the acoustic wave problem with Neumann boundary conditions, using the same setting of Figs. \ref{cond_S_dt01_DIR} and \ref{cond_S_dt001_DIR}, respectively.
We note  that if $p$ is fixed, the condition numbers  {\textsf{cond}}($\cal K$)  are almost always independent  of $h$, in particular they do not increase when $h$ decreases, except in some cases with high $p$ for the explicit case with $\Delta t=0.1$, where the condition numbers scale at most as $h^{-2}$,  as predicted by estimates (\ref{Kk0})-(\ref{Kkmp1}). For the $p$- refinement, with fixed $h$, similarly to the Dirichlet case, the results are again better than the estimates (\ref{Kk0})-(\ref{Kkmp1}). In fact, the condition numbers  {\textsf{cond}}($\cal K$) scale as $ p^{-1}4^{\frac32 p}$ in the case of minimal regularity $k=1$, whereas for maximal regularity $k=p-1$ they  range between $ p^{-1}4^p$  and $ p^{-1}4^{\frac32 p}$  for $\Delta t=0.1$, range  between $ p^{-1}4^{\frac12 p}$  and $ p^{-1}4^p$ in the explicit case with $\Delta t=0.01$, and scale as $ p^{-1}4^p$ in the implicit case with $\Delta t=0.01$.

\begin{figure} %[!t]
\vspace{4mm}
\centerline{
\begin{tabular}{ccc}
\includegraphics[scale=0.40]{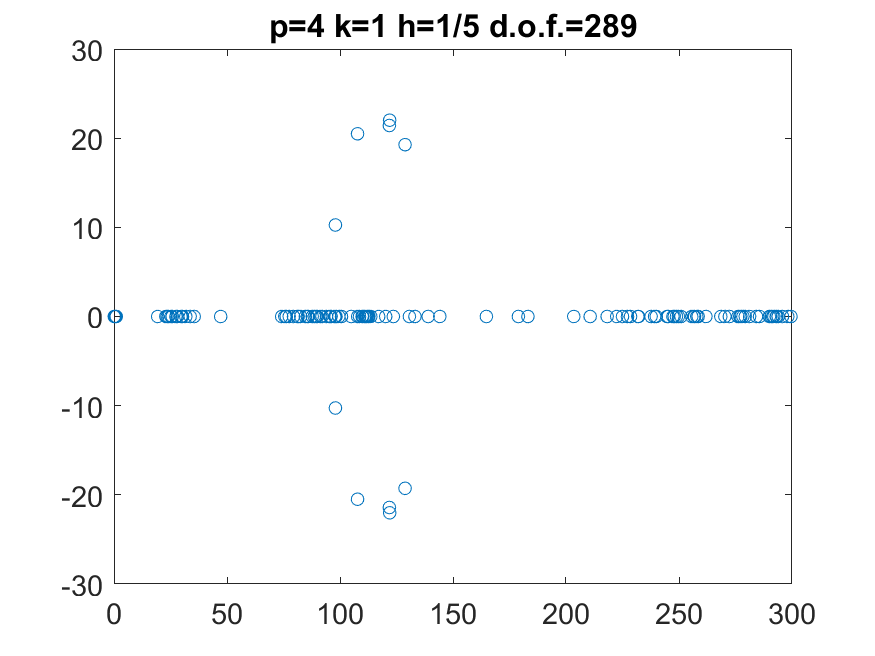} &
\hspace{-5mm}\includegraphics[scale=0.40]{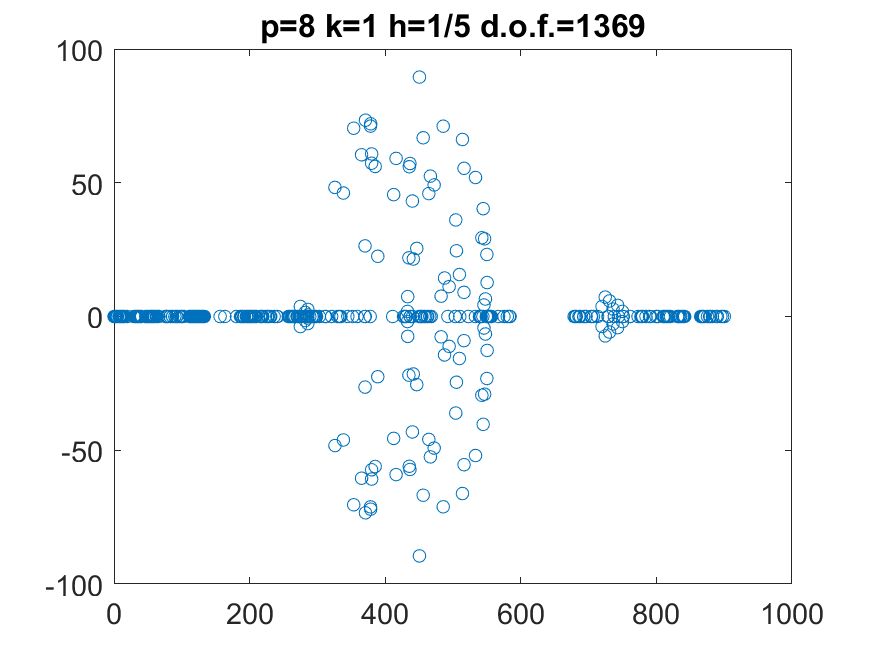} &
\hspace{-5mm}\includegraphics[scale=0.40]{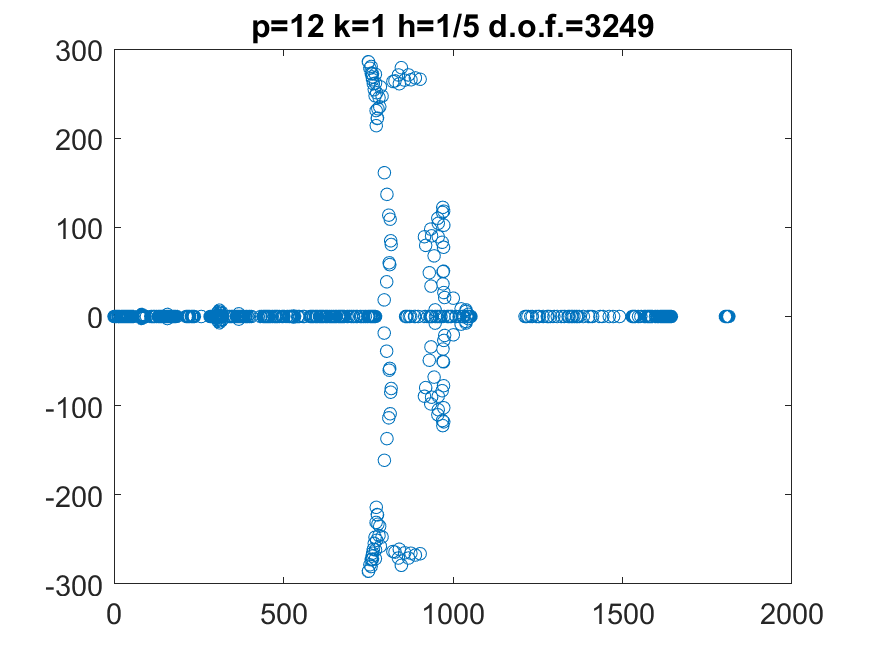} \\
\includegraphics[scale=0.40]{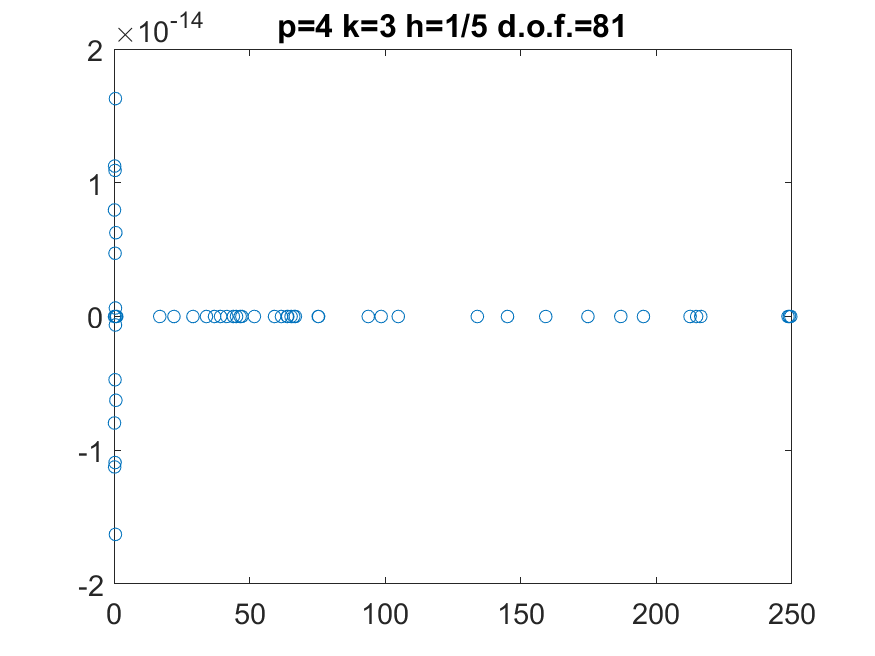} & 
\hspace{-5mm}\includegraphics[scale=0.40]{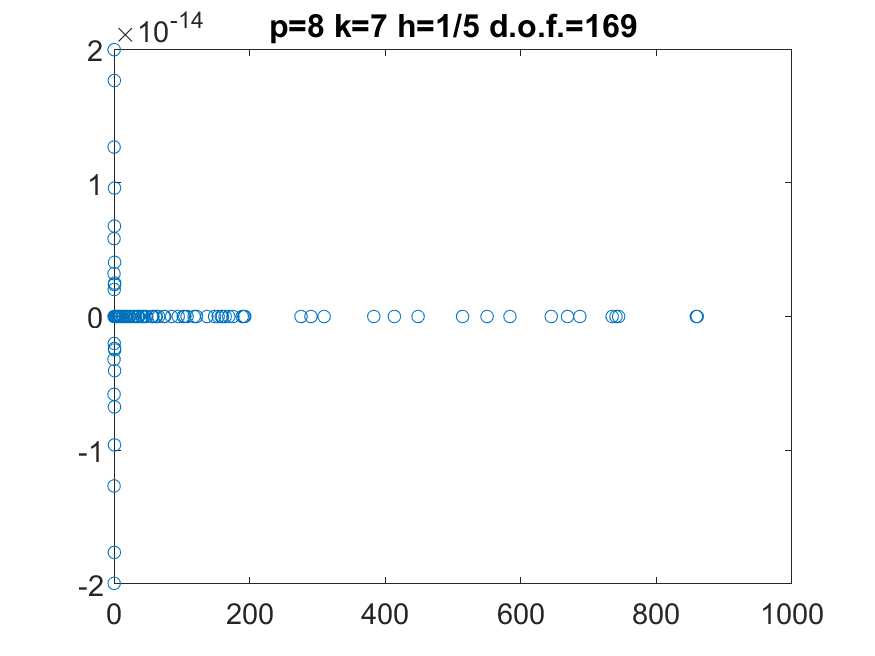} & 
\hspace{-5mm}\includegraphics[scale=0.40]{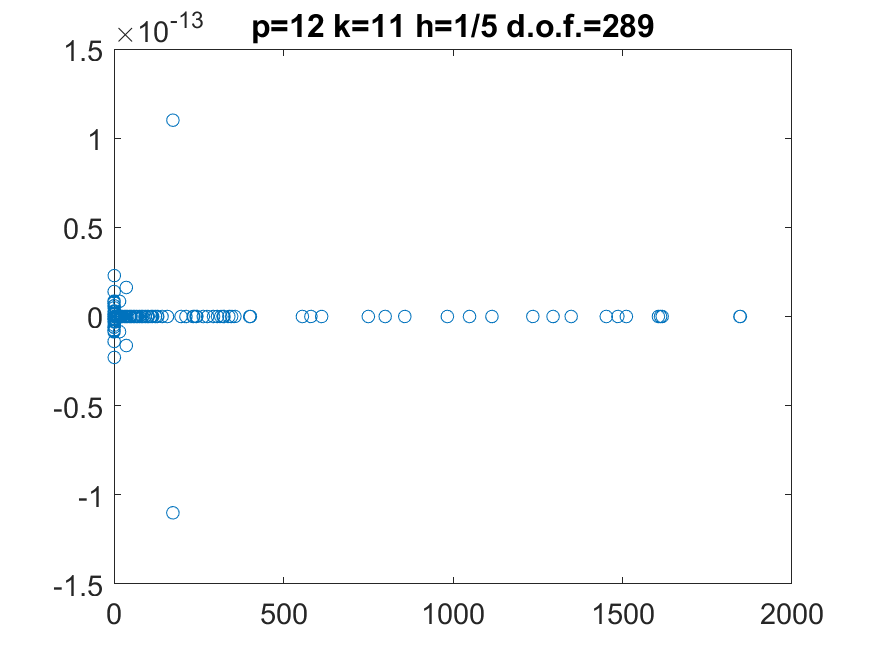}  \\
\end{tabular}
}
\vspace{-2mm}
\caption{Stiffness matrix eigenvalue distribution in the complex plane for the acoustic wave problem with Dirichlet boundary conditions: $p = 4$ (left),
 $p=8$ (center), $p=12$ (right), with $k=1$ (top) or $k=p-1$ (bottom), fixed $h=1/5$, $\Delta t=0.01$, $\gamma=0.5$, $\beta=0.5$.
\label{eig_S_h5_VSp_Laplacian_DIR}}
\end{figure}

\begin{figure} %[!t]
\vspace{4mm}
\centerline{
\begin{tabular}{ccc}
\includegraphics[scale=0.40]{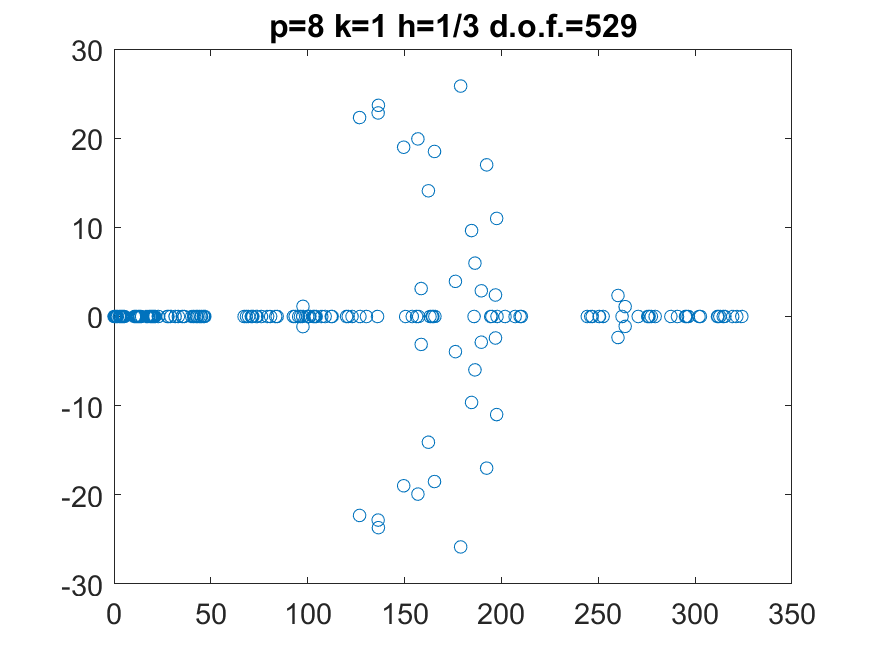} &
\hspace{-5mm}\includegraphics[scale=0.40]{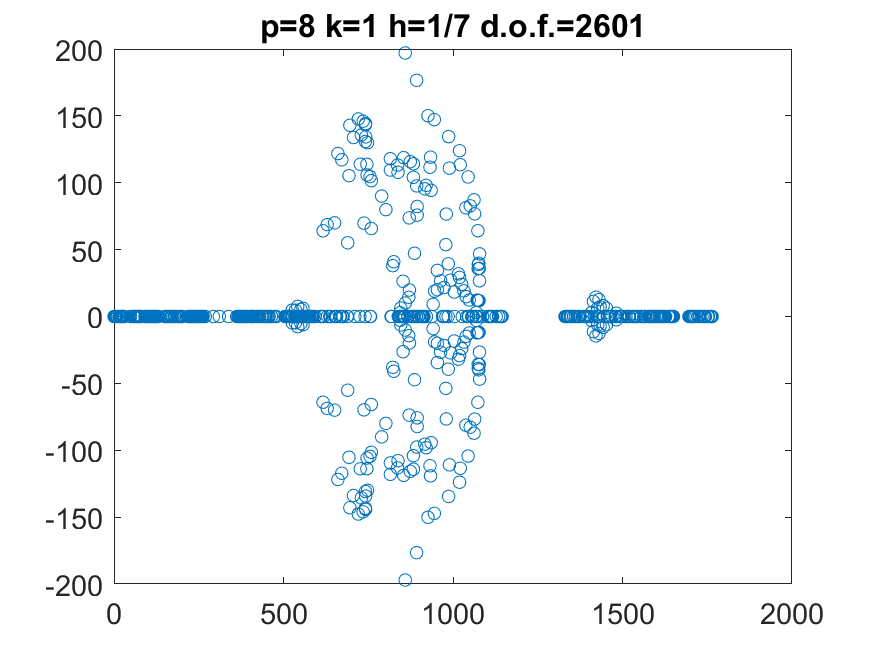} &
\hspace{-5mm}\includegraphics[scale=0.40]{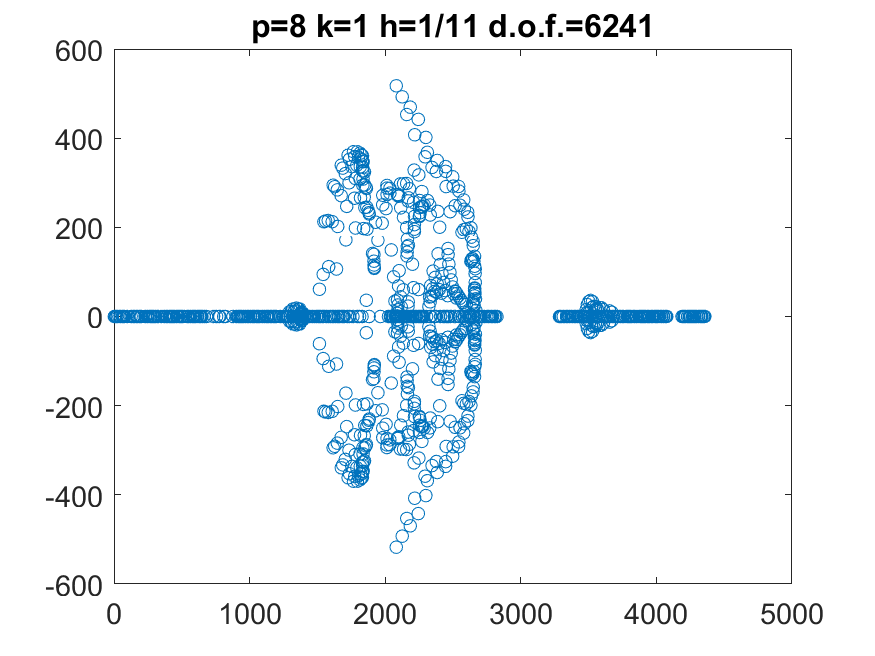} \\
\includegraphics[scale=0.40]{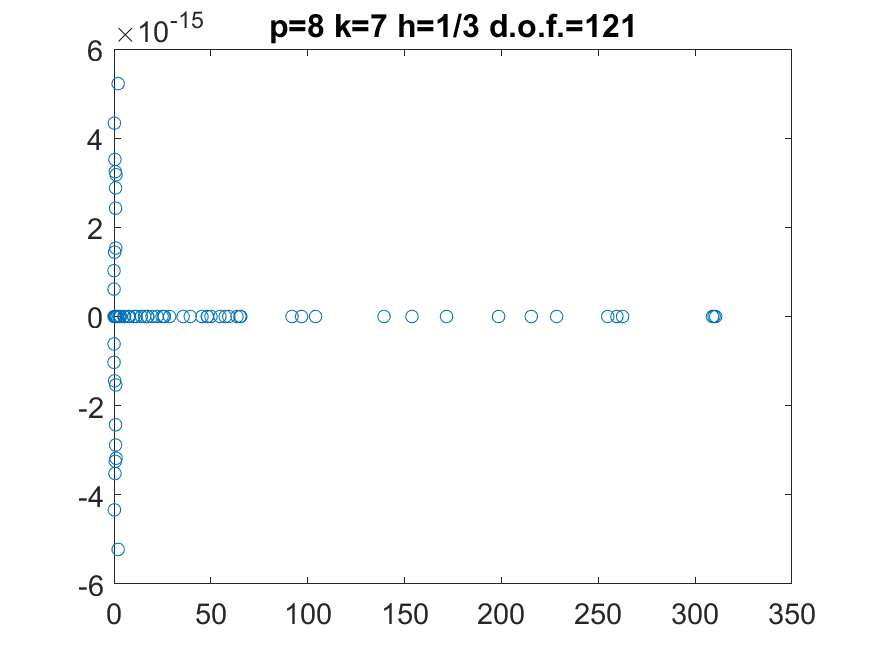} & 
\hspace{-5mm}\includegraphics[scale=0.40]{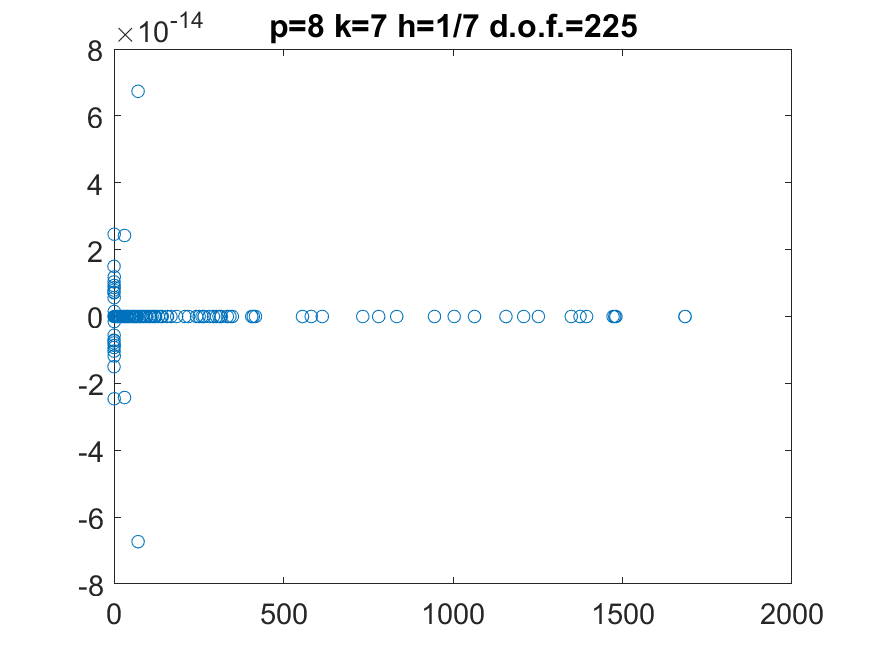} & 
\hspace{-5mm}\includegraphics[scale=0.40]{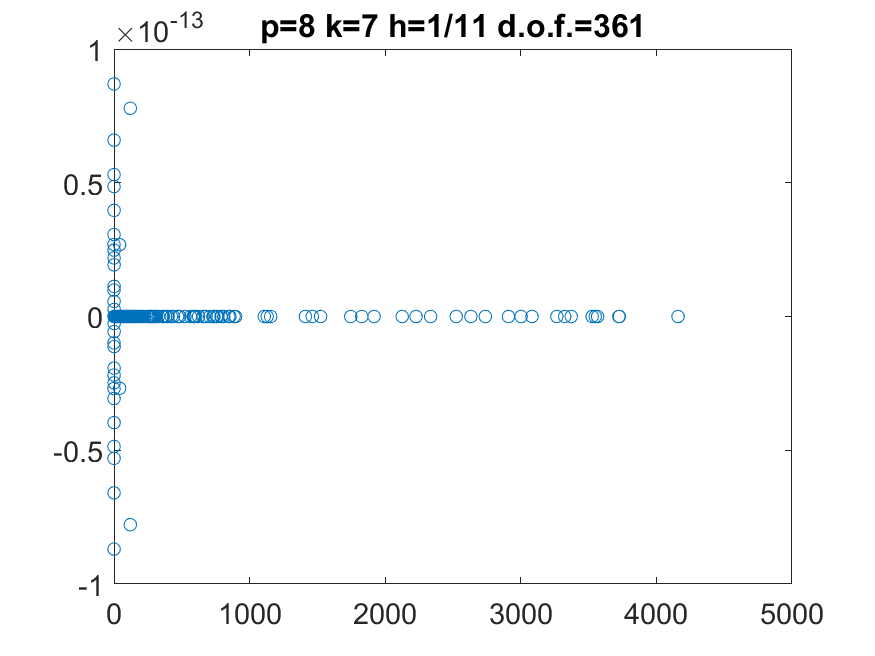}  \\
\end{tabular}
}
\vspace{-2mm}
\caption{Stiffness matrix eigenvalue distribution in the complex plane for the acoustic wave problem with Dirichlet boundary conditions: $h=1/3$ (left),
 $1/h=7$ (center), $1/h=11$ (right), with $k=1$ (top) or $k=p-1$ (bottom), fixed $p=8$, $\Delta t=0.01$, $\gamma=0.5$, $\beta=0.5$.
\label{eig_S_p8_VSh_Laplacian_DIR}}
\end{figure}

Fig. \ref{eig_S_h5_VSp_Laplacian_DIR} reports the stiffness matrix $\cal K$ eigenvalue distribution in the complex plane  for the acoustic wave problem with Dirichlet boundary conditions, for $\Delta t=0.01$, $\gamma=0.5$, $\beta=0.5$. We consider three different values of degree $p=4$ (left), $p=8$ (center), $p=12$ (right), minimal regularity $k=1$ (top) or maximal regularity $k=p-1$ (bottom), for fixed $h=1/5$. 
In Fig. \ref{eig_S_p8_VSh_Laplacian_DIR}, we use the same setting to report the analogous eigenvalue distributions 
%in complex plane of stiffness matrix $\cal K$ for the acoustic wave problem with Dirichlet boundary conditions, and $\Delta $, $\gamma$, $\beta$ as in Fig. \ref{eig_S_h5_VSp_Laplacian_DIR}.
for three different values of  mesh size $h=1/3$ (left), $h=1/7$ (center), $h=1/11$ (right), minimal regularity $k=1$ (top) or maximal regularity $k=p-1$ (bottom), for fixed $p=8$.
%For fixed $p$ and $h$, 
The eigenvalues real parts belong to an interval $[0, r]$ where $r$ increases with both $1/h$ and $p$, but for fixed $p$ and $h$, $r$ is approximately the same for both minimal regularity $k=1$  and maximal regularity $k=p-1$. On the other hand, the eigenvalue imaginary  parts are negligible for maximal regularity $k=p-1$, while for minimal regularity $k=1$ the eigenvalues imaginary parts belong to an interval $[-s, s]$ where $s$ increases with both $1/h$ and $p$.

%(SOLO QUELLI ADDENSATI AL CENTRO SONO COMPLESSI) of eigenvalues having   real parts far from zero and far from their maximum   increase when $p$ increases, fixed $h$, or when $h$ decreases, fixed $p$, while they are almost negligible for maximal regularity $k=p-1$, independently of  $p$ and $h$.

\begin{figure} %[!t]
\vspace{4mm}
\centerline{
\begin{tabular}{ccc}
\includegraphics[scale=0.40]{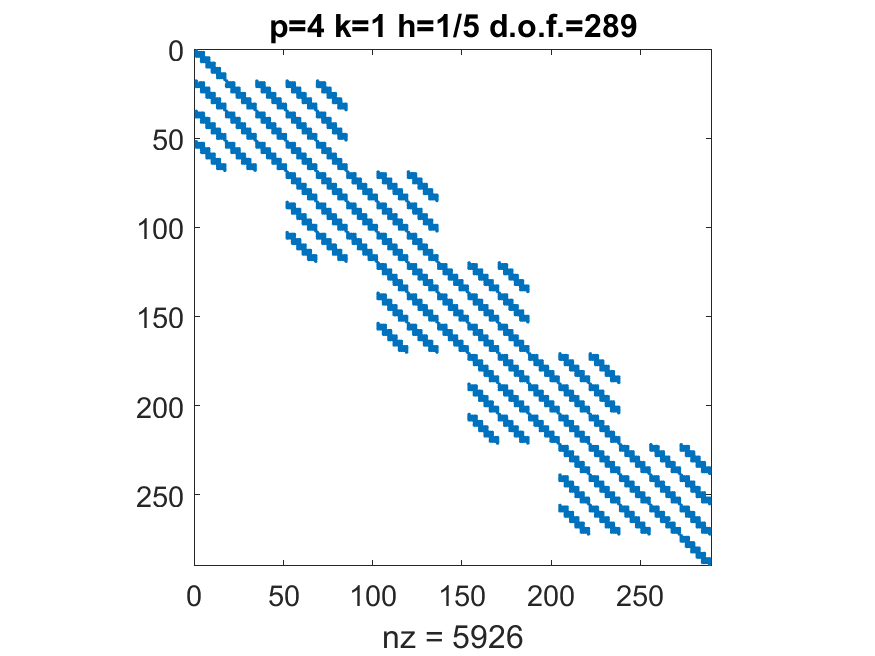} &
\hspace{-5mm}\includegraphics[scale=0.40]{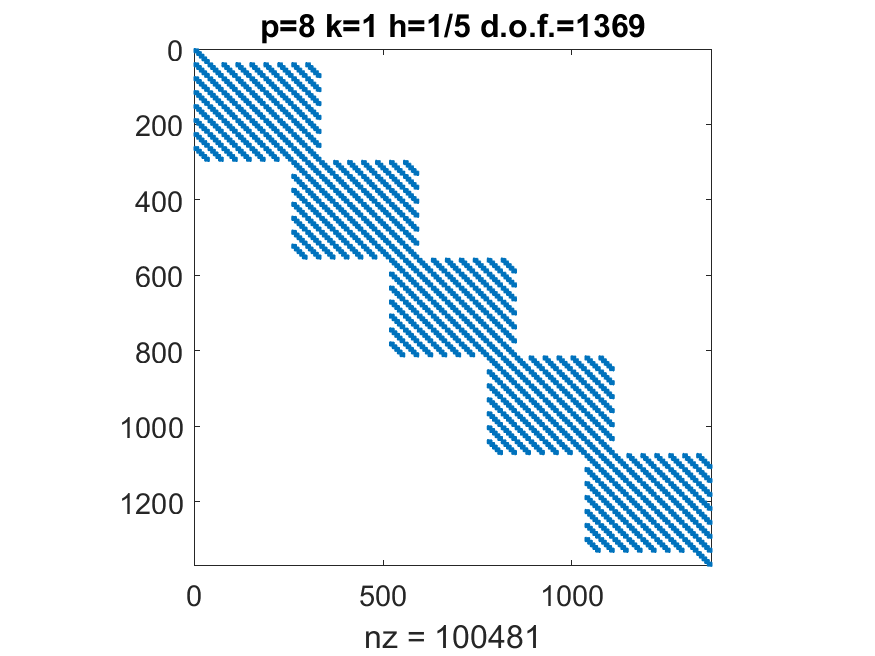} &
\hspace{-5mm}\includegraphics[scale=0.40]{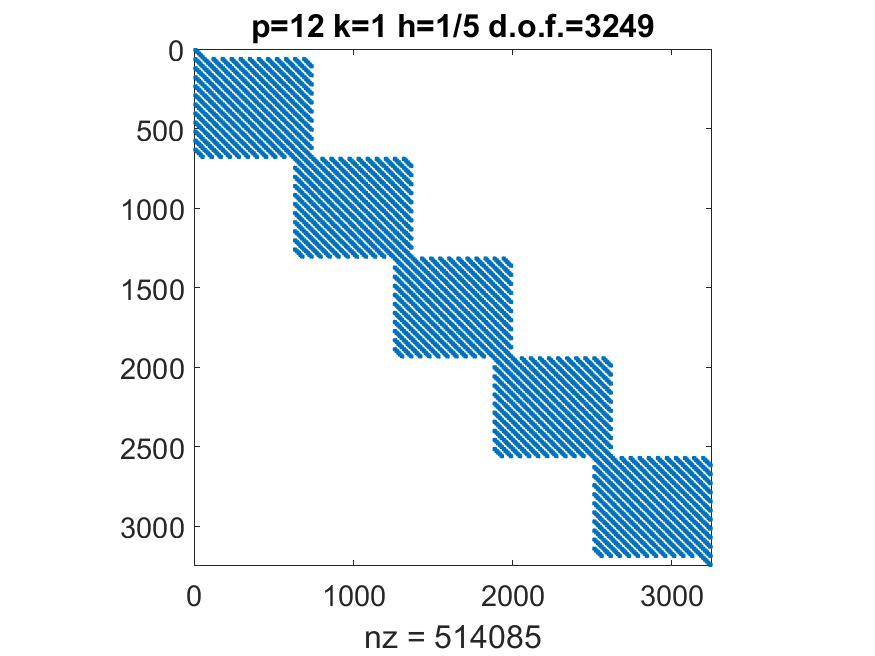} \\
\includegraphics[scale=0.40]{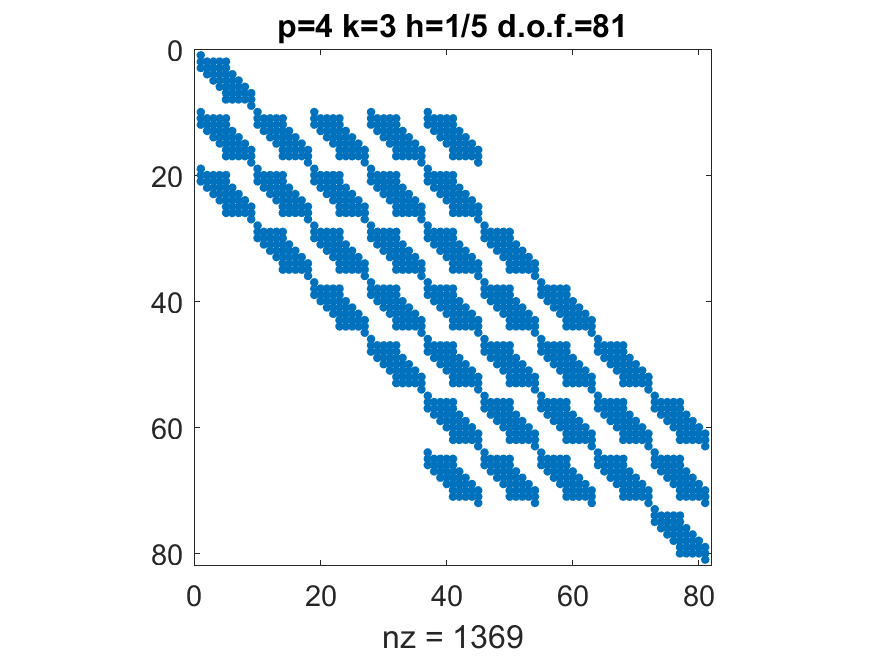} & 
\hspace{-5mm}\includegraphics[scale=0.40]{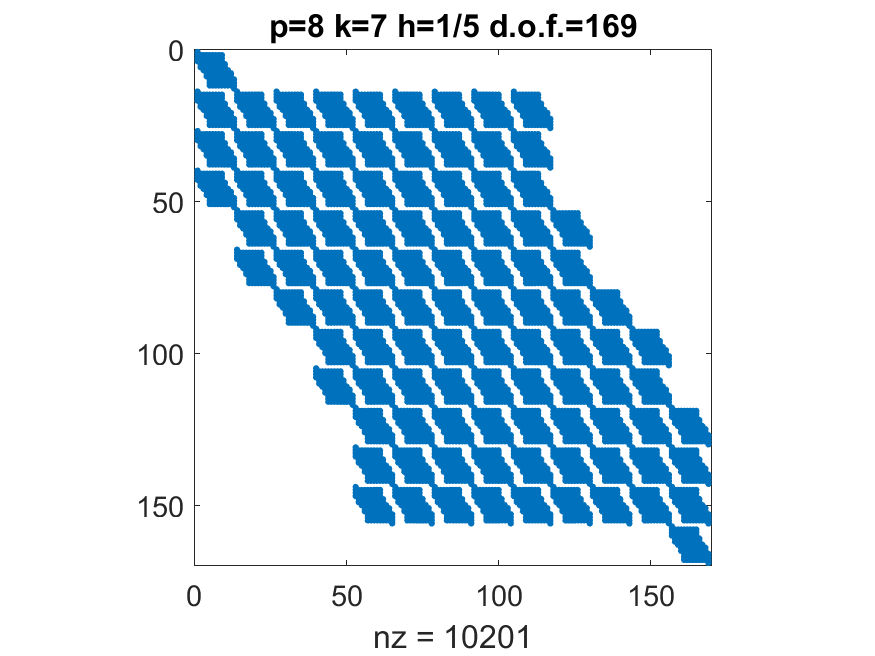} & 
\hspace{-5mm}\includegraphics[scale=0.40]{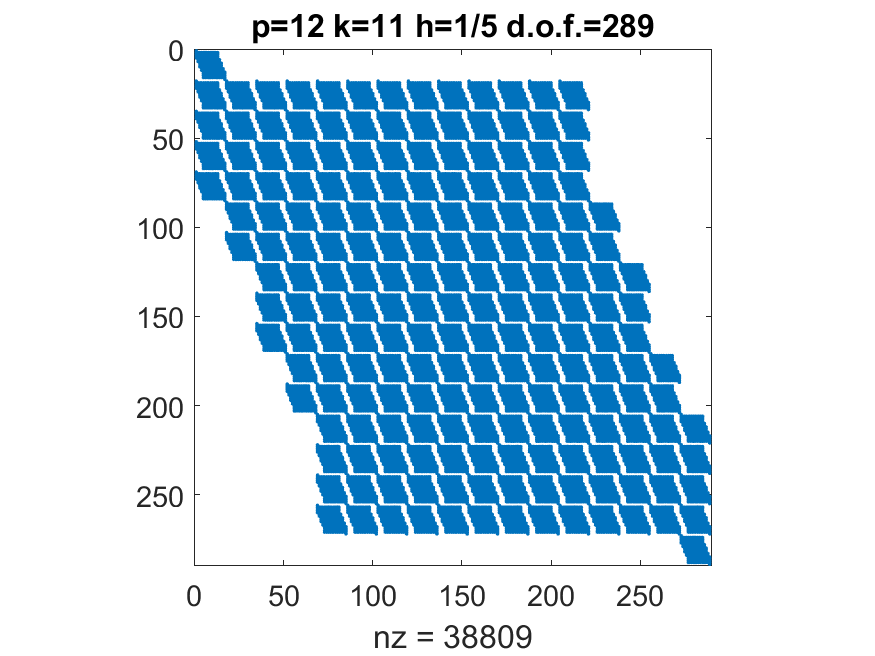}  \\
\end{tabular}
}
\vspace{-2mm}
\caption{Stiffness matrix sparsity pattern for the acoustic wave problem with Dirichlet boundary conditions: $p = 4$ (left), 
 $p=8$ (center), $p=12$ (right), with $k=1$ (top) or $k=p-1$ (bottom), fixed $h=1/5$, $\Delta t=0.01$, $\gamma=0.5$, $\beta=0.5$.
\label{spy_S_h5_VSp_Laplacian_DIR}}
\end{figure}

\begin{figure} %[!t]
\vspace{4mm}
\centerline{
\begin{tabular}{ccc}
\includegraphics[scale=0.40]{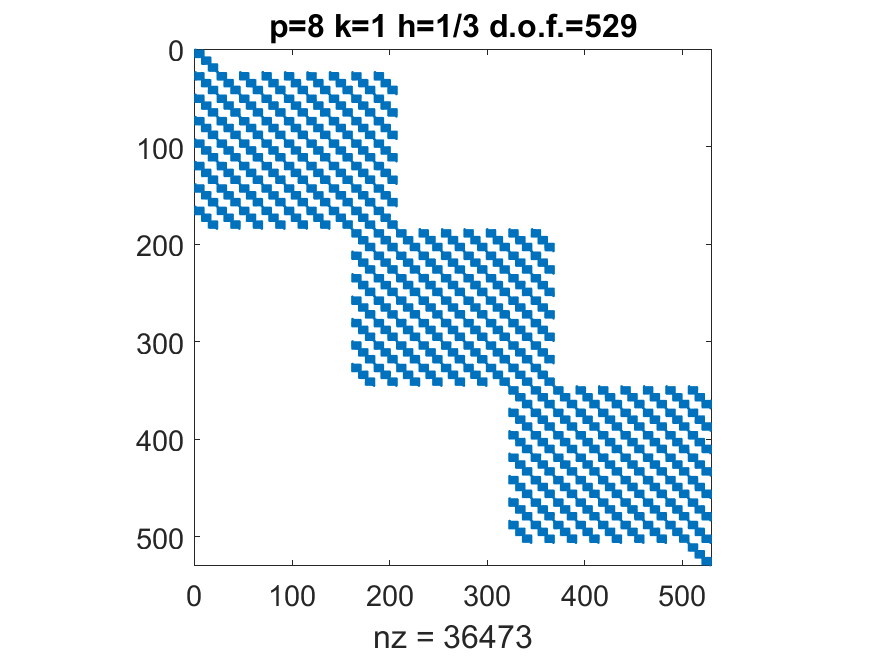} &
\hspace{-5mm}\includegraphics[scale=0.40]{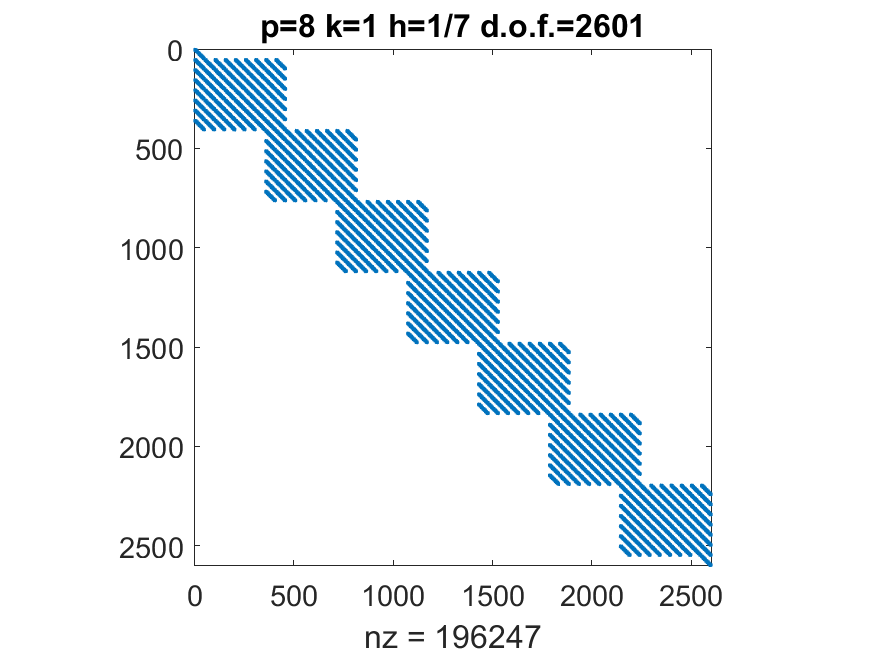} &
\hspace{-5mm}\includegraphics[scale=0.40]{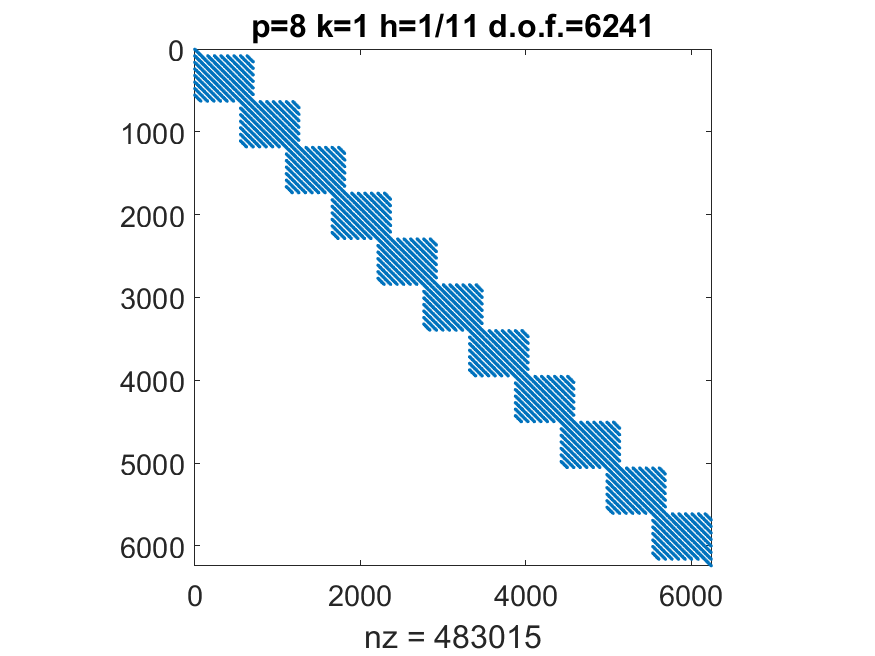} \\
\includegraphics[scale=0.40]{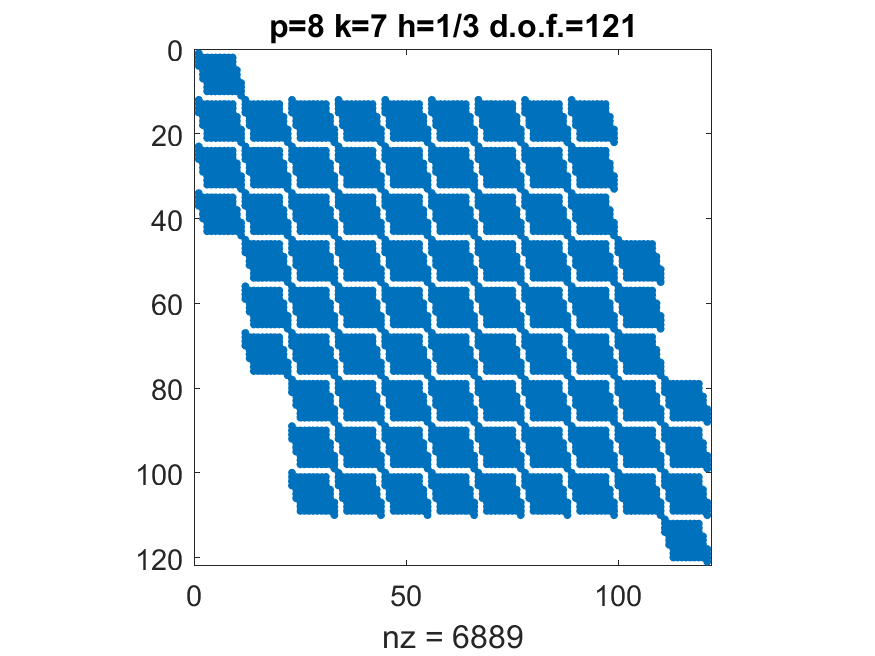} & 
\hspace{-5mm}\includegraphics[scale=0.40]{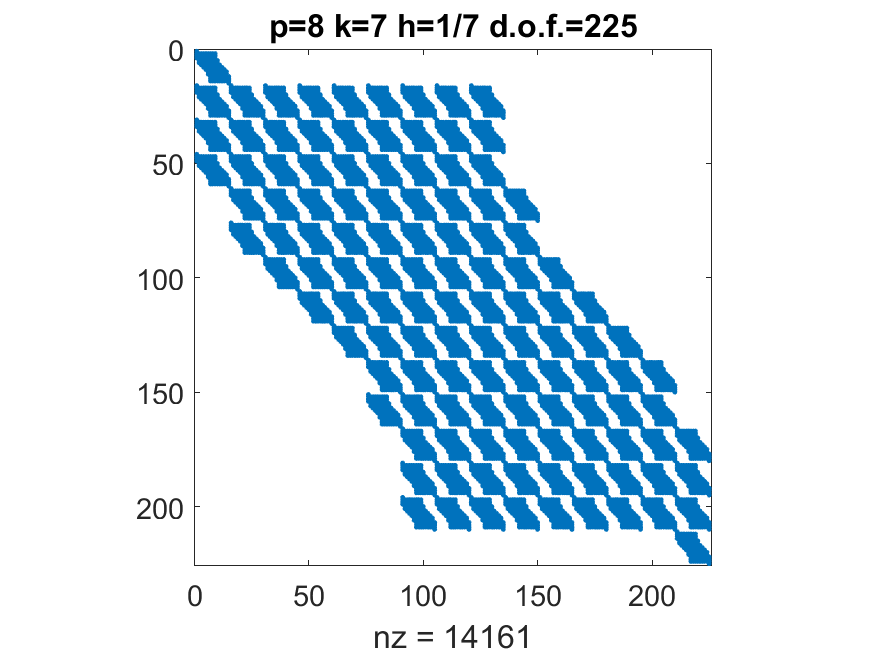} & 
\hspace{-5mm}\includegraphics[scale=0.40]{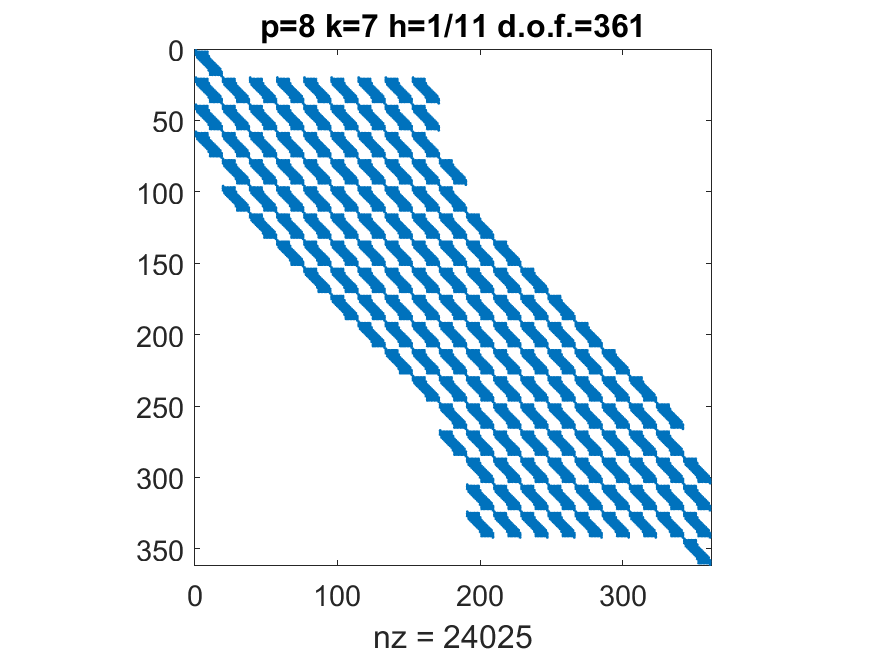}  \\
\end{tabular}
}
\vspace{-2mm}
\caption{Stiffness matrix sparsity pattern for the acoustic wave problem with Dirichlet boundary conditions: for $h=1/3$ (left),
 $1/h=7$ (center), $1/h=11$ (right), with $k=1$ (top) or $k=p-1$ (bottom), fixed $p=8$, $\Delta t=0.01$, $\gamma=0.5$, $\beta=0.5$.
\label{spy_S_p8_VSh_Laplacian_DIR}}
\end{figure}

Figs. \ref{spy_S_h5_VSp_Laplacian_DIR} and \ref{spy_S_p8_VSh_Laplacian_DIR}  shows the sparsity pattern of the stiffness matrix $\cal K$ for the same problem and parameter settings as in Figs.  \ref{eig_S_h5_VSp_Laplacian_DIR} and \ref{eig_S_p8_VSh_Laplacian_DIR}, respectively. The results are analogous to those for the mass matrices in Figs. \ref{spy_M_h5_VSp} and \ref{spy_M_p8_VSh}, with block-diagonal matrices for minimal regularity and almost full matrices in the case of maximal regularity. Again, both {\textsf { d.o.f.}} and {\textsf  {nz}} decrease for increasing regularity.

{\bf{Eigenvalues and condition number of the stiffness matrix with absorbing boundary conditions. }}
%\label{stiff_ABC}

\begin{figure} %[!t]
\vspace{-10mm}
\centerline{
\begin{tabular}{cc}
\includegraphics[scale=0.485]{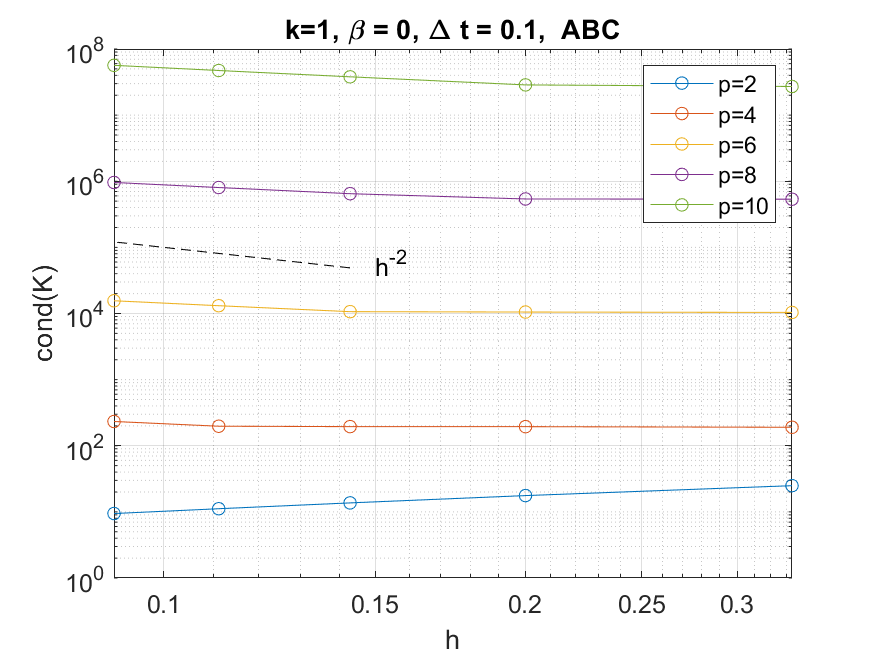} &
\includegraphics[scale=0.485]{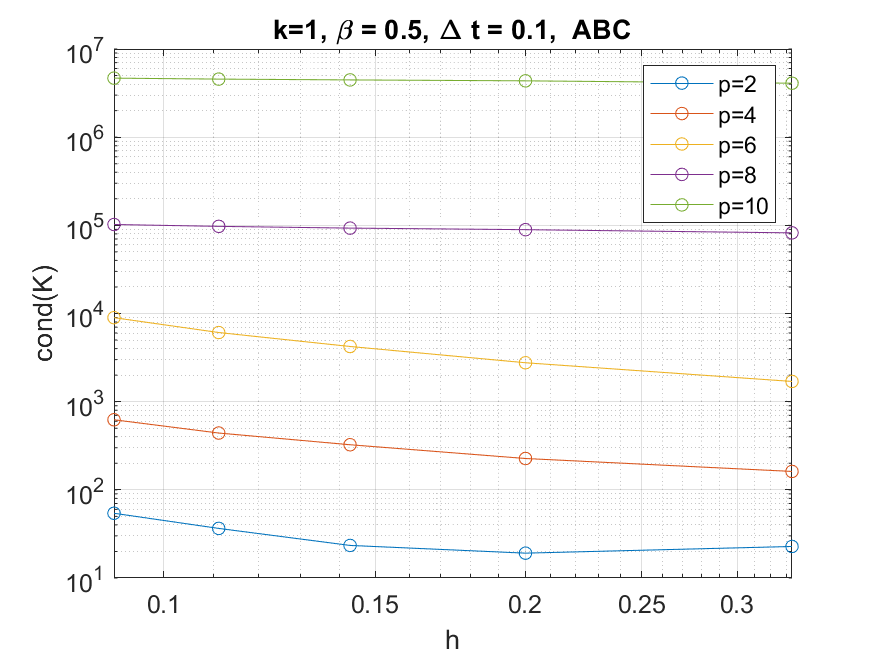} \\
\includegraphics[scale=0.485]{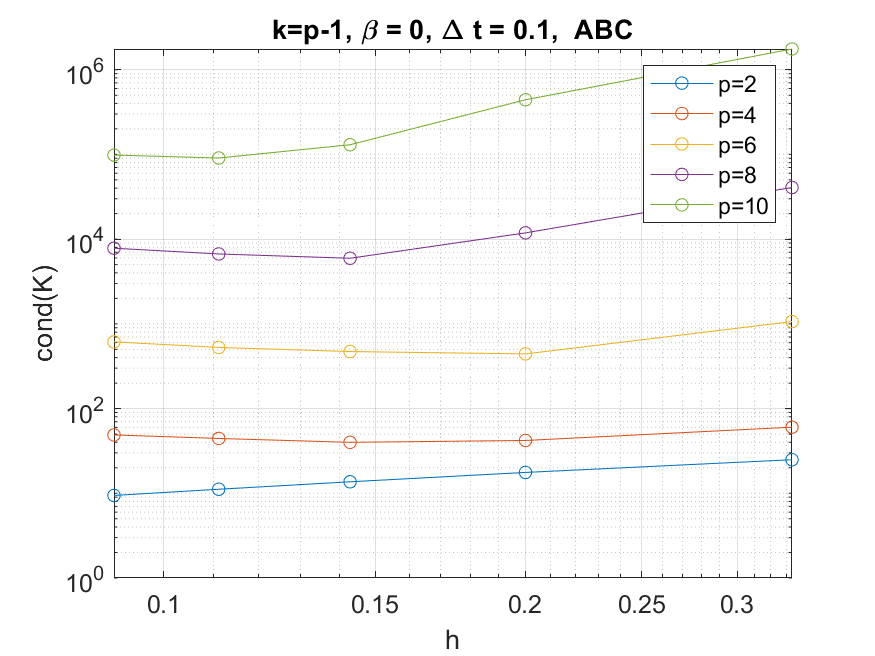} & \includegraphics[scale=0.485]{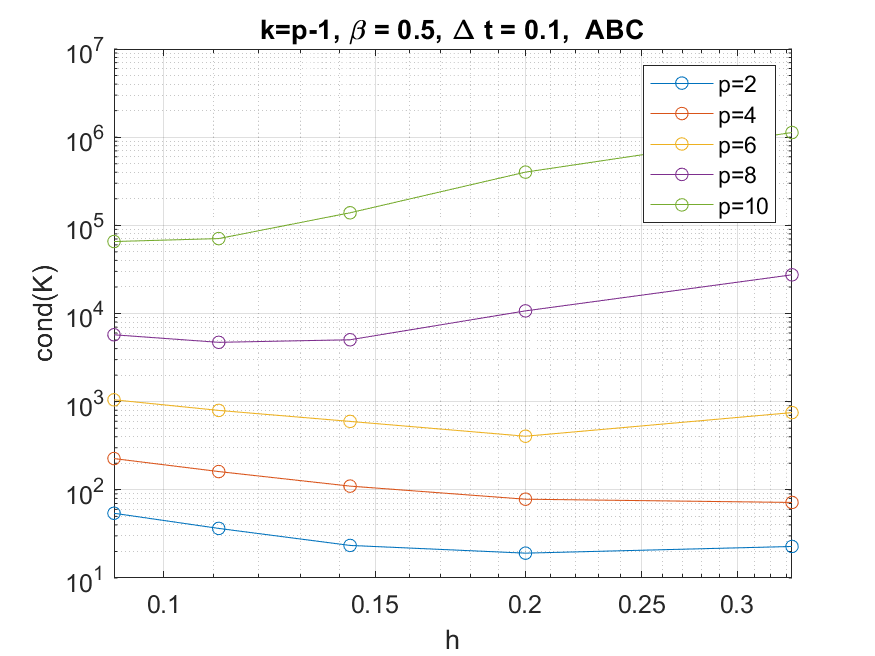} \\
\includegraphics[scale=0.485]{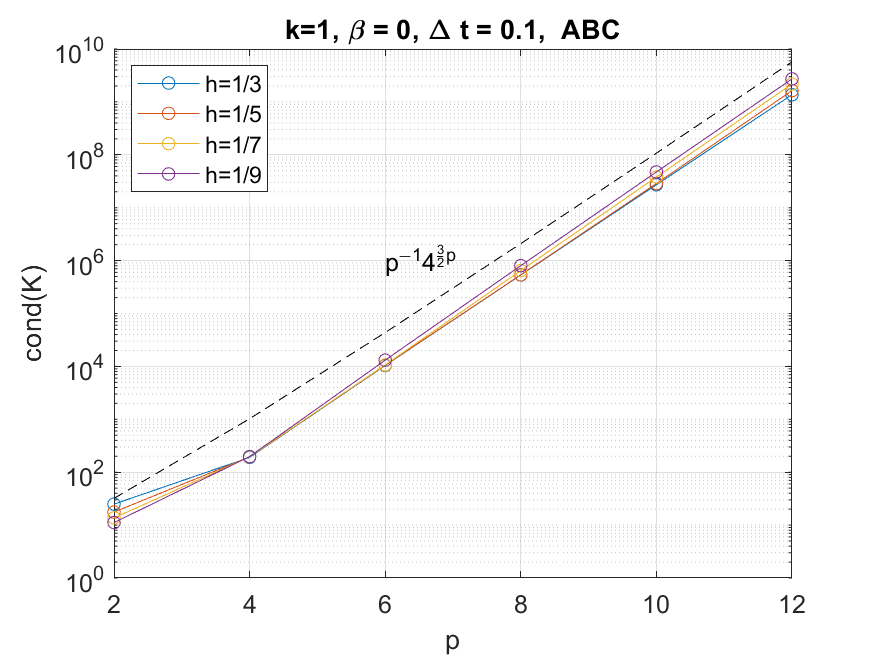} & 
\includegraphics[scale=0.485]{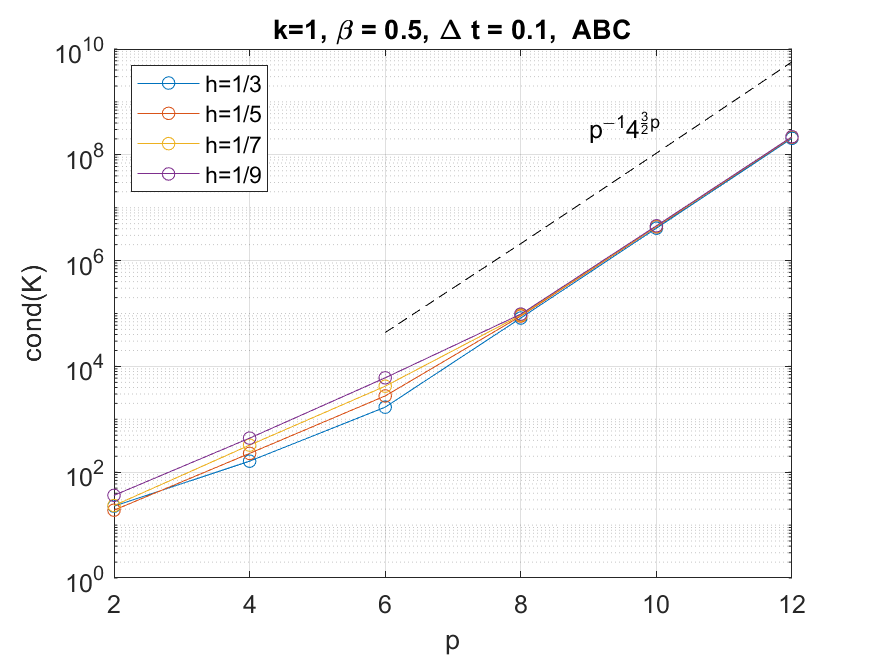} \\
\includegraphics[scale=0.485]{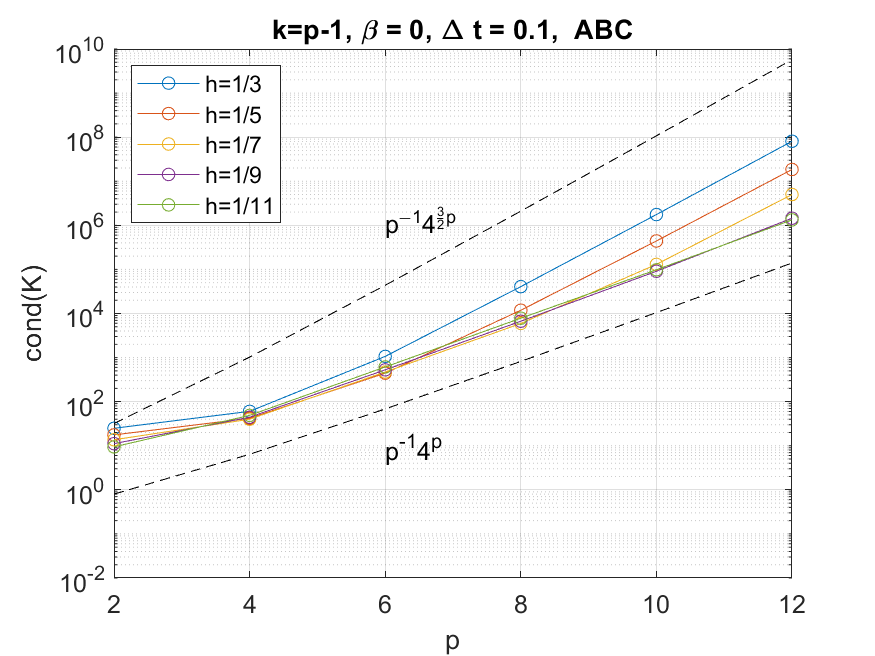} &
\includegraphics[scale=0.485]{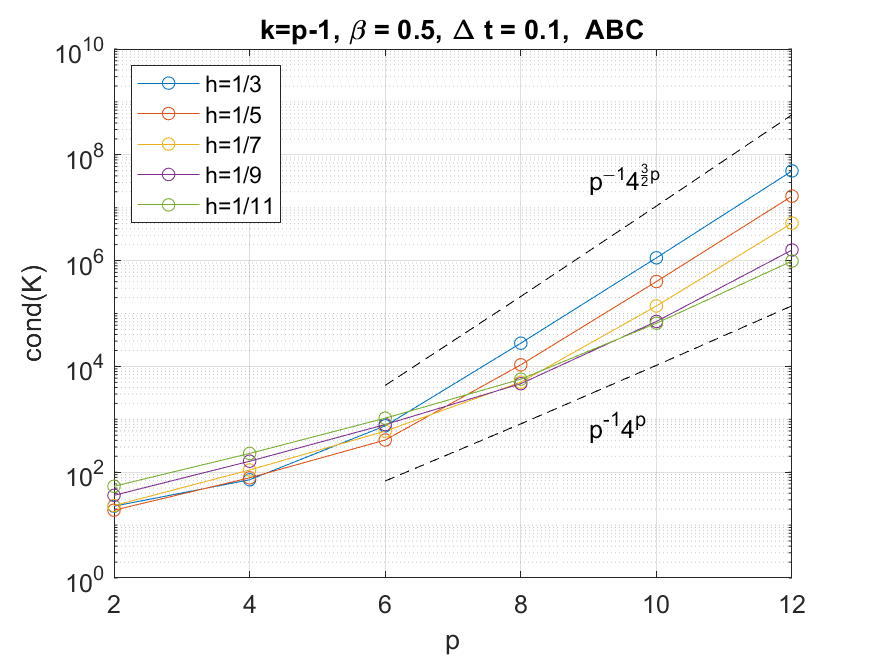} \\
\end{tabular}
}
\caption{Condition number of the stiffness matrix for the acoustic wave problem with absorbing boundary conditions (\ref{ABC}),  $\Delta t = 0.1$, $\gamma=0.5$, $\beta = 0$ (explicit Newmark, left), $\beta = 0.5$ (implicit Newmark, right). From the top to the bottom, vs.: (1)  $h$, for $p=2,\ 4,\ 6,8,\ 10$,  $k=1$; (2) $h$, for $p=2,\ 4,\ 6,8,\ 10$,  $k=p-1$; (3)  $p$, for $h=1/3,\ 1/5,\ 1/7,\ 1/9$,   $k=1$; (4) $p$, for $h=1/3,\ 1/5,\ 1/7,\ 1/9,\ 1/11$, $k=p-1$.
\label{cond_S_dt01_ABC}}
\end{figure}

\begin{figure} %[!t]
\vspace{-10mm}
\centerline{
\begin{tabular}{cc}
\includegraphics[scale=0.485]{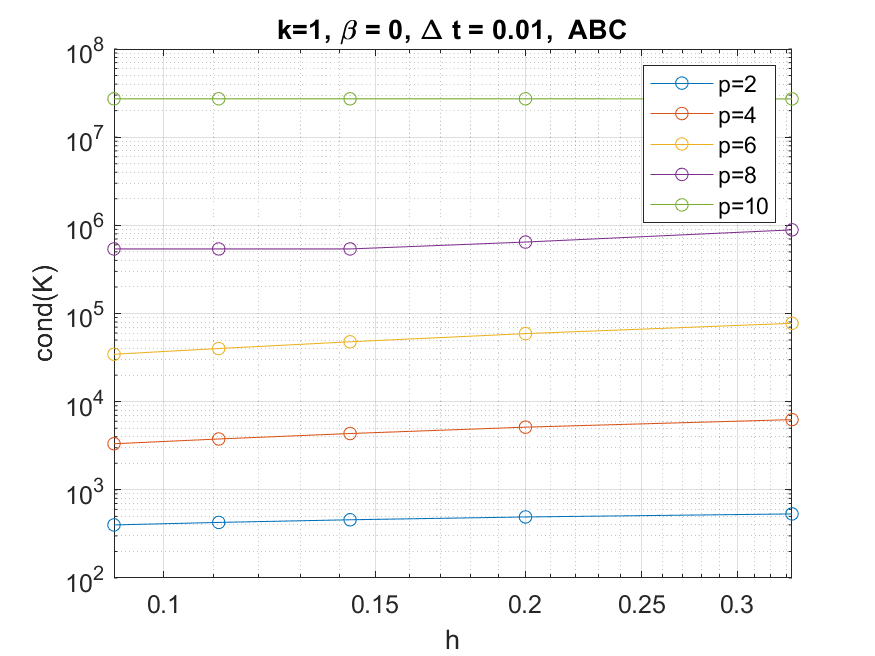} &
\includegraphics[scale=0.485]{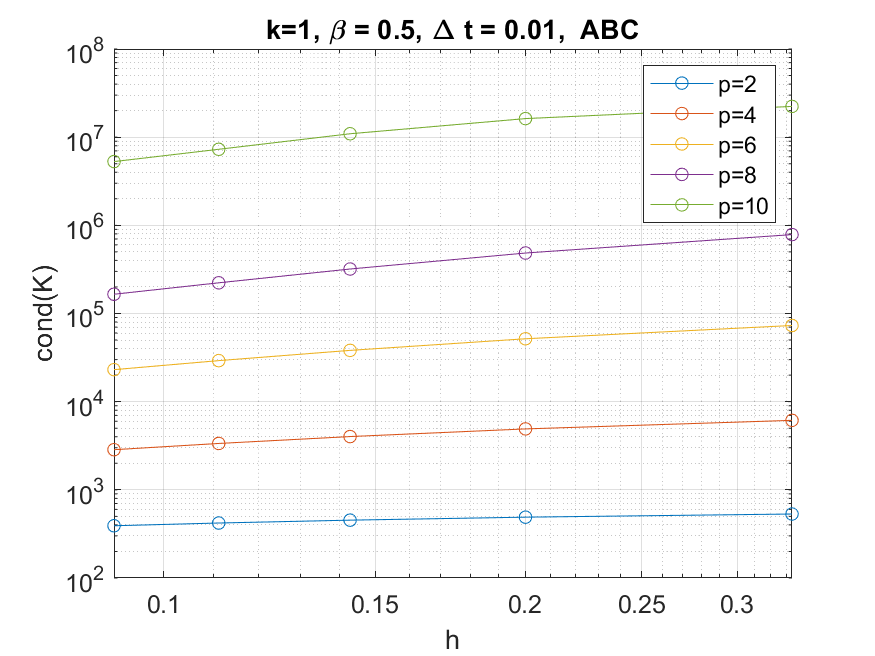} \\
\includegraphics[scale=0.485]{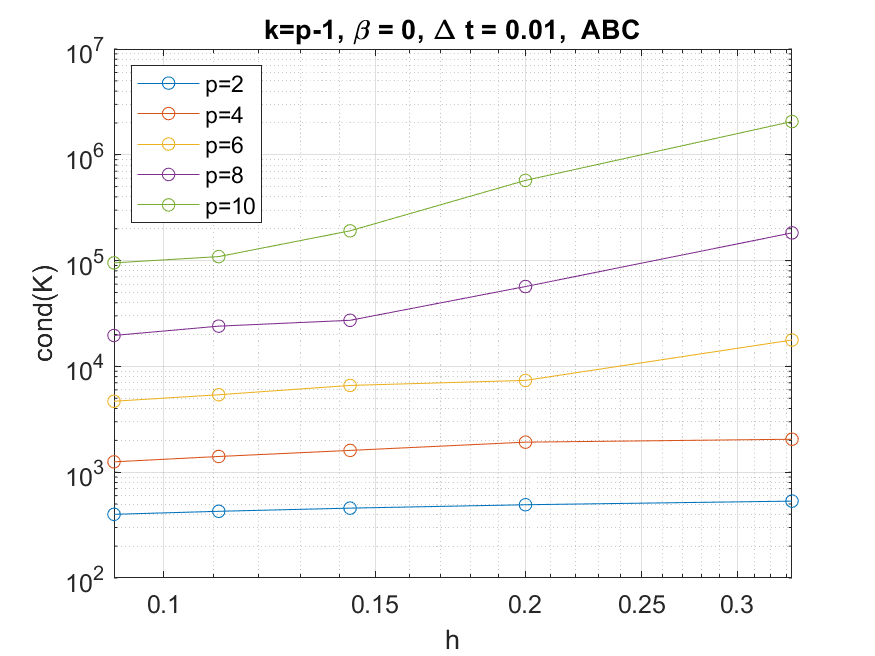} & \includegraphics[scale=0.485]{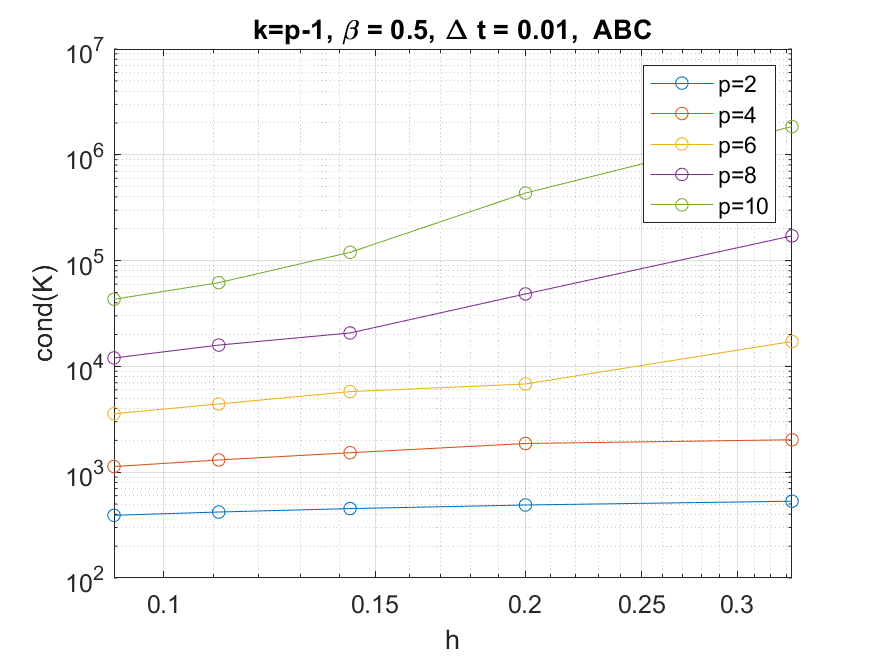} \\
\includegraphics[scale=0.485]{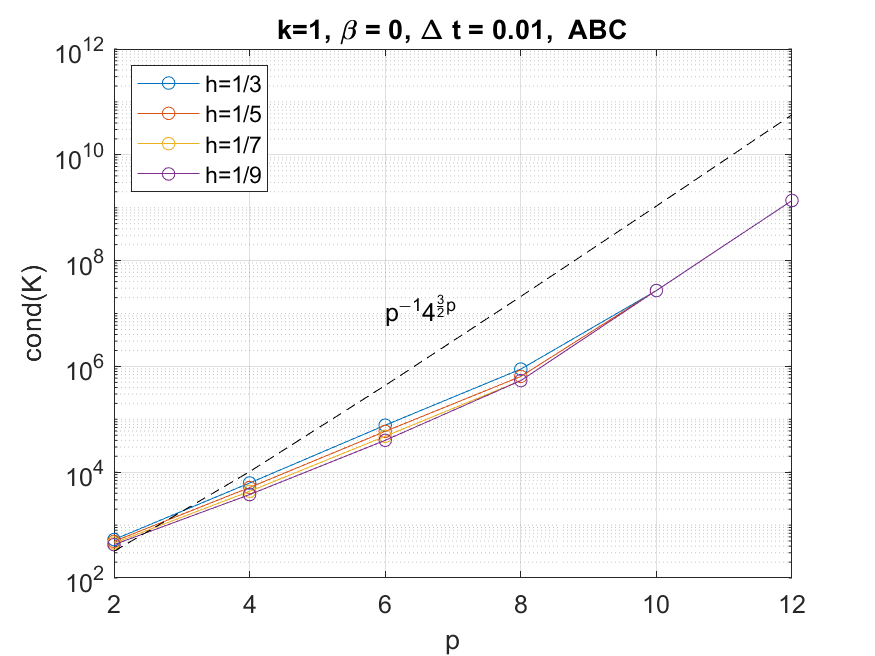} & 
\includegraphics[scale=0.485]{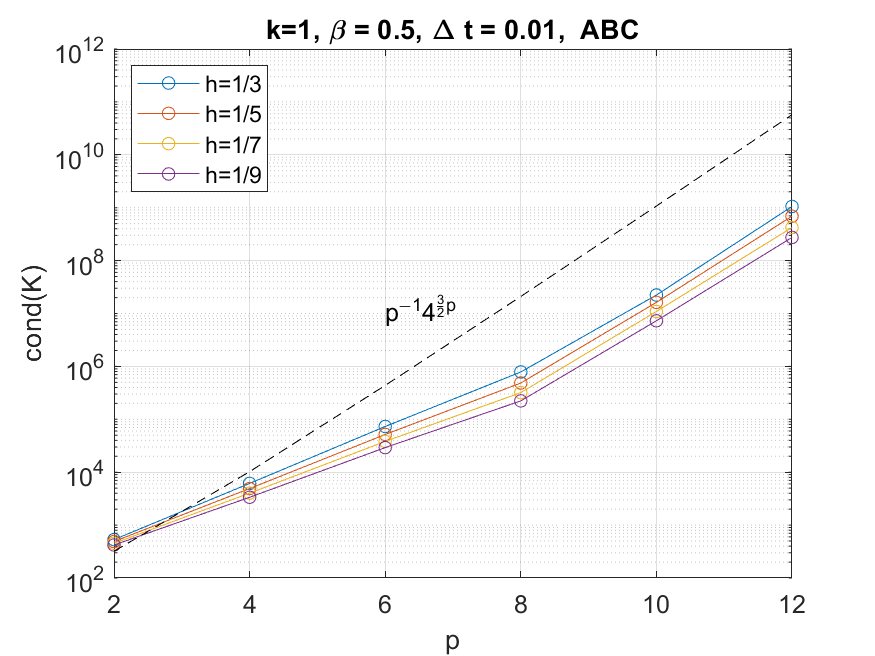} \\
\includegraphics[scale=0.485]{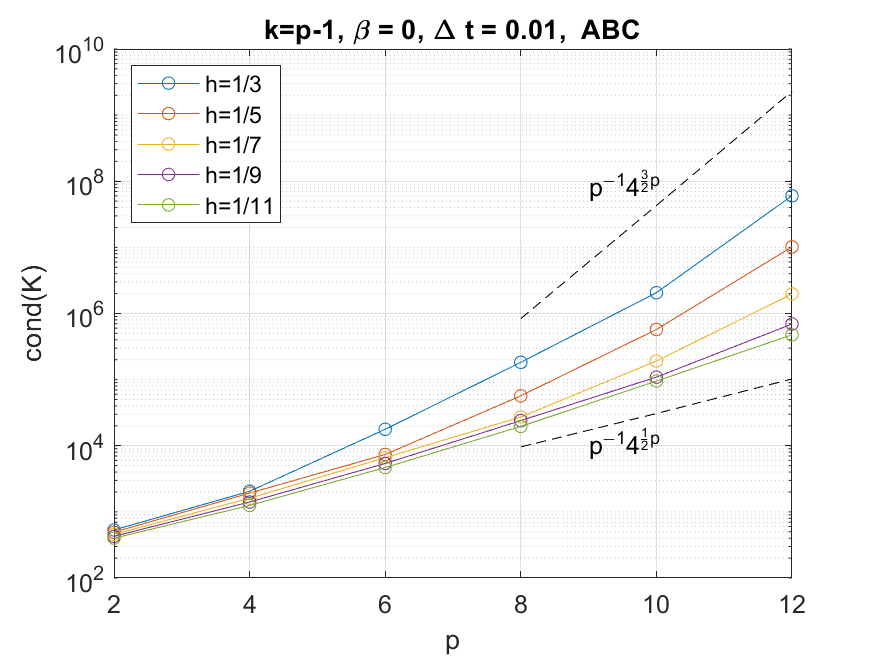} &
\includegraphics[scale=0.485]{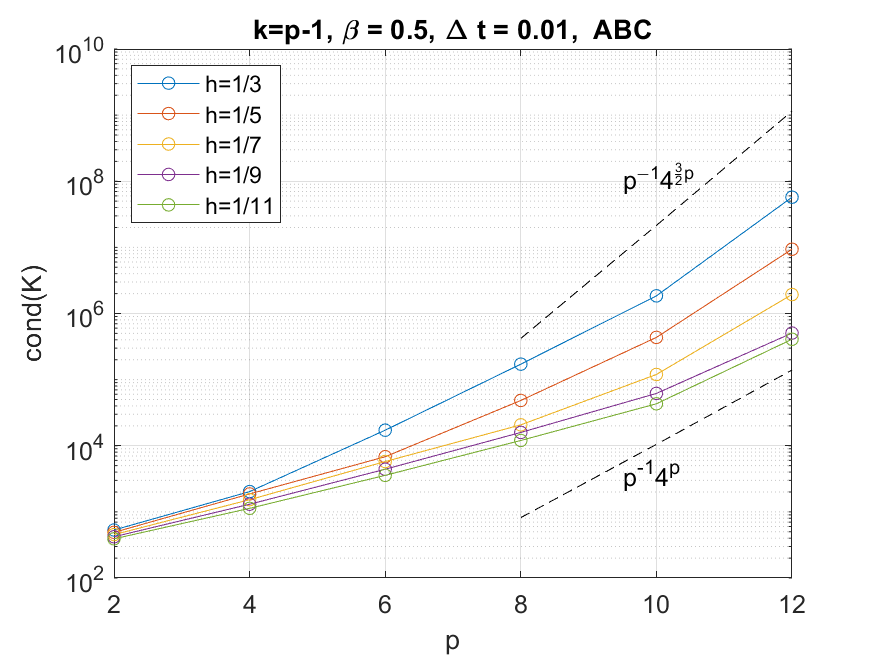} \\
\end{tabular}
}
\caption{Condition number of the stiffness matrix for the acoustic wave problem with absorbing boundary conditions (\ref{ABC}),  $\Delta t = 0.01$, $\gamma=0.5$, $\beta = 0$ (explicit Newmark, left), $\beta = 0.5$ (implicit Newmark, right). From the top to the bottom, vs.: (1)  $h$, for $p=2,\ 4,\ 6,8,\ 10$,  $k=1$; (2) $h$, for $p=2,\ 4,\ 6,8,\ 10$,  $k=p-1$; (3)  $p$, for $h=1/3,\ 1/5,\ 1/7,\ 1/9$,   $k=1$; (4) $p$, for $h=1/3,\ 1/5,\ 1/7,\ 1/9,\ 1/11$, $k=p-1$.
\label{cond_S_dt001_ABC}}
\end{figure}

In Fig.  \ref{cond_S_dt01_ABC} and \ref{cond_S_dt001_ABC}  we report the condition numbers {\textsf{cond}}($\cal K$) of the stiffness matrix for the acoustic wave problem with absorbing boundary conditions (\ref{ABC}), using the same setting of Figs. \ref{cond_S_dt01_DIR} and \ref{cond_S_dt001_DIR}, respectively.
  We observe that the results are analogous to those for  Neumann conditions.
 
\begin{figure} %[!t]
\vspace{-10mm}
\centerline{
\begin{tabular}{cc}
\includegraphics[scale=0.5]{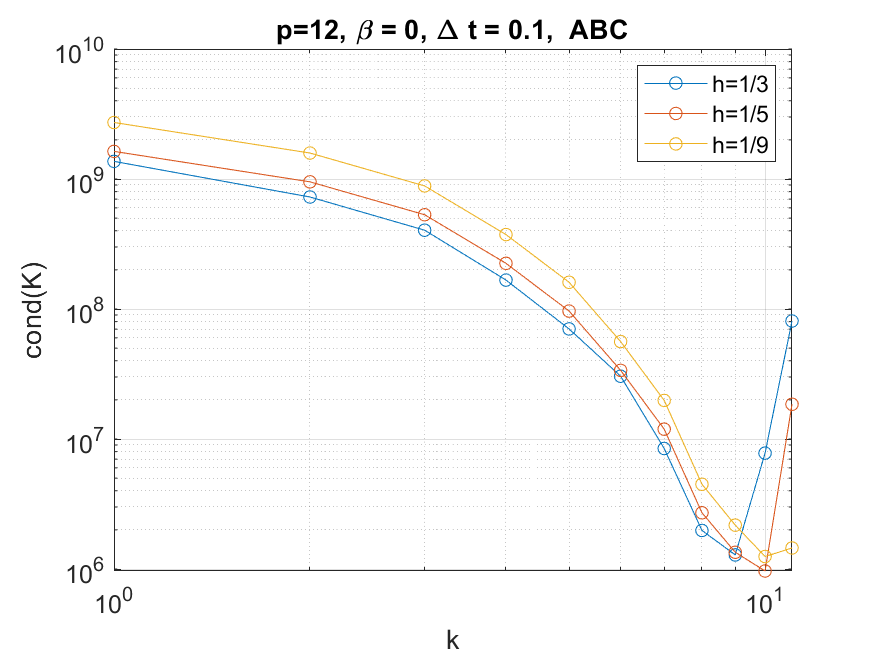} &
\includegraphics[scale=0.5]{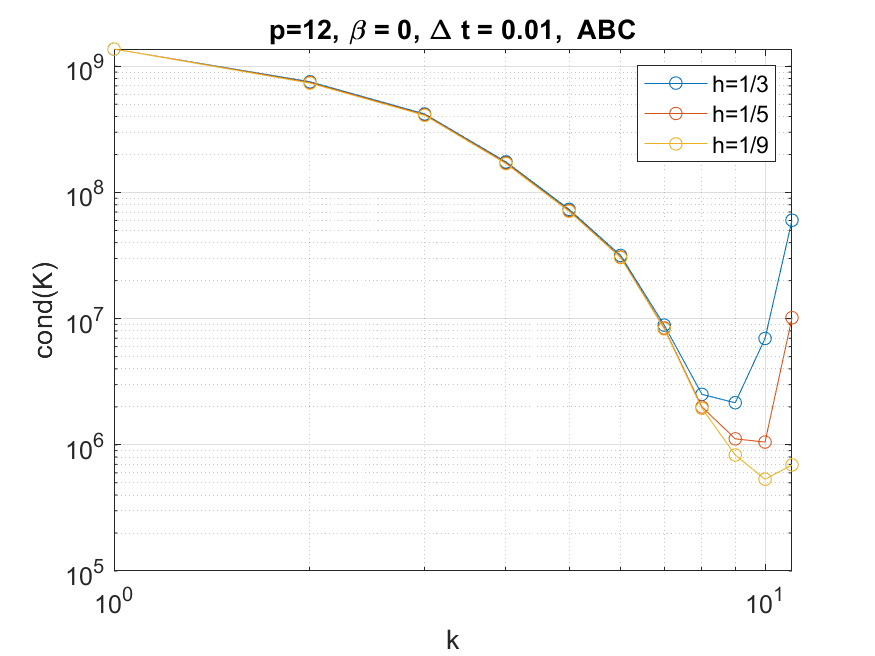} \\
\includegraphics[scale=0.5]{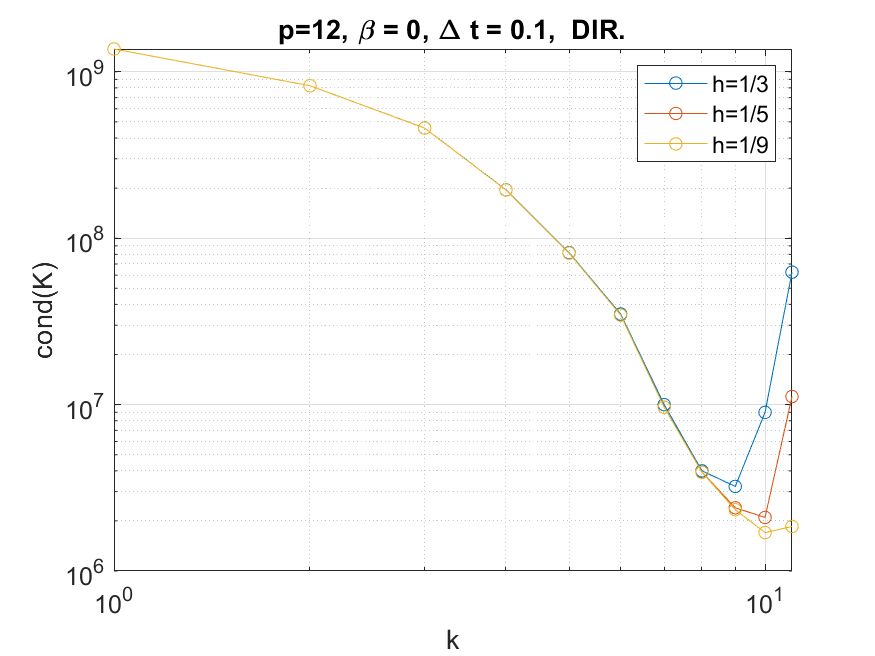} & \includegraphics[scale=0.5]{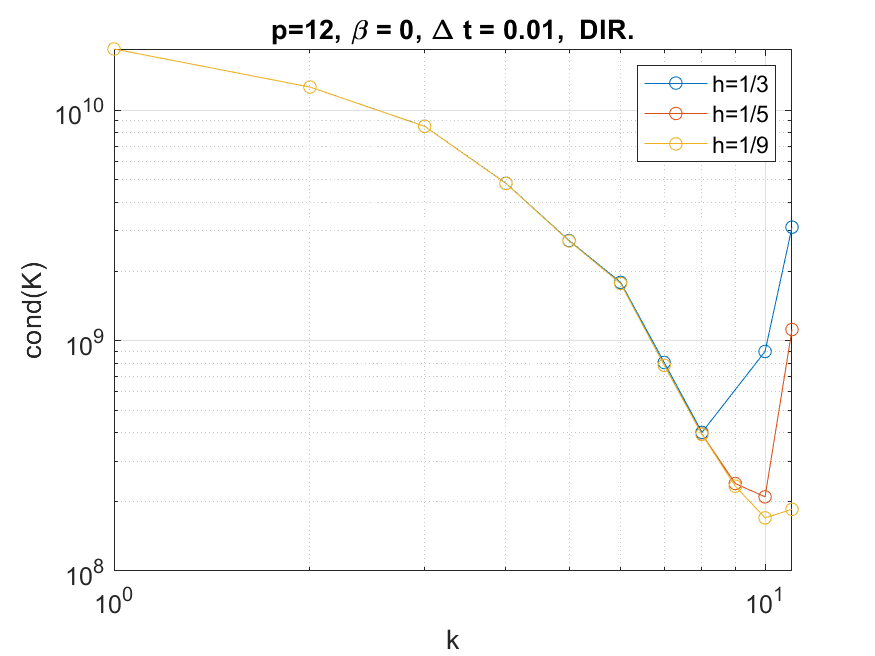} \\
\includegraphics[scale=0.5]{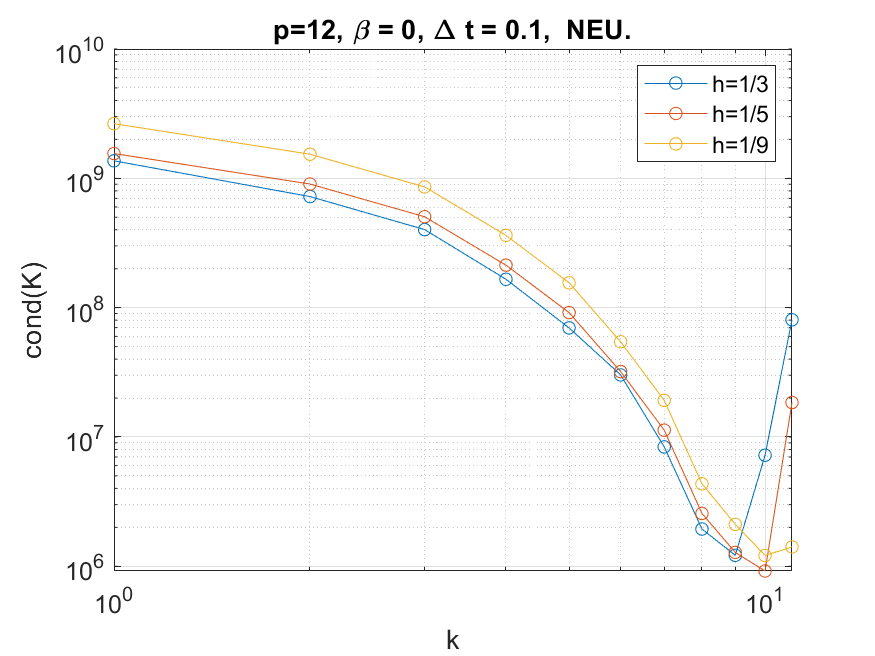} & 
\includegraphics[scale=0.5]{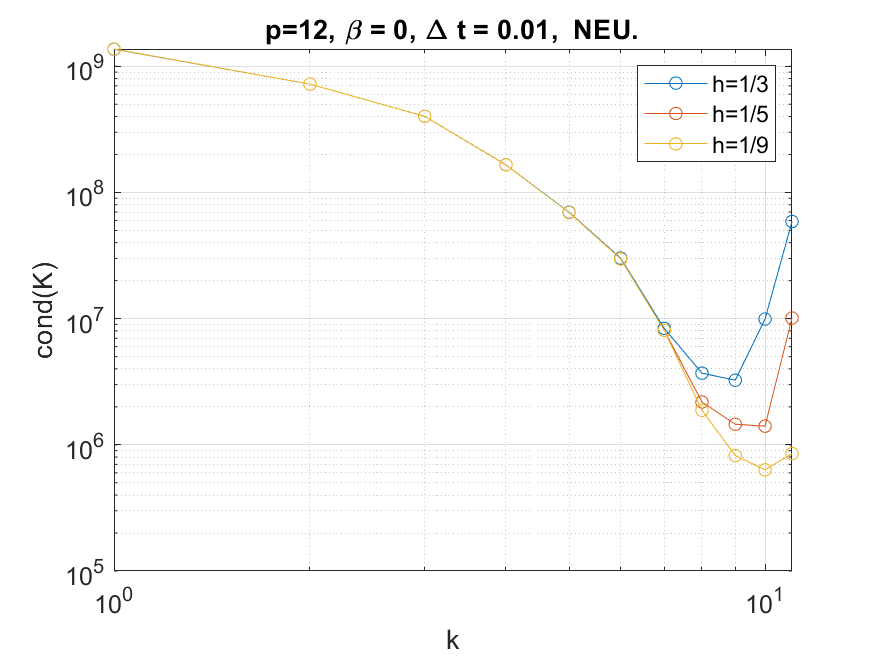} \\
\end{tabular}
}
\caption{Condition numbers {\textsf{cond}}($\cal K$) of the stiffness matrix for the acoustic wave problem with different types of boundary conditions: absorbing (top), Dirichlet (center), Neumann (bottom), as a function of  $k$, $\beta=0$ (explicit Newmark), $\gamma=0.5$, degree $p=12$, $\Delta t=0.1$ (left) and $\Delta t=0.01$ (right),  $h=1/3,\ 1/5, \ 1/9$.
\label{cond_S_k}}
\end{figure}

In Fig.  \ref{cond_S_k}  we report the condition numbers {\textsf{cond}}($\cal K$) of the stiffness matrix for the acoustic wave problem with different types of boundary conditions: absorbing (top), Dirichlet (center), Neumann (bottom). We vary the values of regularity $k$, and choose  $\beta=0$ (explicit Newmark), $\gamma=0.5$, degree $p=12$, $\Delta t=0.1$ (left) or $\Delta t=0.01$ (right), and three different values of  mesh size $h$. The results are analogous to those for the mass matrices in Fig. \ref{cond_M_k} for each type of boundary conditions.

\begin{figure} %[!t]
\vspace{4mm}
\centerline{
\begin{tabular}{ccc}
\includegraphics[scale=0.40]{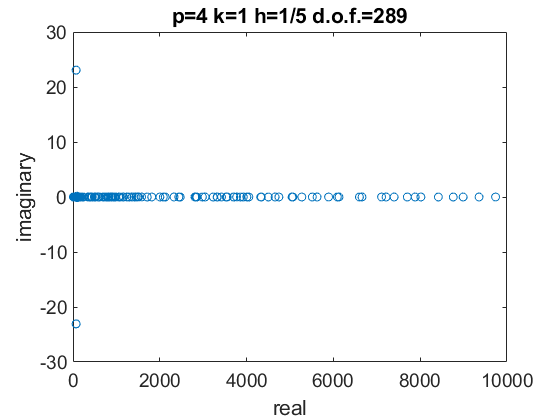} &
\hspace{-5mm}\includegraphics[scale=0.40]{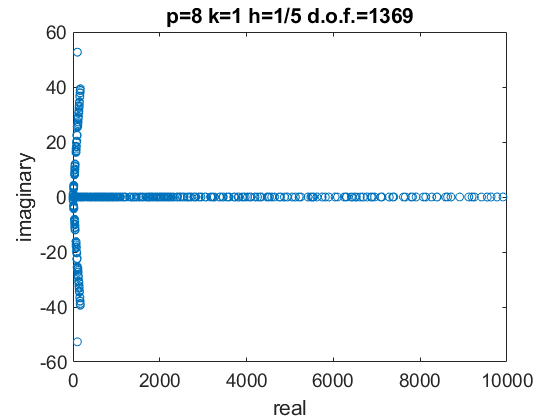} &
\hspace{-5mm}\includegraphics[scale=0.40]{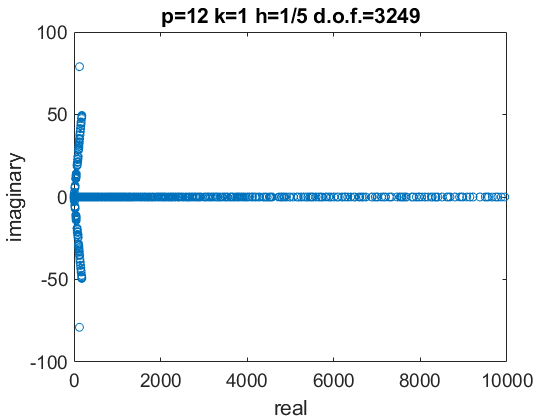} \\
\includegraphics[scale=0.40]{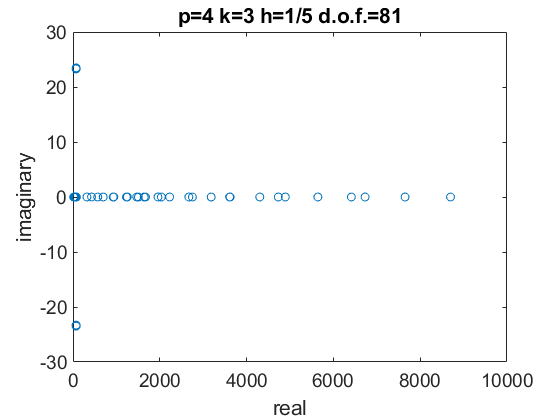} & 
\hspace{-5mm}\includegraphics[scale=0.40]{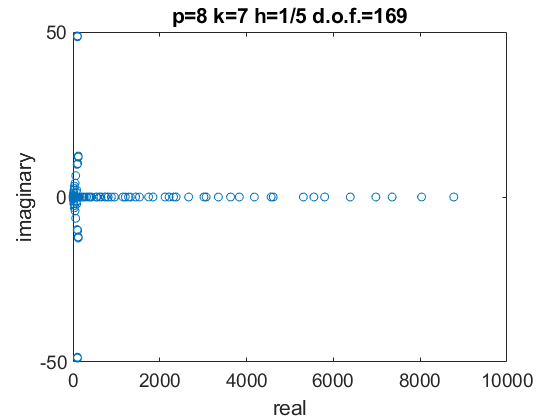} & 
\hspace{-5mm}\includegraphics[scale=0.40]{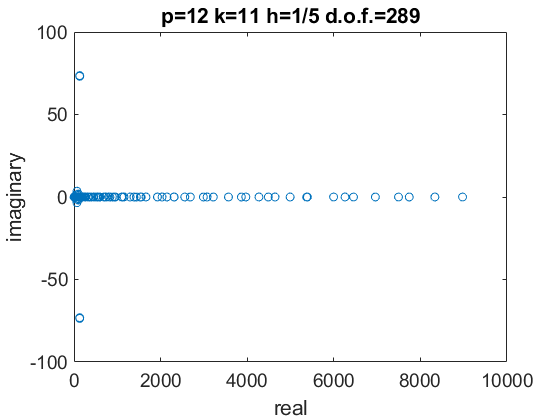}  \\
\end{tabular}
}
\vspace{-2mm}
\caption{Stiffness matrix eigenvalue distribution in the complex plane for the acoustic wave problem with  absorbing boundary conditions and explicit Newmark scheme ($\beta=0)$: $p = 4$ (left),
 $p=8$ (center), $p=12$ (right), with $k=1$ (top) or $k=p-1$ (bottom), fixed $h=1/5$, $\Delta t=0.01$, $\gamma=0.5$.
\label{eig_S_h5_VSp_explicit}}
\end{figure}

\begin{figure} %[!t]
\vspace{4mm}
\centerline{
\begin{tabular}{ccc}
\includegraphics[scale=0.40]{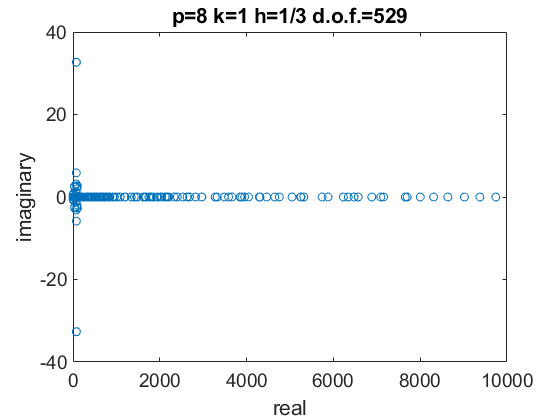} &
\hspace{-5mm}\includegraphics[scale=0.40]{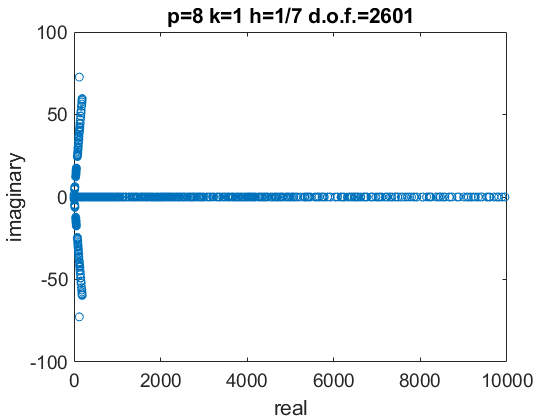} &
\hspace{-5mm}\includegraphics[scale=0.40]{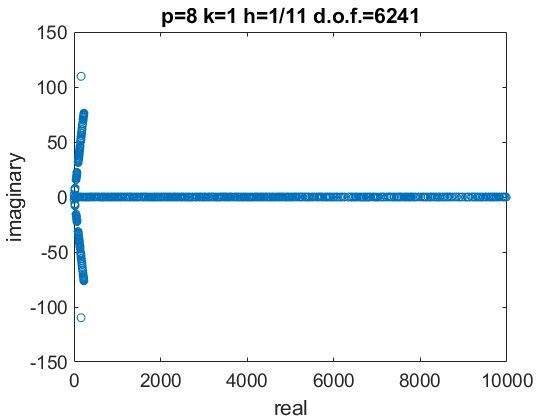} \\
\includegraphics[scale=0.40]{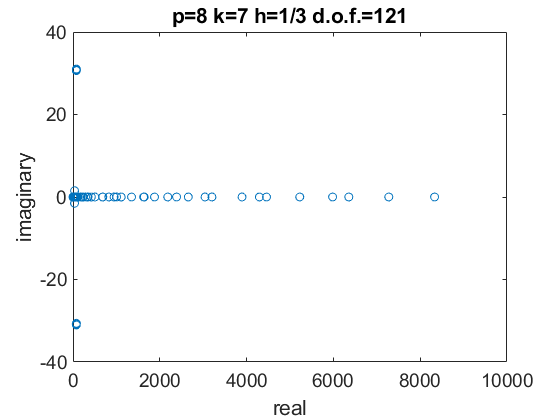} & 
\hspace{-5mm}\includegraphics[scale=0.40]{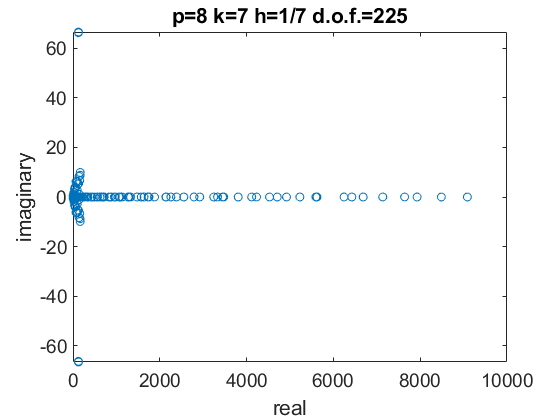} & 
\hspace{-5mm}\includegraphics[scale=0.40]{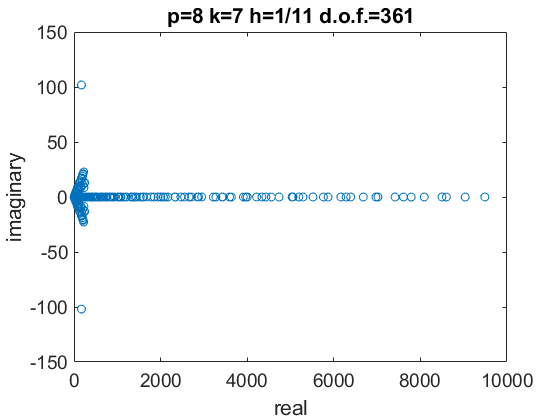}  \\
\end{tabular}
}
\vspace{-2mm}
\caption{Eigenvalue distribution in the complex plane of stiffness matrix for the acoustic wave problem with  absorbing boundary conditions and explicit Newmark scheme ($\beta=0)$: $h=1/3$ (left),
 $1/h=7$ (center), $1/h=11$ (right), with $k=1$ (top) or $k=p-1$ (bottom), fixed $p=8$, $\Delta t=0.01$, $\gamma=0.5$.
\label{eig_S_p8_VSh_explicit}}
\end{figure}

In Fig. \ref{eig_S_h5_VSp_explicit} and \ref{eig_S_p8_VSh_explicit} we report the stiffness matrix $\cal K$ eigenvalue distribution in complex plane for the acoustic wave problem with absorbing boundary conditions,  and explicit Newmark scheme ($\beta=0)$, for $\Delta t=0.01$ and $\gamma=0.5$. Results refer to the same  parameter settings as in Figs. \ref{eig_S_h5_VSp_Laplacian_DIR} and \ref{eig_S_p8_VSh_Laplacian_DIR}.  The eigenvalues real parts belong to an interval $[0, r]$ with $r$   independent of all parameters  $1/h$, $p$ and $k$. Here the value $r$ is much larger than  in the case of Dirichlet boundary conditions. Furthermore, the  real  parts of complex eigenvalues are almost  close to zero.

\begin{figure} %[!t]
\vspace{4mm}
\centerline{
\begin{tabular}{ccc}
\includegraphics[scale=0.40]{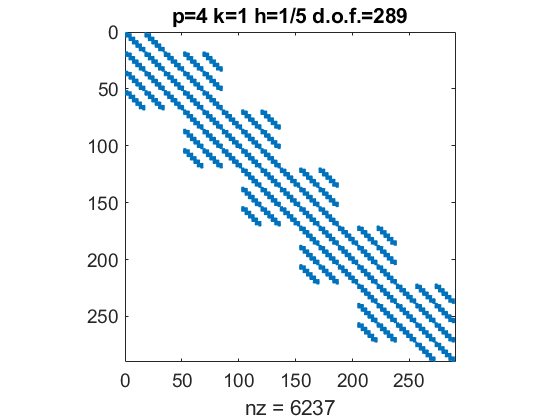} &
\hspace{-5mm}\includegraphics[scale=0.40]{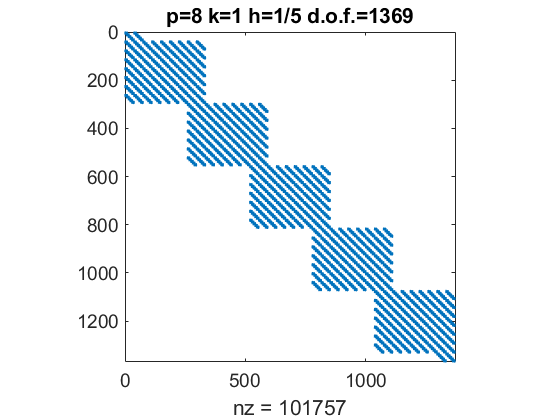} &
\hspace{-5mm}\includegraphics[scale=0.40]{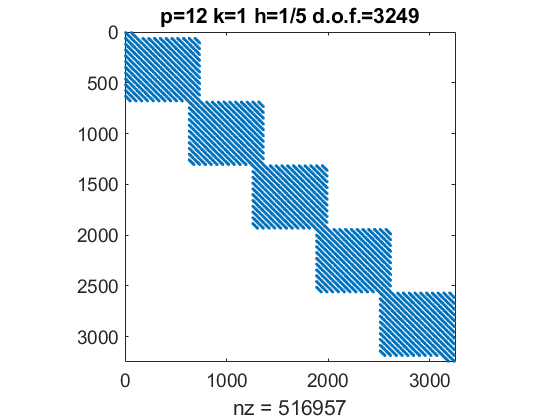} \\
\includegraphics[scale=0.40]{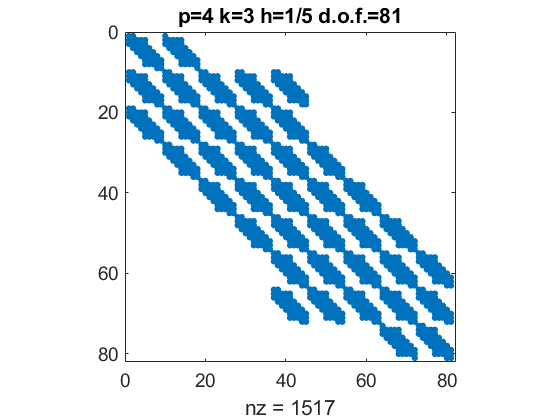} & 
\hspace{-5mm}\includegraphics[scale=0.40]{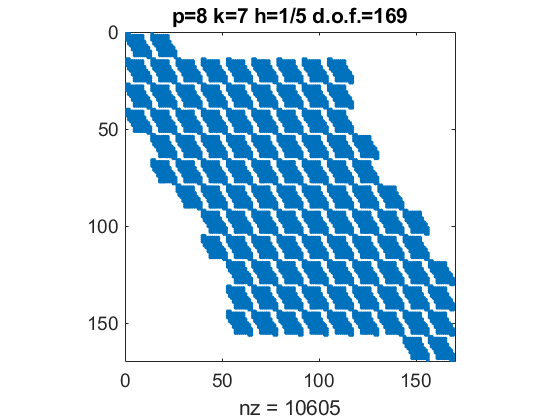} & 
\hspace{-5mm}\includegraphics[scale=0.40]{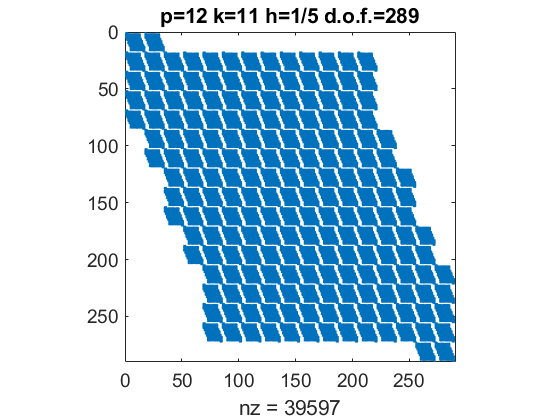}  \\
\end{tabular}
}
\vspace{-2mm}
\caption{Stiffness matrix sparsity pattern for the acoustic wave problem with  absorbing boundary conditions and explicit Newmark scheme ($\beta=0)$: $p = 4$ (left),
 $p=8$ (center), $p=12$ (right), with $k=1$ (top) or $k=p-1$ (bottom), fixed $h=1/5$, $\Delta t=0.01$, $\gamma=0.5$.
\label{spy_S_h5_VSp_explicit}}
\end{figure}

\begin{figure} %[!t]
\vspace{4mm}
\centerline{
\begin{tabular}{ccc}
\includegraphics[scale=0.40]{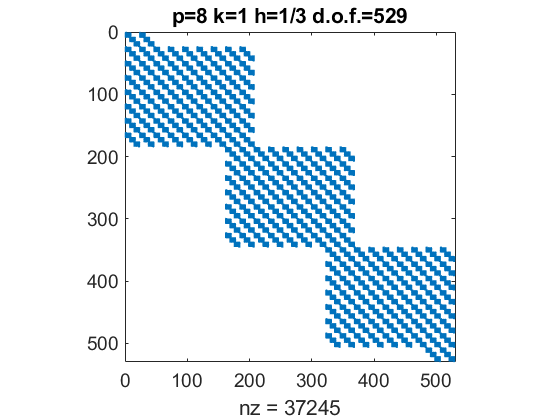} &
\hspace{-5mm}\includegraphics[scale=0.40]{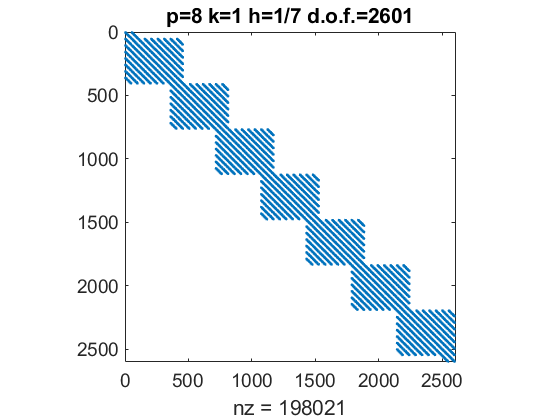} &
\hspace{-5mm}\includegraphics[scale=0.40]{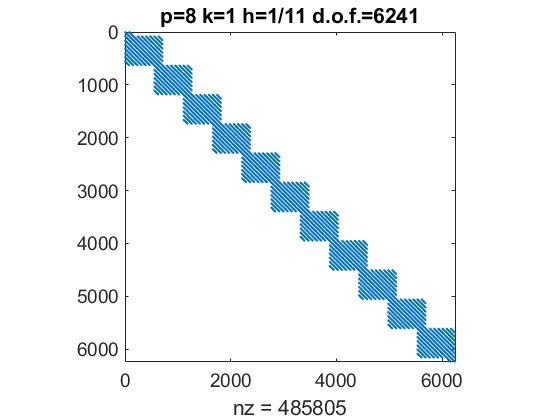} \\
\includegraphics[scale=0.40]{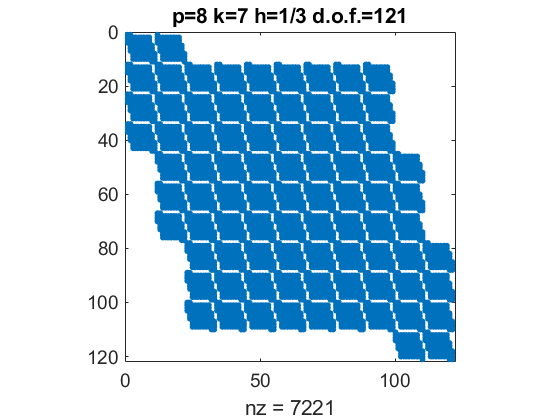} & 
\hspace{-5mm}\includegraphics[scale=0.40]{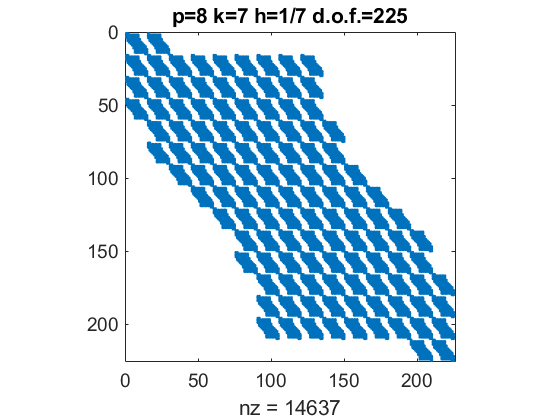} & 
\hspace{-5mm}\includegraphics[scale=0.40]{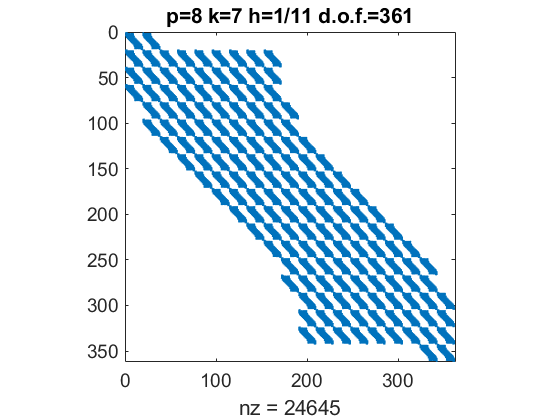}  \\
\end{tabular}
}
\vspace{-2mm}
\caption{Stiffness matrix sparsity pattern for the acoustic wave problem with  absorbing boundary conditions and explicit Newmark scheme ($\beta=0)$: $h=1/3$ (left),
 $1/h=7$ (center), $1/h=11$ (right), with $k=1$ (top) or $k=p-1$ (bottom), fixed $p=8$, $\Delta t=0.01$, $\gamma=0.5$.
\label{spy_S_p8_VSh_explicit}}
\end{figure}

\begin{figure} %[!t]
\vspace{4mm}
\centerline{
\begin{tabular}{ccc}
\includegraphics[scale=0.40]{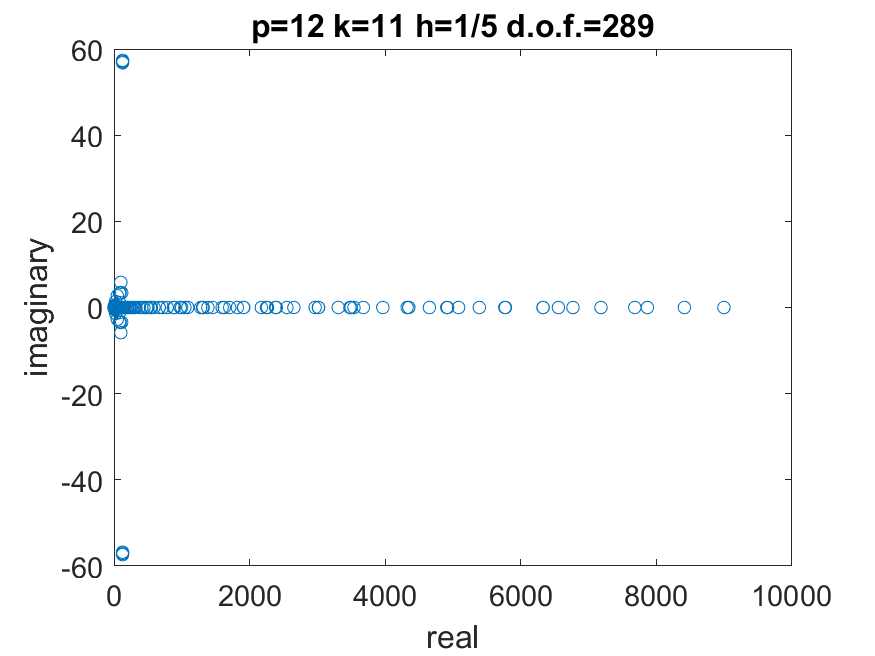} &
\hspace{-5mm}\includegraphics[scale=0.40]{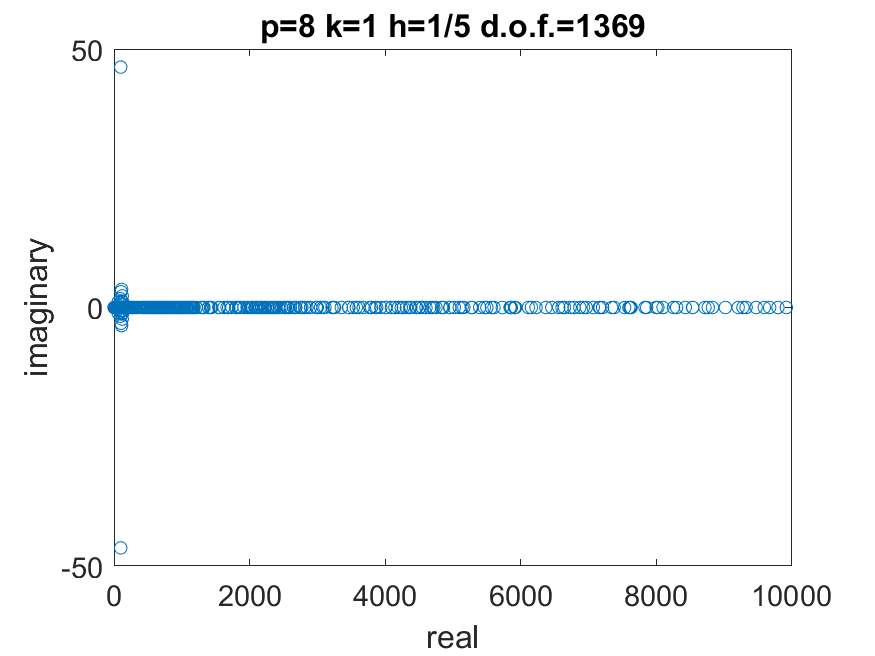} &
\hspace{-5mm}\includegraphics[scale=0.40]{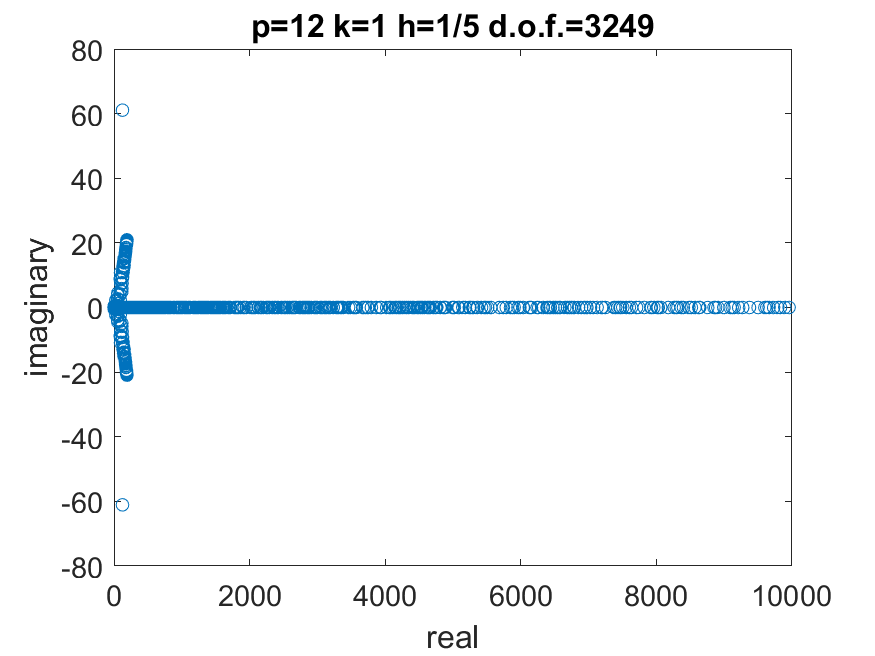} \\
\includegraphics[scale=0.40]{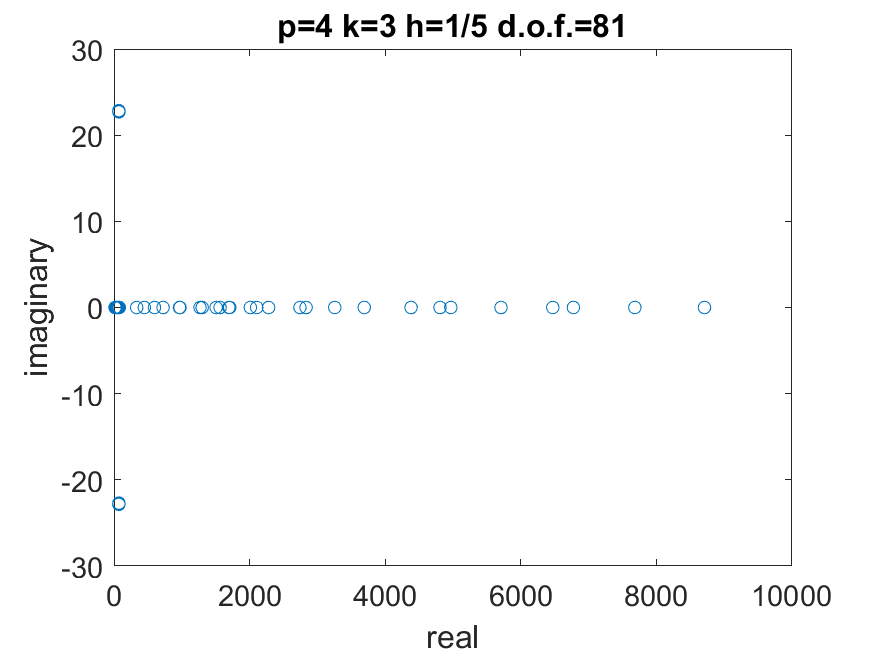} & 
\hspace{-5mm}\includegraphics[scale=0.40]{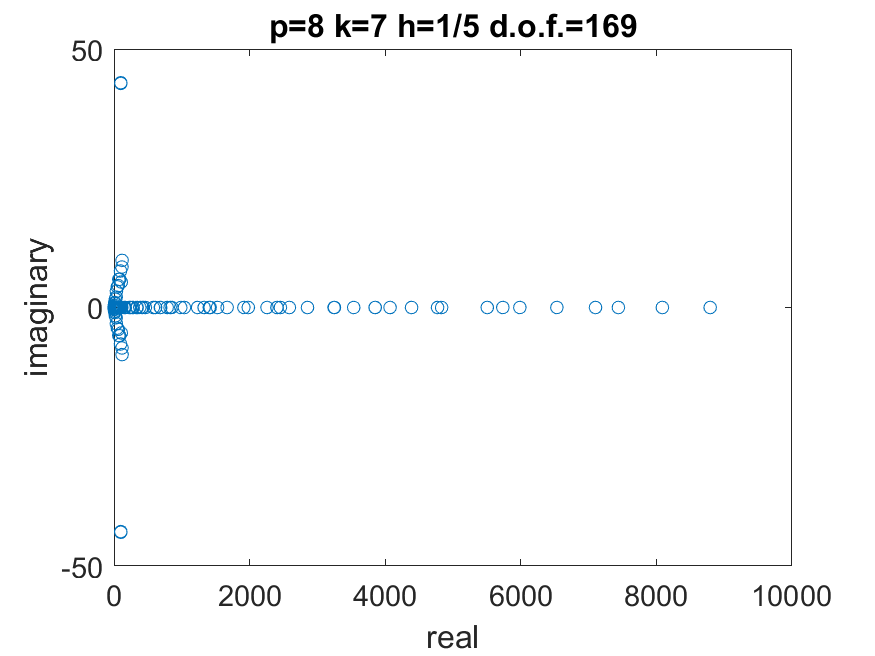} & 
\hspace{-5mm}\includegraphics[scale=0.40]{spy_eig_stiff/EigStiffVSpIMPLICIT/EIG_Stiff_Impl_ABC_p12k11h5dof289.png}  \\
\end{tabular}
}
\vspace{-2mm}
\caption{Stiffness matrix eigenvalue distribution in the complex plane for the acoustic wave problem with  absorbing boundary conditions and implicit Newmark scheme ($\beta=0.5)$: $p = 4$ (left),
 $p=8$ (center), $p=12$ (right), with $k=1$ (top) or $k=p-1$ (bottom), fixed $h=1/5$, $\Delta t=0.01$, $\gamma=0.5$.
\label{eig_S_h5_VSp_implicit}}
\end{figure}

\begin{figure} %[!t]
\vspace{4mm}
\centerline{
\begin{tabular}{ccc}
\includegraphics[scale=0.40]{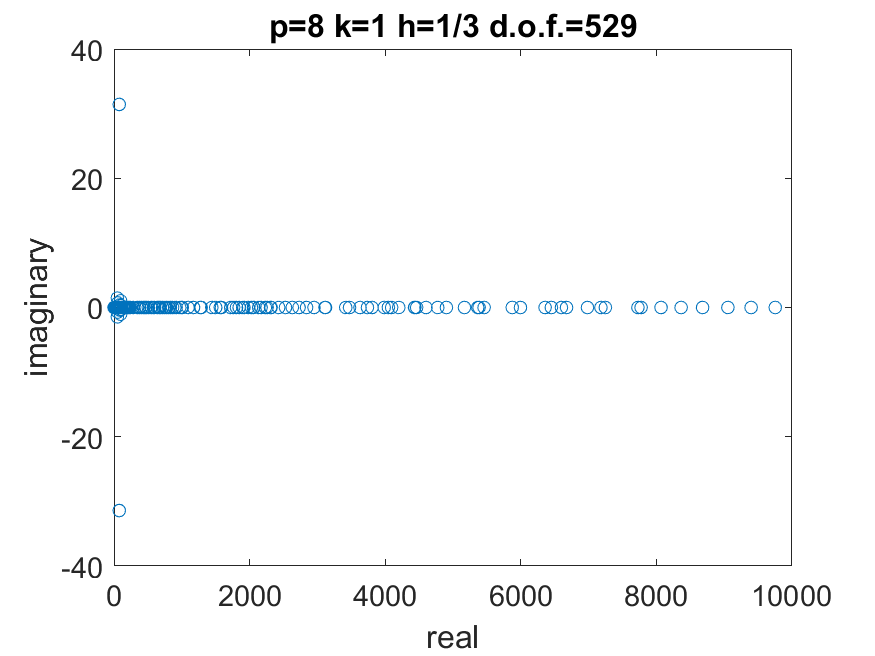} &
\hspace{-5mm}\includegraphics[scale=0.40]{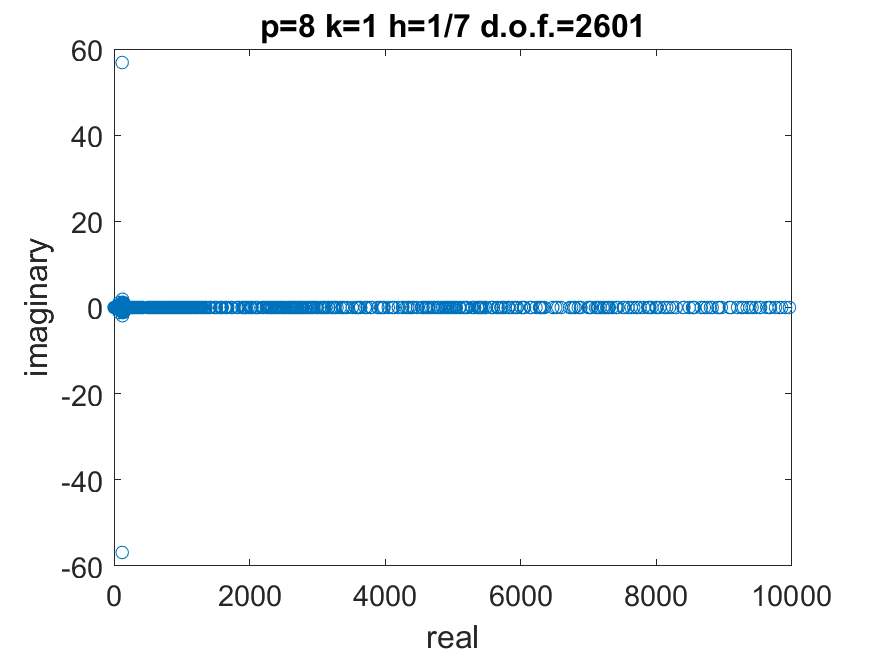} &
\hspace{-5mm}\includegraphics[scale=0.40]{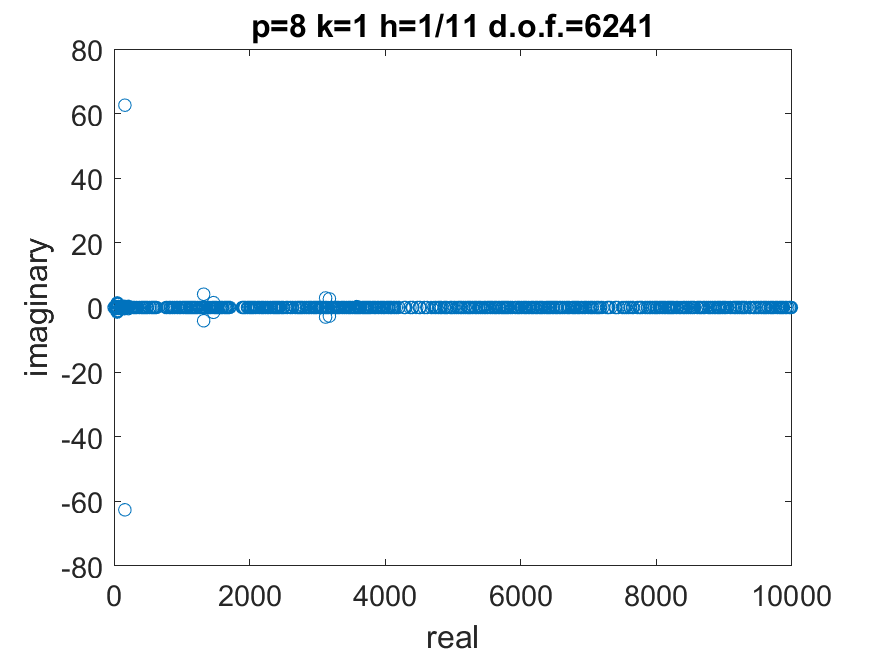} \\
\includegraphics[scale=0.40]{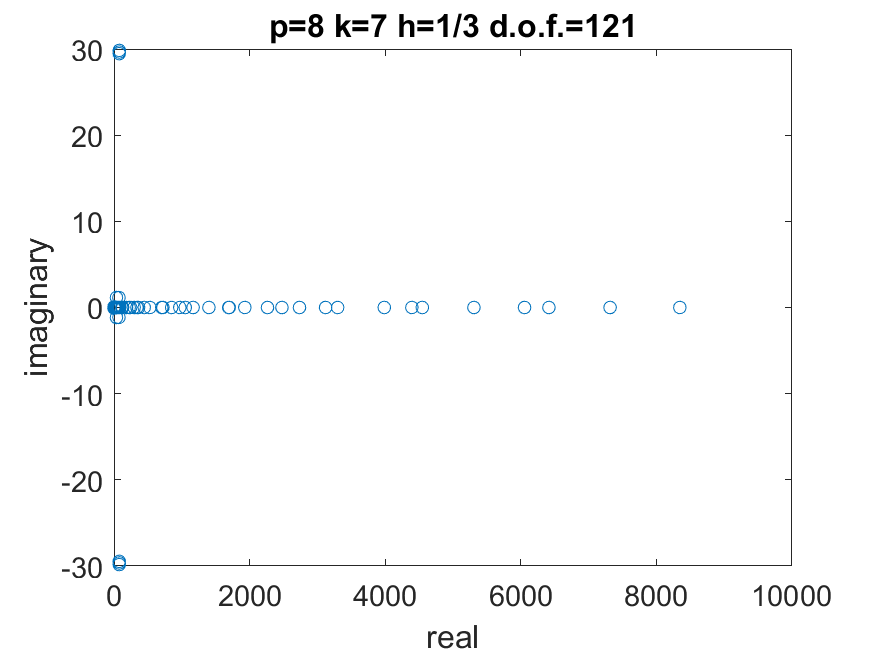} & 
\hspace{-5mm}\includegraphics[scale=0.40]{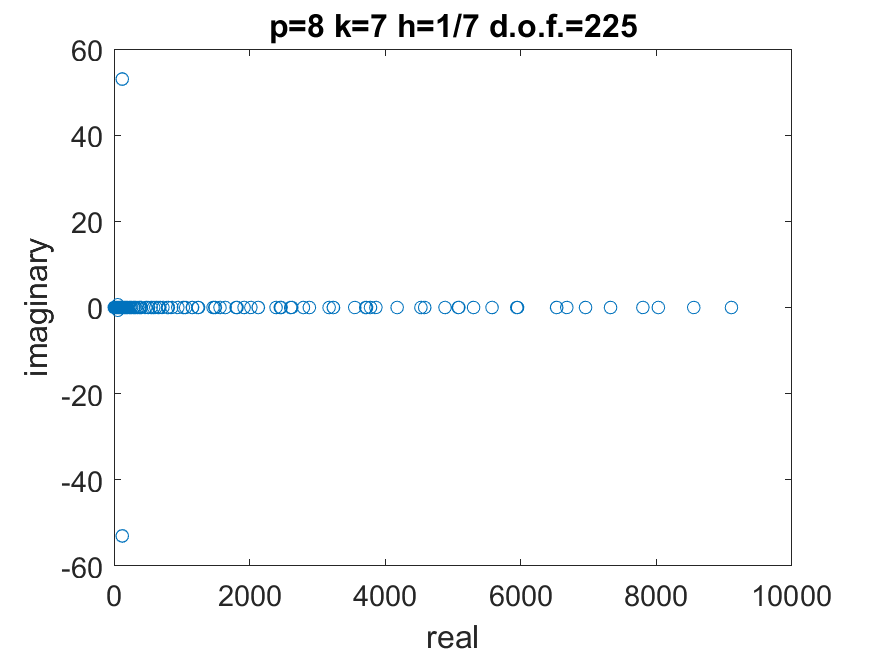} & 
\hspace{-5mm}\includegraphics[scale=0.40]{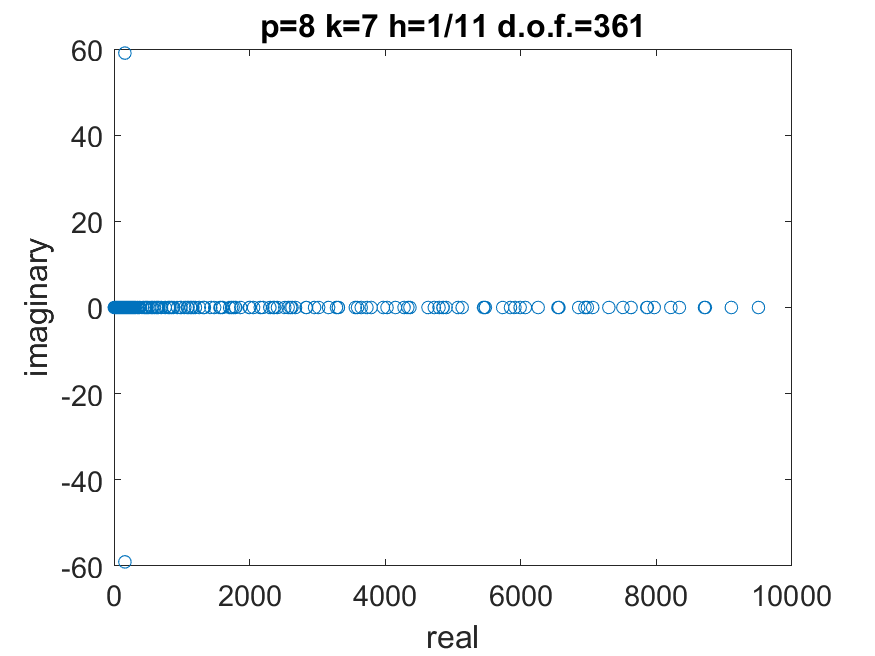}  \\
\end{tabular}
}
\vspace{-2mm}
\caption{Stiffness matrix eigenvalue distribution in the complex plane for the acoustic wave problem with  absorbing boundary conditions and implicit Newmark scheme ($\beta=0.5)$: $h=1/3$ (left),
 $1/h=7$ (center), $1/h=11$ (right), with $k=1$ (top) or $k=p-1$ (bottom), fixed $p=8$, $\Delta t=0.01$, $\gamma=0.5$.
\label{eig_S_p8_VSh_implicit}}
\end{figure}

\begin{figure} %[!t]
\vspace{4mm}
\centerline{
\begin{tabular}{ccc}
\includegraphics[scale=0.40]{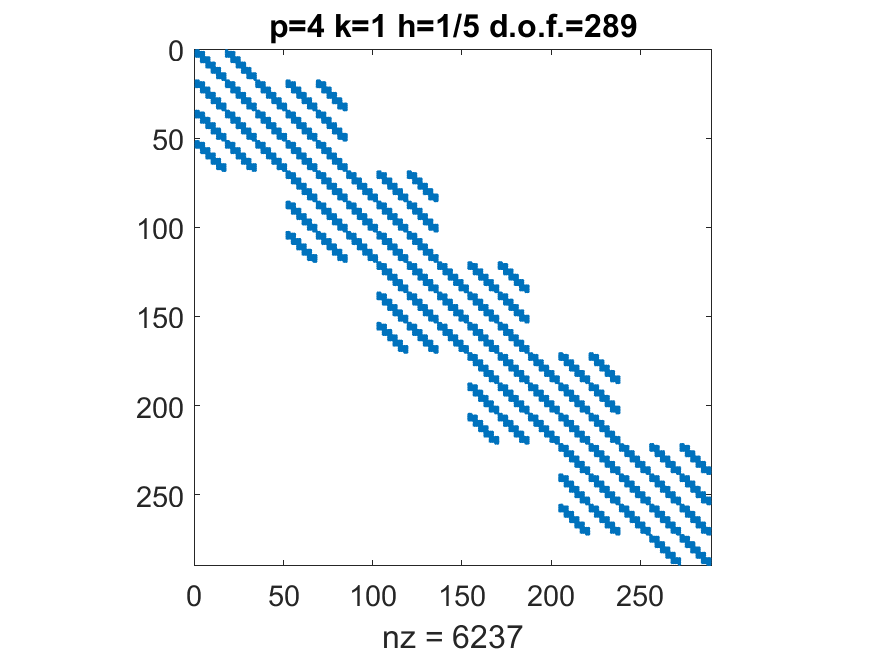} &
\hspace{-5mm}\includegraphics[scale=0.40]{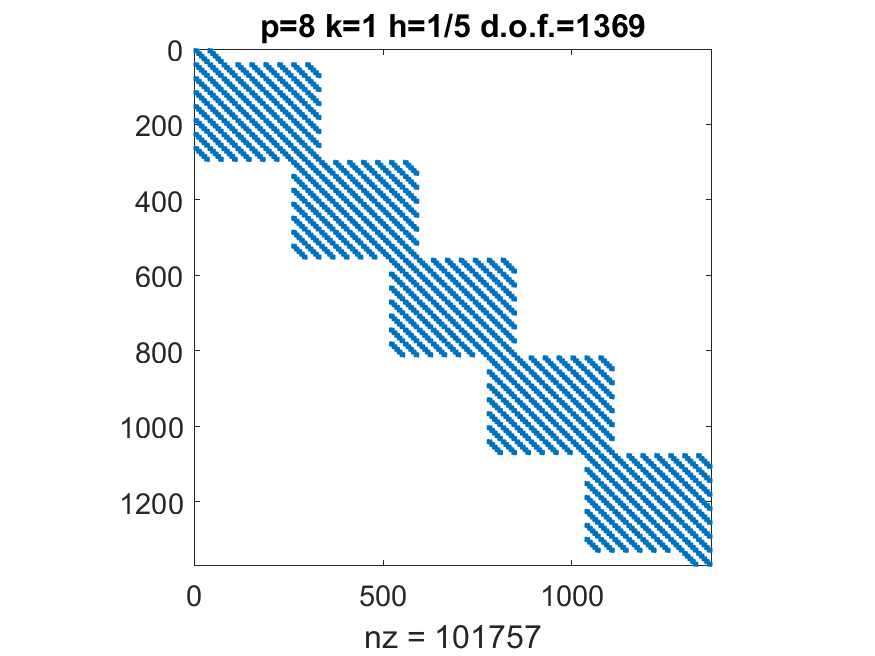} &
\hspace{-5mm}\includegraphics[scale=0.40]{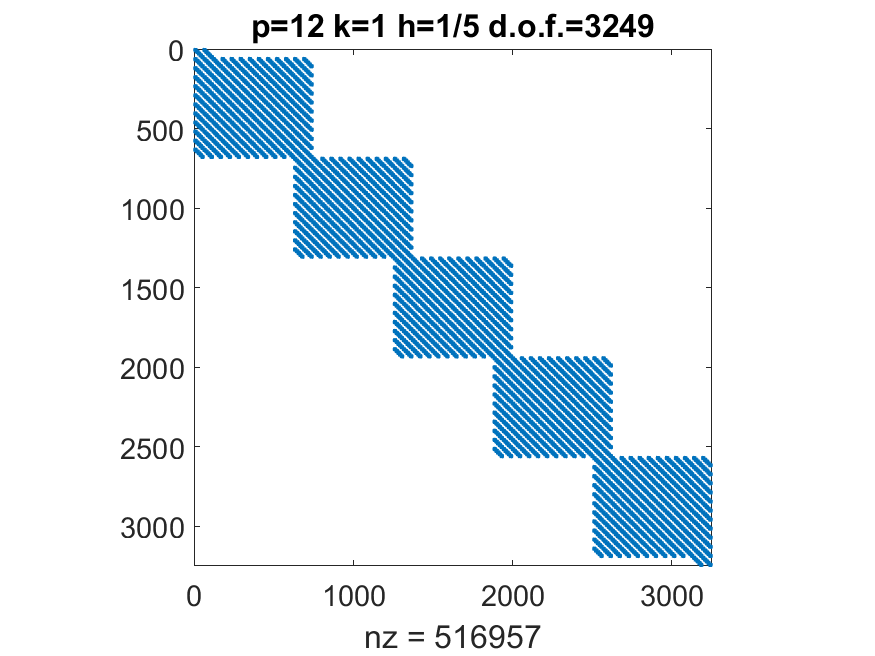} \\
\includegraphics[scale=0.40]{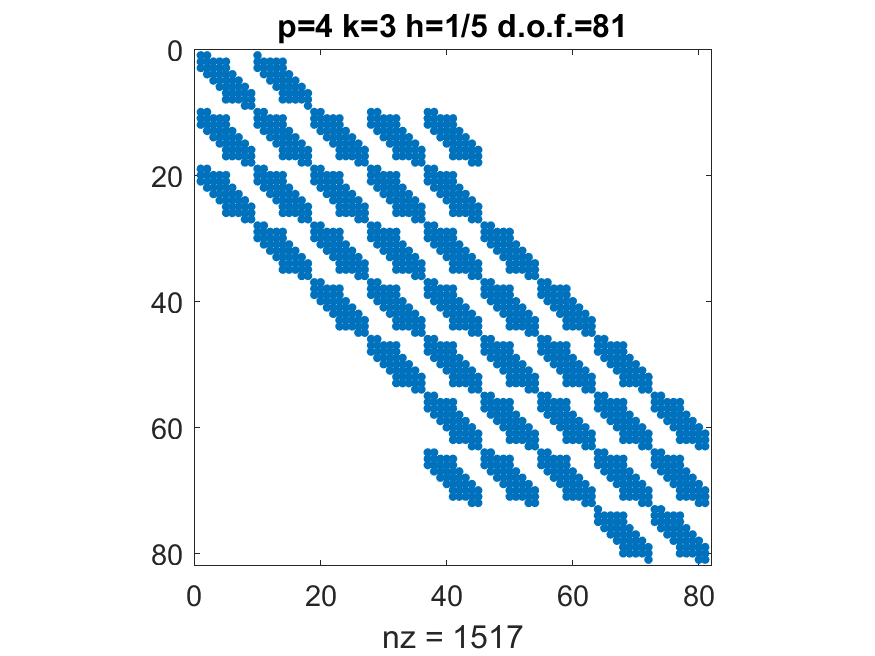} & 
\hspace{-5mm}\includegraphics[scale=0.40]{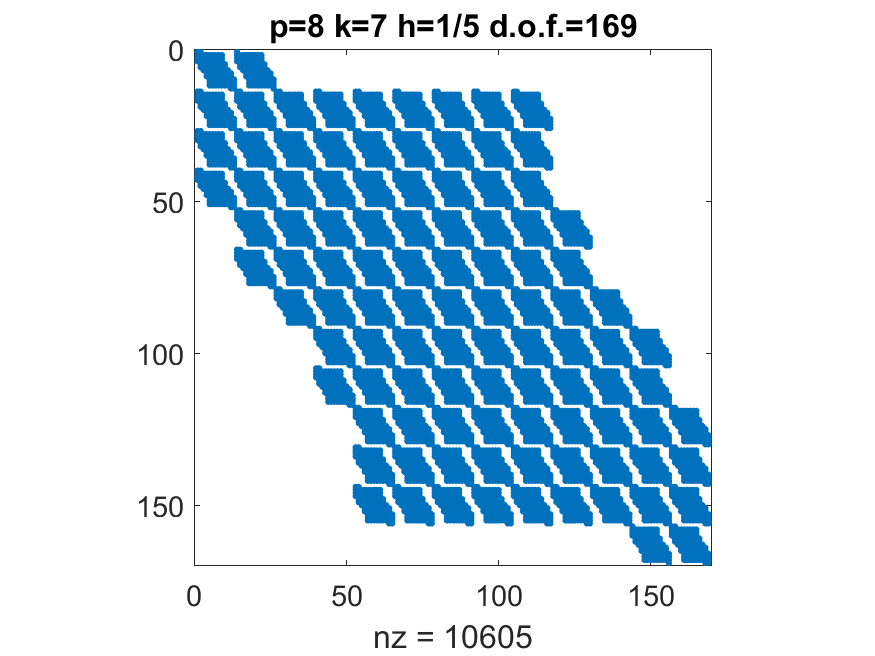} & 
\hspace{-5mm}\includegraphics[scale=0.40]{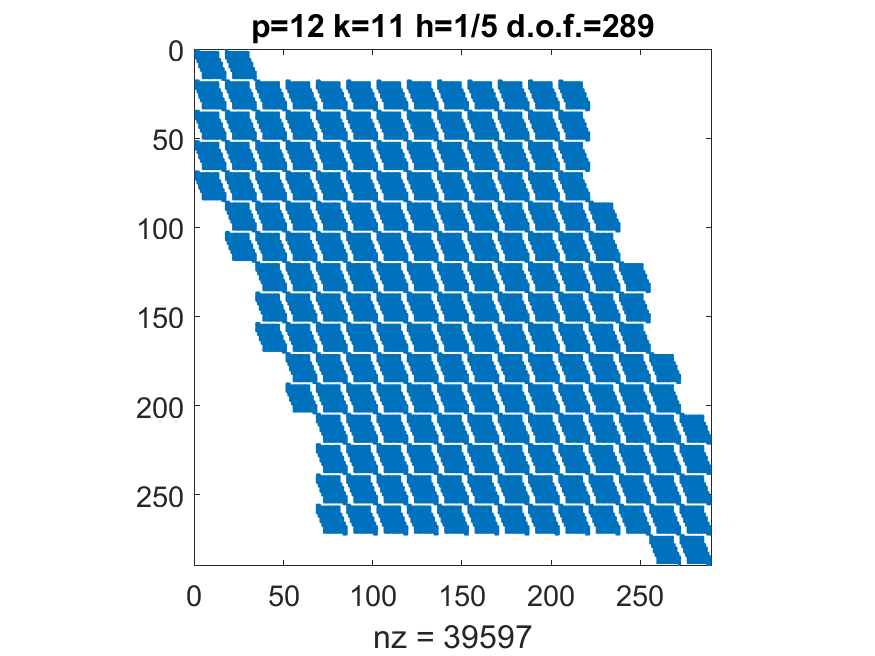}  \\
\end{tabular}
}
\vspace{-2mm}
\caption{Stiffness matrix sparsity pattern for the acoustic wave problem with  absorbing boundary conditions and implicit Newmark scheme ($\beta=0.5)$: $p = 4$ (left),
 $p=8$ (center), $p=12$ (right), with $k=1$ (top) or $k=p-1$ (bottom), fixed $h=1/5$, $\Delta t=0.01$, $\gamma=0.5$.
\label{spy_S_h5_VSp_implicit}}
\end{figure}

\begin{figure} %[!t]
\vspace{4mm}
\centerline{
\begin{tabular}{ccc}
\includegraphics[scale=0.40]{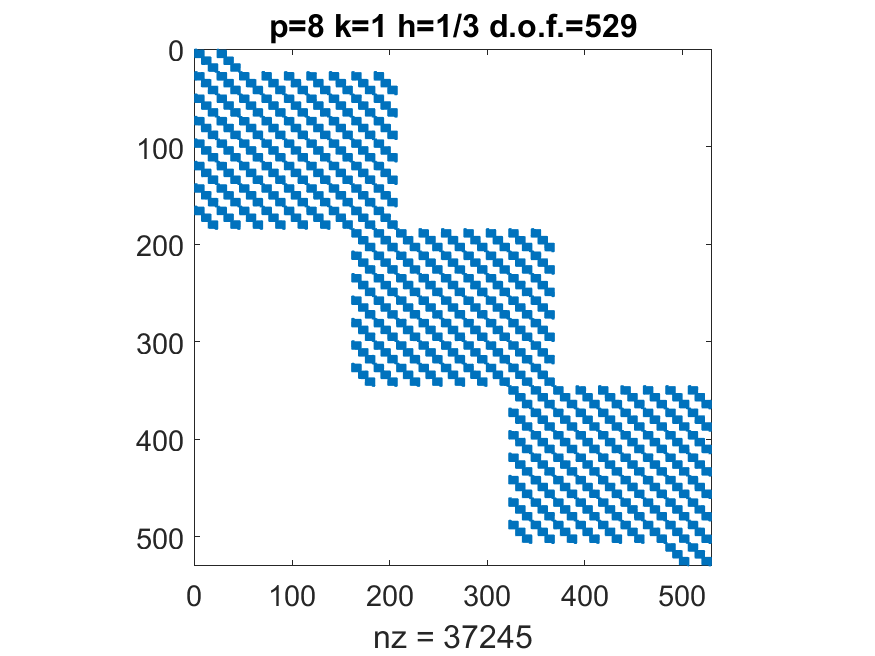} &
\hspace{-5mm}\includegraphics[scale=0.40]{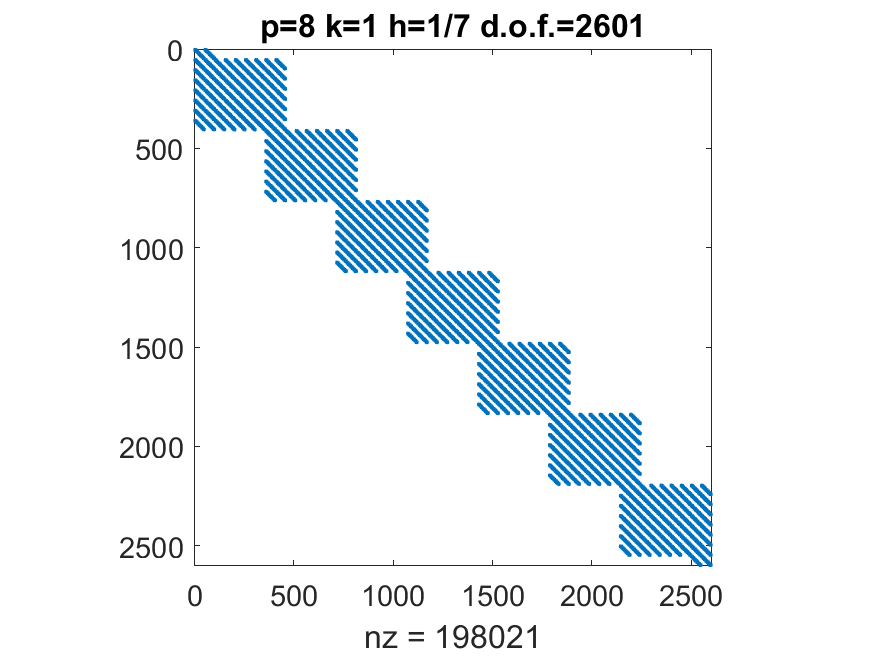} &
\hspace{-5mm}\includegraphics[scale=0.40]{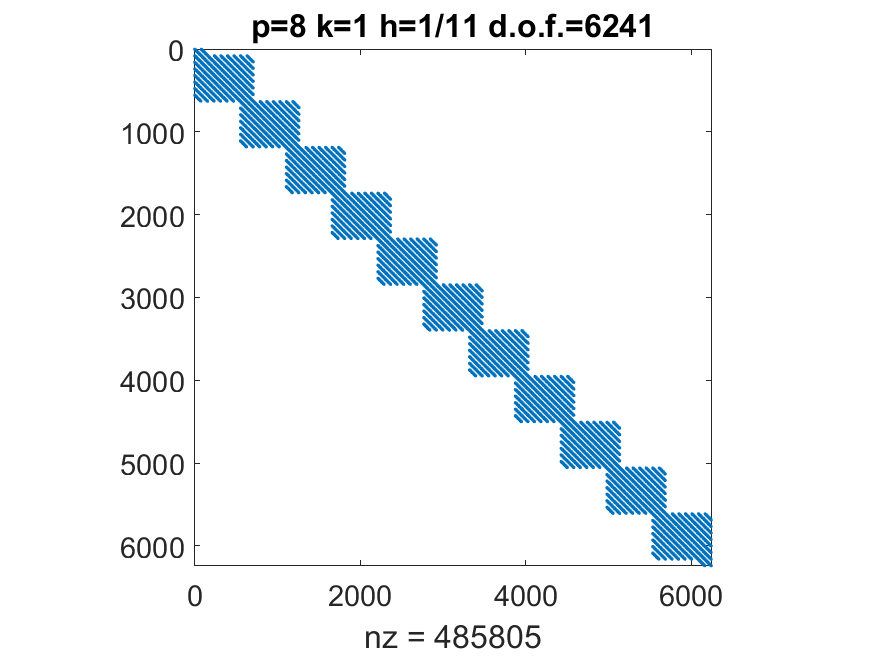} \\
\includegraphics[scale=0.40]{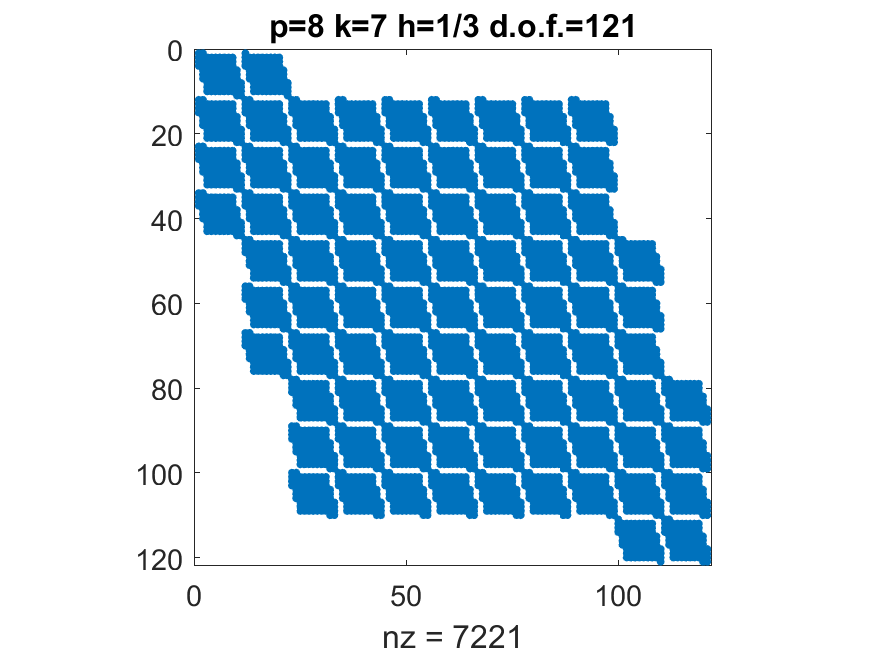} & 
\hspace{-5mm}\includegraphics[scale=0.40]{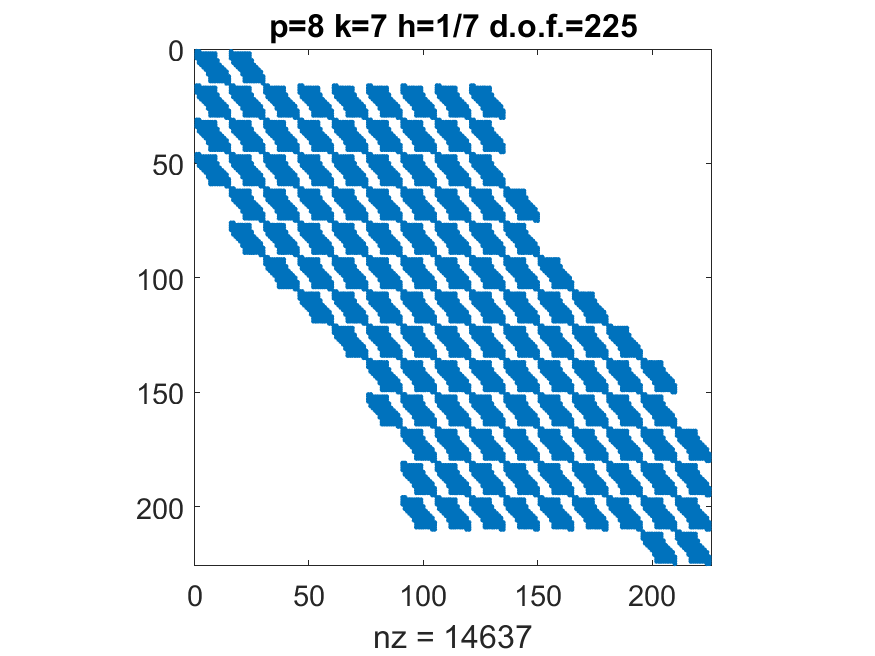} & 
\hspace{-5mm}\includegraphics[scale=0.40]{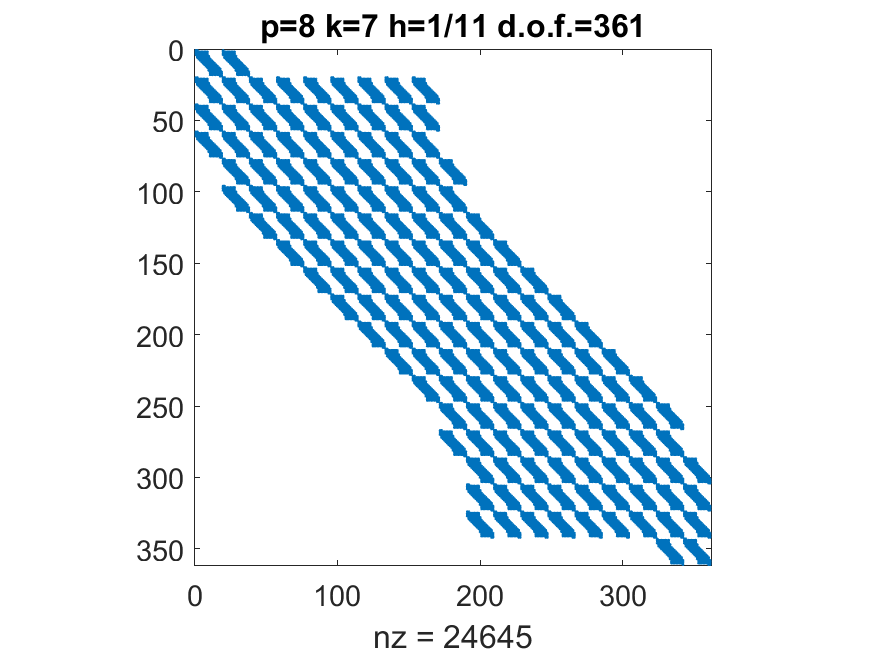}  \\
\end{tabular}
}
\vspace{-2mm}
\caption{Stiffness matrix sparsity pattern for the acoustic wave problem with  absorbing boundary conditions and implicit Newmark scheme ($\beta=0.5)$: $h=1/3$ (left),
 $1/h=7$ (center), $1/h=11$ (right), with $k=1$ (top) or $k=p-1$ (bottom), fixed $p=8$, $\Delta t=0.01$, $\gamma=0.5$.
\label{spy_S_p8_VSh_implicit}}
\end{figure}

 Figs. \ref{spy_S_h5_VSp_explicit} and  \ref{spy_S_p8_VSh_explicit}  report the sparsity pattern of the stiffness matrix $\cal K$ for the the acoustic wave problem with absorbing boundary conditions, explicit Newmark scheme ($\beta=0)$, for the same parameter settings as in Figs.  \ref{eig_S_h5_VSp_Laplacian_DIR} and \ref{eig_S_p8_VSh_Laplacian_DIR}, respectively. The results are similar, with block-diagonal matrices for minimal regularity and almost full matrices in the case of maximal regularity. Again, both {\textsf { d.o.f.}} and {\textsf  {nz}} decrease for increasing regularity. 
In Figs. \ref{eig_S_h5_VSp_implicit},  \ref{eig_S_p8_VSh_implicit},  \ref{spy_S_h5_VSp_implicit}, and  \ref{spy_S_p8_VSh_implicit}, we report the same data as in Figs.  \ref{eig_S_h5_VSp_explicit},  \ref{eig_S_p8_VSh_explicit},  \ref{spy_S_h5_VSp_explicit}, and  \ref{spy_S_p8_VSh_explicit}, respectively, for    implicit Newmark scheme ($\beta=0.5)$, all other values of parameters being unchanged. The sparsity pattern   are similar to the explicit case, but we observe definitely fewer complex eigenvalues. 
%and in some cases they reduce to only one couple.

%In Fig. \ref{eig_S_VSbc} we compare the eigenvalue distribution in complex plane of stiffness matrix %$\cal K$ for the acoustic wave problem with  different types of boundary conditions: Dirichlet (left), %Neumann (center), absorbing (right), and minimal regularity $k=1$ (top) or maximal regularity $k=p-1$ %(bottom). We choose  $\beta=0.5$ (implicit Newmark), $\gamma=0.5$, degree $p=8$, mesh size $h=1/7$, %$\Delta t=0.01$. ( QUI NON SO COSA DIRE)

%In Fig.\ref{spy_S_VSbc}  we show the pattern of stiffness matrix $\cal K$  for the same problem and same %values of parameters as in Fig.  \ref{eig_S_VSbc}. There are not significant differences except for the %first and last row-block.

\section{Conclusions}

In this paper we have considered the spectral properties of the mass and stiffness matrices related to the approximation of the acoustic wave equation with Dirichlet, Neumann and absorbing boundary conditions in the reference square domain by IGA collocation   methods in space and  Newmark advancing schemes in time, both explicit and implicit. Since no theoretical results are yet available in the literature for IGA collocation matrices, we have  addressed a systematic numerical study of the eigenvalue distribution, condition numbers and sparsity of the mass and stiffness matrices varying the discretization parameters   $p$,  $h$,   $k$,   $\Delta t$, $\beta$, {\textsf { d.o.f.}} and {\textsf  {nz}}. This analysis is of interest not only 
in order to estimate the maximum allowable  time step $\Delta t$ for explicit Newmark schemes, but also for the possible investigation of efficient preconditioned iterative solutions of the linear systems arising at each step of the time advancing schemes, since the corresponding matrices  are full and non-symmetric  for both explicit and implicit methods.

Despite the lack of specific collocation bounds and estimates, we have compared some results of our tests with the theory and numerical experiments available for matrices resulting from IGA Galerkin approximation of the Laplacian with Dirichlet boundary conditions.  
Our results show that analogous estimates for the condition numbers of mass and stiffness matrices hold also for IGA collocation discretizations of acoustic wave problems, and in some case the collocation bounds are better than the Galerkin ones. 

{\em Limitations and future work.}
This study was confined to acoustic wave problems in the reference square, but we do not expect different trends and technical complexity  in extending the tests to three dimensional domains, by using the tensor product structure of IGA collocation.\\
Future work will consider the issue of preconditioning the linear systems arising at each time step, as well as the extension of this work to
elastic wave problems.

% UNCOMMENT TO CHECK FOR UNUSED REFERENCES
\nocite{*}

%\centerline{\large\em\bf References}
%\newpage

\end{document}